\documentclass[12pt,reqno]{amsart}%
\usepackage{amsmath,amsthm,amscd,amsthm,upref,indentfirst}
\usepackage{amsfonts,mathrsfs}
\usepackage{amssymb,amsbsy,bm}
\usepackage{graphicx}
\usepackage{amsmath}
\usepackage{amsfonts}
\usepackage{amssymb}%
\usepackage[all]{xy}
\xyoption{tips}
\SelectTips {xy} {12}
\usepackage[UglyObsolete,small,nohug,heads=LaTeX]{diagrams}
\diagramstyle[labelstyle=\scriptstyle]

\usepackage[usenames]{color}

\theoremstyle{plain}
\newtheorem{theorem}{Theorem}
\newtheorem{assertion}[theorem]{Assertion}

\newtheorem{proposition}[theorem]{Proposition}

\theoremstyle{definition}
\newtheorem{definition}[theorem]{Definition}

\theoremstyle{remark}

\newtheorem{example}[theorem]{Example}
\numberwithin{equation}{section}
\numberwithin{theorem}{section}
\renewcommand{\mathfrak}{\textsc}
\renewcommand{\mathbf}{\bm}
\normalfont\upshape
\numberwithin{equation}{section}
\numberwithin{theorem}{section}

\usepackage[normalem]{ulem}
\usepackage[square,comma]{natbib}

\author{Steven Duplij}
\thanks{\emph{Published}: Journal of Kharkov National University,
(ser. Nuclei, Particles and Fields), Vol. 1017, ╣ 3(55) (2012) 28-59.}
\address{Center for Mathematics, Science and Education Rutgers University, 118
Frelinghuysen Rd., Piscataway, NJ 08854-8019}
\curraddr{V.N. Karazin
Kharkov National University, Svoboda Sq. 4, Kharkov 61022, Ukraine}
\email{duplij@math.rutgers.edu, sduplij@gmail.com}
\urladdr{http://homepages.spa.umn.edu/\~{}duplij}
\title[Polyadic systems, representations and quantum groups]{\textbf{POLYADIC SYSTEMS, REPRESENTATIONS AND QUANTUM GROUPS} }
\date{20 May 2012}

\begin{document}
\maketitle

\begin{abstract}
\noindent Polyadic systems and their representations are reviewed and a classification of general polyadic systems is presented.  A new multiplace generalization of associativity preserving homomorphisms, a 'heteromorphism' which connects polyadic systems having unequal arities, is introduced via an explicit formula, together with related definitions for multiplace representations and multiactions. Concrete examples of matrix representations for some ternary groups are then reviewed. Ternary algebras and Hopf algebras are defined, and their properties are studied. At the end some ternary generalizations of quantum groups and the Yang-Baxter equation are presented.
\end{abstract}

\bigskip

\vskip 1cm

\textsc{\tableofcontents}

\newpage

\section{Introduction}

One of the most promising directions in generalizing physical theories is the
consideration of higher arity algebras \cite{ker1}, in other words ternary and
$n$-ary algebras, in which the binary composition law is substituted by a
ternary or $n$-ary one \cite{azc/izq}.

Firstly, ternary algebraic operations (with the arity $n=3$) were introduced
already in the XIX-th century by A. Cayley in 1845 and later by J. J.
Sylvester in 1883. The notion of an $n$-ary group was introduced in 1928 by
\cite{dor3} (inspired by E. N\"{o}ther) and is a natural generalization of the
notion of a group. Even before this, in 1924, a particular case, that is, the
ternary group of idempotents, was used in \cite{pru} to study infinite abelian
groups. The important coset theorem of Post explained the connection between
$n$-ary groups and their covering binary groups \cite{pos}. The next step in
study of $n$-ary groups was the Gluskin-Hossz\'{u} theorem \cite{hos,glu1}.
Another definition of $n$-ary groups can be given as a universal algebra with
additional laws \cite{dud/gla/gle} or identities containing special elements
\cite{rus}.

The representation theory of (binary) groups \cite{weyl,FultonHarris} plays an
important role in their physical applications \cite{cornwell}. It is initially
based on a matrix realization of the group elements with the abstract group
action realized as the usual matrix multiplication \cite{cur/rei,collins}. The
cubic and $n$-ary generalizations of matrices and determinants were made in
\cite{kap/gel/zel,sokolov}, and their physical application appeared in
\cite{kaw2,rau2}. In general, particular questions of $n$-ary group
representations were considered, and matrix representations derived, by the
author \cite{bor/dud/dup3}, and some general theorems connecting
representations of binary and $n$-ary groups were presented in \cite{dud/sha}.
The intention here is to generalize the above constructions of $n$-ary group
representations to more complicated and nontrivial cases.

In physics, the most applicable structures are the nonassociative Grassmann,
Clifford and Lie algebras \cite{loh/paa/sor,lou/abl,georgi}, and so their
higher arity generalizations play the key role in further applications.
Indeed, the ternary analog of Clifford algebra was considered in \cite{abr},
and the ternary analog of Grassmann algebra \cite{abr0} was exploited to
construct various ternary extensions of supersymmetry \cite{abr/ker/roy}.

The construction of realistic physical models is based on Lie algebras, such
that the fields take their values in a concrete binary Lie algebra
\cite{georgi}. In the higher arity studies, the standard Lie bracket is
replaced by a linear $n$-ary bracket, and the algebraic structure of the
corresponding model is defined by the additional characteristic identity for
this generalized bracket, corresponding to the Jacobi identity \cite{azc/izq}.
There are two possibilities to construct the generalized Jacobi identity: 1)
The Lie bracket is a derivation by itself; 2) A double Lie bracket vanishes,
when antisymmetrized with respect to its entries. The first case leads to the
so called Filippov algebras \cite{fil0} (or $n$-Lie algebra) and second case
corresponds to generalized Lie algebras \cite{mil/vin} (or higher order Lie algebras).

The infinite-dimensional version of $n$-Lie algebras are the Nambu algebras
\cite{nam0,tak3}, and their $n$-bracket is given by the Jacobian determinant
of $n$ functions, the Nambu bracket, which in fact satisfies the Filippov
identity \cite{fil0}. Recently, the ternary Filippov algebras were
successfully applied to a three-dimensional superconformal gauge theory
describing the effective worldvolume theory of coincident $M2$-branes of
$M$-theory \cite{bag/lam2,bag/lam1,gus1}. The infinite-dimensional Nambu
bracket realization \cite{ho/ima/mat/shi} gave the possibility to describe a
condensate of nearly coincident $M2$-branes \cite{low1}.

From another side, Hopf algebras \cite{abe,sweedler,montgomery} play a
fundamental role in quantum group theory \cite{kassel,shn/ste}. Previously,
their Von Neumann generalization was introduced in
\cite{dup/li3,dup/sin1,dup/li2}, their actions on the quantum plane were
classified in \cite{dup/sin3}, and ternary Hopf algebras were defined and
studied in \cite{dup26,bor/dud/dup1}.

The goal of this paper is to give a comprehensive review of polyadic systems
and their representations. First, we classify general polyadic systems and
introduce $n$-ary semigroups and groups. Then we consider their homomorphisms
and multiplace generalizations, paying attention to their associativity. We
define multiplace representations and multiactions, and give examples of
matrix representations for some ternary groups. We define and investigate
ternary algebras and Hopf algebras, study their properties and give some
examples. At the end we consider some ternary generalizations of quantum
groups and the Yang-Baxter equation.

\section{Preliminaries}

Let $G$ be a non-empty set (underlying set, universe, carrier), its elements
we denote by lower-case Latin letters $g_{i}\in G$. The $n$\textit{-tuple} (or
\textit{polyad}) $g_{1},\ldots,g_{n}$ of elements from $G$ is denoted by
$\left(  g_{1},\ldots,g_{n}\right)  $. The Cartesian product\footnote{We place
the sign for the Cartesian product $\left(  \times\right)  $ into the power,
because the same abbreviation will also be used below for other types of
product.} $\overset{n}{\overbrace{G\times\ldots\times G}}=G^{\times n}$
consists of all $n$-tuples $\left(  g_{1},\ldots,g_{n}\right)  $, such that
$g_{i}\in G$, $i=1,\ldots,n$. For all equal elements $g\in G$, we denote
$n$-tuple (polyad) by power $\left(  g^{n}\right)  $. If the number of
elements in the $n$-tuple is clear from the context or is not important, we
denote it with one bold letter $\left(  \mathbf{g}\right)  $, in other cases
we use the power in brackets $\left(  \mathbf{g}^{\left(  n\right)  }\right)
$. We now introduce two important constructions on sets.

\begin{definition}
The $i$\textit{-projection} of the Cartesian product $G^{\times n}$ on its
$i$-th \textquotedblleft axis\textquotedblright\ is the map $\mathsf{Pr}%
_{i}^{\left(  n\right)  }:G^{\times n}\rightarrow G$ such that $\left(
g_{1},\ldots g_{i},\ldots,g_{n}\right)  \longmapsto g_{i}$.
\end{definition}

\begin{definition}
The $i$\textit{-diagonal} $\mathsf{Diag}_{n}:G\rightarrow G^{\times n}$ sends
one element to the equal element $n$-tuple $g\longmapsto\left(  g^{n}\right)
$.
\end{definition}

The one-point set $\left\{  \bullet\right\}  $ can be treated as a unit for
the Cartesian product, since there are bijections between $G$ and
$G\times\left\{  \bullet\right\}  ^{\times n}$, where $G$ can be on any place.
On the Cartesian product $G^{\times n}$ one can define a polyadic ($n$-ary,
$n$-adic, if it is necessary to specify $n$, its arity or rank) operation
$\mu_{n}:G^{\times n}\rightarrow G$. For operations we use small Greek letters
and place arguments in square brackets $\mu_{n}\left[  \mathbf{g}\right]  $.
The operations with $n=1,2,3$ are called \textit{unary, binary and ternary}.
The case $n=0$ is special and corresponds to fixing a distinguished element of
$G$, a \textquotedblleft constant\textquotedblright\ $c\in G$, and it is
called a \textit{0-ary operation} $\mu_{0}^{\left(  c\right)  }$, which maps
the one-point set $\left\{  \bullet\right\}  $ to $G$, such that $\mu
_{0}^{\left(  c\right)  }:\left\{  \bullet\right\}  \rightarrow G$, and
formally has the value $\mu_{0}^{\left(  c\right)  }\left[  \left\{
\bullet\right\}  \right]  =c\in G$. The 0-ary operation \textquotedblleft
kills\textquotedblright\ arity, which can be seen from the following
\cite{bergman1}: the composition of $n$-ary and $m$-ary operations $\mu
_{n}\circ\mu_{m}$ gives $\left(  n+m-1\right)  $-ary operation by%
\begin{equation}
\mu_{n+m-1}\left[  \mathbf{g},\mathbf{h}\right]  =\mu_{n}\left[
\mathbf{g},\mu_{m}\left[  \mathbf{h}\right]  \right]  . \label{m}%
\end{equation}
Then, if to compose $\mu_{n}$ with the 0-ary operation $\mu_{0}^{\left(
c\right)  }$, we obtain%
\begin{equation}
\mu_{n-1}^{\left(  c\right)  }\left[  \mathbf{g}\right]  =\mu_{n}\left[
\mathbf{g},c\right]  , \label{mc}%
\end{equation}
because $\mathbf{g}$ is a polyad of length $\left(  n-1\right)  $. So, it is
necessary to make a clear distinction between the 0-ary operation $\mu
_{0}^{\left(  c\right)  }$ and its value $c$ in $G$, as will be seen and will
become important below.

\begin{definition}
A \textit{polyadic system} $\mathfrak{G}$ is a set $G$ which is closed under
polyadic operations.
\end{definition}

We will write $\mathfrak{G}=\left\langle \text{set}|\text{operations}%
\right\rangle $ or $\mathfrak{G}=\left\langle \text{set}|\text{operations}%
|\text{relations}\right\rangle $, where \textquotedblleft
relations\textquotedblright\ are some additional properties of operations
(e.g., associativity conditions for semigroups or cancellation properties). In
such a definition it is not necessary to list the images of 0-ary operations
(e.g. the unit or zero in groups), as is done in various other definitions.
Here, we mostly consider concrete polyadic systems with one \textquotedblleft
chief\textquotedblright\ (fundamental) $n$-ary operation $\mu_{n}$, which is
called \textit{polyadic multiplication} (or $n$\textit{-ary multiplication}).

\begin{definition}
A $n$\textit{-ary system} $\mathfrak{G}_{n}=\left\langle G\mid\mu
_{n}\right\rangle $ is a set $G$ closed under one $n$-ary operation $\mu_{n}$
(without any other additional structure).
\end{definition}

Note that a set with one closed binary operation without any other relations
was called a groupoid by Hausmann and Ore \cite{hau/ore} (see, also
\cite{cli/pre1}). However, nowadays the term \textquotedblleft
groupoid\textquotedblright\ is widely used in category theory and homotopy
theory for a different construction with binary multiplication, the so-called
Brandt groupoid \cite{bra1} (see, also, \cite{bruck}). Alternatively, and much
later on, Bourbaki \cite{bourbaki1} introduced the term \textquotedblleft
magma\textquotedblright\ for binary systems. Then, the above terms were
extended to the case of one fundamental $n$-ary operation as well.
Nevertheless, we will use some neutral notations \textquotedblleft polyadic
system\textquotedblright\ and \textquotedblleft$n$-ary
system\textquotedblright\ (when arity $n$ is fixed/known/important), which
adequately indicates all of their main properties.

Let us consider the \textit{changing arity problem}:

\begin{definition}
\label{def-charity}For a given $n$-ary system $\left\langle G\mid\mu
_{n}\right\rangle $ to construct another polyadic system $\left\langle
G\mid\mu_{n^{\prime}}^{\prime}\right\rangle $ over the same set $G$, which has
multiplication with a different arity $n^{\prime}$.
\end{definition}

The formulas (\ref{m}) and (\ref{mc}) give us the simplest examples of how to
change the arity of a polyadic system. In general, there are 3 ways:

\begin{enumerate}
\item \textsl{Iterating}. Using composition of the operation $\mu_{n}$ with
itself, one can increase the arity from $n$ to $n_{iter}^{\prime}$ (as in
(\ref{m})) without changing the signature of the system. We denote the number
of iterating multiplications by $\ell_{\mu}$, and use the bold Greek letters
$\mathbf{\mu}_{n}^{\ell_{\mu}}$ for the resulting composition of $n$-ary
multiplications, such that%
\begin{equation}
\mu_{n^{\prime}}^{\prime}=\mathbf{\mu}_{n}^{\ell_{\mu}}\overset{def}%
{=}\overset{\ell_{\mu}}{\overbrace{\mu_{n}\circ\left(  \mu_{n}\circ
\ldots\left(  \mu_{n}\times\operatorname*{id}\nolimits^{\times\left(
n-1\right)  }\right)  \ldots\times\operatorname*{id}\nolimits^{\times\left(
n-1\right)  }\right)  }}, \label{mn}%
\end{equation}
where%
\begin{equation}
n^{\prime}=n_{iter}=\ell_{\mu}\left(  n-1\right)  +1, \label{n}%
\end{equation}
which gives the length of a polyad $\left(  \mathbf{g}\right)  $ in the
notation $\mathbf{\mu}_{n}^{\ell_{\mu}}\left[  \mathbf{g}\right]  $. Without
assuming associativity there many variants for placing $\mu_{n}$'s among
$\operatorname*{id}${}'s in the r.h.s. of (\ref{mn}). The operation
$\mathbf{\mu}_{n}^{\ell_{\mu}}$ is named a \textit{long product }\cite{dor3}
or \textit{derived} \cite{dud07}.

\item \textsl{Reducing (Collapsing)}. Using $n_{c}$ distinguished elements or
constants (or $n_{c}$ additional $0$-ary operations $\mu_{0}^{\left(
c_{i}\right)  }$, $i=1,\ldots n_{c}$), one can decrease arity from $n$ to
$n_{red}^{\prime}$ (as in (\ref{mc})), such that\footnote{In \cite{dud/mic2}
$\mu_{n}^{\left(  c_{1}\ldots c_{n_{c}}\right)  }$ is named a retract (which
term is already busy and widely used in category theory for another
construction).}%
\begin{equation}
\mu_{n^{\prime}}^{\prime}=\mu_{n^{\prime}}^{\left(  c_{1}\ldots c_{n_{c}%
}\right)  }\overset{def}{=}\mu_{n}\circ\left(  \overset{n_{c}}{\overbrace
{\mu_{0}^{\left(  c_{1}\right)  }\times\ldots\times\mu_{0}^{\left(  c_{n_{c}%
}\right)  }}}\times\operatorname*{id}\nolimits^{\times\left(  n-n_{c}\right)
}\right)  ,
\end{equation}
where%
\begin{equation}
n^{\prime}=n_{red}=n-n_{c}, \label{nr}%
\end{equation}
and the $0$-ary operations $\mu_{0}^{\left(  c_{i}\right)  }$ can be on any places.

\item \textsl{Mixing}. Changing (increasing or decreasing) arity may be done
by combining iterating and reducing (maybe with additional operations of
different arity). If we do not use additional operations, the final arity can
be presented in a general form using (\ref{n}) and (\ref{nr}). It will depend
on the order of iterating and reducing, and so we have two subcases:

\begin{enumerate}
\item \textsl{Iterating}$\rightarrow$\textsl{Reducing}. We have%
\begin{equation}
n^{\prime}=n_{iter\rightarrow red}=\ell_{\mu}\left(  n-1\right)  -n_{c}+1.
\end{equation}
The maximal number of constants (when $n_{iter\rightarrow red}^{\prime}=2$) is
equal to%
\begin{equation}
n_{c}^{\max}=\ell_{\mu}\left(  n-1\right)  -1
\end{equation}
and can be increased by increasing the number of multiplications $\ell_{\mu}$.

\item \textsl{Reducing}$\rightarrow$\textsl{Iterating}. We obtain%
\begin{equation}
n^{\prime}=n_{red\rightarrow iter}=\ell_{\mu}\left(  n-1-n_{c}\right)  +1.
\end{equation}

Now the maximal number of constants is%
\begin{equation}
n_{c}^{\max}=n-2
\end{equation}
and this is achieved only when $\ell_{\mu}=1$.
\end{enumerate}
\end{enumerate}

To give examples of the third (mixed) case we put $n=4$, $\ell_{\mu}=3$,
$n_{c}=2$ for both subcases of opposite ordering:

\begin{enumerate}
\item \textsl{Iterating}$\rightarrow$\textsl{Reducing}. We can put%
\begin{equation}
\mu_{8}^{\left(  c_{1},c_{2}\right)  \prime}\left[  \mathbf{g}^{\left(
8\right)  }\right]  =\mu_{4}\left[  g_{1},g_{2},g_{3},\mu_{4}\left[
g_{4},g_{5},g_{6},\mu_{4}\left[  g_{7},g_{8},c_{1},c_{2}\right]  \right]
\right]  ,
\end{equation}
which corresponds to the following commutative diagram%
\begin{equation}
\begin{diagram} G^{\times 8}& \rTo^{\epsilon}&G^{\times 8}\times\left\{ \bullet \right\}^2 & \rTo^{\operatorname*{id}^{\times 6}\times\mu_{4} } & G^{\times 7}&\rTo^{\operatorname*{id}^{\times 3}\times\mu_{4}}& G^{\times 4}\\&\;\;\;\rdTo(9,3)_{\mu_{8}^{\left( c_{1},c_{2}\right) \prime}} &&&&& \dTo_{\mu_{4}}\\ &&&&& & \\ &&& &&&G \\ \end{diagram}
\end{equation}

\item \textsl{Reducing}$\rightarrow$\textsl{Iterating}. We can have%
\begin{equation}
\mu_{4}^{\left(  c_{1},c_{2}\right)  \prime}\left[  \mathbf{g}^{\left(
4\right)  }\right]  =\mu_{4}\left[  g_{1},c_{1},c_{2},\mu_{4}\left[
g_{2},c_{1},c_{2},\mu_{4}\left[  g_{3},c_{1},c_{2},g_{4}\right]  \right]
\right]  ,
\end{equation}

such that the diagram
\begin{equation}
\begin{diagram} G^{\times 4}& \rTo^{\epsilon}&\left(G\times\left\{ \bullet \right\}^2\right)^{\times 3}\times G & \rTo^{\operatorname*{id}^{\times 6}\times\mu_{4} } & G^{\times 7}&\rTo^{\operatorname*{id}^{\times 3}\times\mu_{4}}& G^{\times 4}\\&\;\;\rdTo(11,3)_{\mu_{4}^{\left( c_{1},c_{2}\right) \prime}} &&&&& \;\;\;\;\dTo_{\mu_{4}}\\ &&&&& & \\ &&& &&&\,G \\ \end{diagram}
\end{equation}
is commutative.
\end{enumerate}

It is important to find conditions where iterating and reducing compensate
each other, i.e. they do not change arity overall. Indeed, let the number of
the iterating multiplications $\ell_{\mu}$ be fixed, then we can find such a
number of reducing constants $n_{c}^{\left(  0\right)  }$, such that the final
arity will coincide with the initial arity $n$. The result will depend on the
order of operations. There are two cases:

\begin{enumerate}
\item \textsl{Iterating}$\rightarrow$\textsl{Reducing}. For the number of
reducing constants $n_{c}^{\left(  0\right)  }$ we obtain from (\ref{n}) and
(\ref{nr})%
\begin{equation}
n_{c}^{\left(  0\right)  }=\left(  n-1\right)  \left(  \ell_{\mu}-1\right)  ,
\label{nc1}%
\end{equation}
such that there is no restriction on $\ell_{\mu}$.

\item \textsl{Reducing}$\rightarrow$\textsl{Iterating}. For $n_{c}^{\left(
0\right)  }$ we get%
\begin{equation}
n_{c}^{\left(  0\right)  }=\dfrac{\left(  n-1\right)  \left(  \ell_{\mu
}-1\right)  }{\ell_{\mu}}, \label{nc2}%
\end{equation}
and now $\ell_{\mu}\leq n-1$. The requirement that $n_{c}^{\left(  0\right)
}$ should be an integer gives two further possibilities%
\begin{equation}
n_{c}^{\left(  0\right)  }=\left\{
\begin{array}
[c]{l}%
\dfrac{n-1}{2},\ \ \ \ell_{\mu}=2,\\
n-2,\ \ \ \ell_{\mu}=n-1.
\end{array}
\right.
\end{equation}

\end{enumerate}

The above relations can be useful in the study of various $n$-ary
multiplication structures and their presentation in special form is needed in
concrete problems.

\section{Special elements and properties of polyadic systems}

Let us recall the definitions of some standard algebraic systems and their
special elements, which will be considered in this paper, using our notation.

\begin{definition}
A \textit{zero} of a polyadic system is a distinguished element $z$ (and the
corresponding 0-ary operation $\mu_{0}^{\left(  z\right)  }$) such that for
any $\left(  n-1\right)  $-tuple (polyad) $\mathbf{g\in}G^{\times\left(
n-1\right)  }$ we have%
\begin{equation}
\mu_{n}\left[  \mathbf{g},z\right]  =z, \label{z}%
\end{equation}
where $z$ can be on any place in the l.h.s. of (\ref{z}).
\end{definition}

There is only one zero (if its place is not fixed) which can be possible in a
polyadic system. As in the binary case, an analog of positive powers of an
element \cite{pos} should coincide with the number of multiplications
$\ell_{\mu}$ in the iterating (\ref{mn}).

\begin{definition}
A (positive) \textit{polyadic power} of an element is%
\begin{equation}
g^{\left\langle \ell_{\mu}\right\rangle }=\mathbf{\mu}_{n}^{\ell_{\mu}}\left[
g^{\ell_{\mu}\left(  n-1\right)  +1}\right]  . \label{pp}%
\end{equation}

\end{definition}

\begin{definition}
An element of a polyadic system $g$ is called $\ell_{\mu}$-\textit{nilpotent}
(or simply \textit{nilpotent} for $\ell_{\mu}=1$), if there exist such
$\ell_{\mu}$ that%
\begin{equation}
g^{\left\langle \ell_{\mu}\right\rangle }=z. \label{mz}%
\end{equation}

\end{definition}

\begin{definition}
\label{def-nil}A polyadic system with zero $z$ is called $\ell_{\mu}%
$-\textit{nilpotent}, if there exists $\ell_{\mu}$ such that for any $\left(
\ell_{\mu}\left(  n-1\right)  +1\right)  $-tuple (polyad) $\mathbf{g}$ we have%
\begin{equation}
\mathbf{\mu}_{n}^{\ell_{\mu}}\left[  \mathbf{g}\right]  =z. \label{mg}%
\end{equation}

\end{definition}

Therefore, the \textit{index of nilpotency} (number of elements whose product
is zero) of an $\ell_{\mu}$-nilpotent $n$-ary system is $\left(  \ell_{\mu
}\left(  n-1\right)  +1\right)  $, while its polyadic power is $\ell_{\mu}$ .

\begin{definition}
A \textit{polyadic (}$n$-\textit{ary) identity} (or neutral element) of a
polyadic system is a distinguished element $e$ (and the corresponding 0-ary
operation $\mu_{0}^{\left(  e\right)  }$) such that for any element $g\in G$
we have%
\begin{equation}
\mu_{n}\left[  g,e^{n-1}\right]  =g, \label{e}%
\end{equation}
where $g$ can be on any place in the l.h.s. of (\ref{e}).
\end{definition}

In binary groups the identity is the only neutral element, while in polyadic
systems, there exist \textit{neutral polyads} $\mathbf{n}$ consisting of
elements of $G$ satisfying%
\begin{equation}
\mu_{n}\left[  g,\mathbf{n}\right]  =g, \label{mng}%
\end{equation}
where $g$ can be also on any place. The neutral polyads are not determined
uniquely. It follows from (\ref{e}) that the sequence of polyadic identities
$e^{n-1}$ is a neutral polyad.

\begin{definition}
\label{def-mi}An element of a polyadic system $g$ is called $\ell_{\mu}%
$-\textit{idempotent} (or simply \textit{idempotent} for $\ell_{\mu}=1$), if
there exist such $\ell_{\mu}$ that%
\begin{equation}
g^{\left\langle \ell_{\mu}\right\rangle }=g. \label{mi}%
\end{equation}

\end{definition}

Both zero and the identity are $\ell_{\mu}$-idempotents with arbitrary
$\ell_{\mu}$. We define \textit{(total) associativity} as the invariance of
the composition of two $n$-ary multiplications%
\begin{equation}
\mathbf{\mu}_{n}^{2}\left[  \mathbf{g},\mathbf{h},\mathbf{u}\right]  =\mu
_{n}\left[  \mathbf{g},\mu_{n}\left[  \mathbf{h}\right]  ,\mathbf{u}\right]
=invariant\label{ghu}%
\end{equation}
under placement of the internal multiplication in r.h.s. with a fixed order of
elements in the whole polyad of $\left(  2n-1\right)  $ elements
$\mathbf{t}^{\left(  2n-1\right)  }=\left(  \mathbf{g},\mathbf{h}%
,\mathbf{u}\right)  $. Informally, \textquotedblleft internal
brackets/multiplication can be moved on any place\textquotedblright, which
gives $n$ relations%
\begin{equation}
\mu_{n}\circ\left(  \mu_{n}\times\operatorname*{id}\nolimits^{\times\left(
n-1\right)  }\right)  =\ldots=\mu_{n}\circ\left(  \operatorname*{id}%
\nolimits^{\times\left(  n-1\right)  }\times\mu_{n}\right)  .
\end{equation}
There are many other particular kinds of associativity which were introduced
in \cite{thu49} and studied in \cite{belousov,sok1}. Here we will confine
ourselves the most general, total associativity (\ref{ghu}). In this case, the
iteration does not depend on the placement of internal multiplications in the
r.h.s. of (\ref{mn}).

\begin{definition}
A \textit{polyadic semigroup} ($n$-\textit{ary semigroup}) is a $n$-ary system
in which the operation is associative, or $\mathfrak{G}_{n}^{semigrp}%
=\left\langle G\mid\mu_{n}\mid\text{associativity}\right\rangle $.
\end{definition}

In a polyadic system with zero (\ref{z}) one can have \textit{trivial
associativity}, when all $n$ terms are (\ref{ghu}) are equal to zero, i.e.%
\begin{equation}
\mathbf{\mu}_{n}^{2}\left[  \mathbf{g}\right]  =z \label{mgz}%
\end{equation}
for any $\left(  2n-1\right)  $-tuple $\mathbf{g}$. Therefore, we state that

\begin{assertion}
Any $2$-nilpotent $n$-ary system (having index of nilpotency $\left(
2n-1\right)  $) is a polyadic semigroup.
\end{assertion}

In the case of changing arity one should use in (\ref{mgz}) not the changed
final arity $n^{\prime}$, but the \textquotedblleft real\textquotedblright%
\ arity which is $n$ for the reducing case and $\ell_{\mu}\left(  n-1\right)
+1$ for all other cases. Let us give some examples.

\begin{example}
In the mixed (interacting-reducing) case with $n=2$, $\ell_{\mu}=3$, $n_{c}%
=1$, we have a ternary system $\left\langle G\mid\mu_{3}\right\rangle $
iterated from a binary system $\left\langle G\mid\mu_{2},\mu_{0}^{\left(
c\right)  }\right\rangle $ with one distinguished element $c$ (or an
additional $0$-ary operation)\footnote{This construction is named the
$b$-derived groupoid in \cite{dud/mic2}.}%
\begin{equation}
\mu_{3}^{\left(  c\right)  }\left[  g,h,u\right]  =\left(  g\cdot\left(
h\cdot\left(  u\cdot c\right)  \right)  \right)  ,
\end{equation}
where for binary multiplication we denote $g\cdot h=\mu_{2}\left[  g,h\right]
$. Thus, if the ternary system $\left\langle G\mid\mu_{3}^{\left(  c\right)
}\right\rangle $ is nilpotent of index $7$ (see \ref{def-nil}), then it is a
ternary semigroup (because $\mu_{3}^{\left(  c\right)  }$ is trivially
associative) independently of the associativity of $\mu_{2}$ (see, e.g.
\cite{bor/dud/dup3}).
\end{example}

It is very important to find the \textit{associativity preserving} conditions
(constructions), where an associative initial operation $\mu_{n}$ leads to an
associative final operation $\mu_{n^{\prime}}^{\prime}$ during the change of arity.

\begin{example}
An associativity preserving reduction can be given by the construction of a
binary associative operation using $\left(  n-2\right)  $-tuple $\mathbf{c}$
consisting of $n_{c}=n-2$ different constants%
\begin{equation}
\mu_{2}^{\left(  \mathbf{c}\right)  }\left[  g,h\right]  =\mu_{n}\left[
g,\mathbf{c},h\right]  . \label{mgh}%
\end{equation}

\end{example}

Associativity preserving mixing constructions with different arities and
places were considered in \cite{dud/mic2,mic2,sok1}.

\begin{definition}
An associative polyadic system with identity (\ref{e}) is called a
\textit{polyadic monoid}.
\end{definition}

The structure of any polyadic monoid is fixed \cite{pop/pop}: it can be
obtained by iterating a binary operation \cite{cup/trp} (for polyadic groups
this was shown in \cite{dor3}).

In polyadic systems, there are several analogs of binary commutativity. The
most straightforward one comes from commutation of the multiplication with permutations.

\begin{definition}
A polyadic system is $\sigma$-\textit{commutative}, if $\mu_{n}=\mu_{n}%
\circ\sigma$, or%
\begin{equation}
\mu_{n}\left[  \mathbf{g}\right]  =\mu_{n}\left[  \sigma\circ\mathbf{g}%
\right]  , \label{ms}%
\end{equation}
where $\sigma\circ\mathbf{g}=\left(  g_{\sigma\left(  1\right)  }%
,\ldots,g_{\sigma\left(  n\right)  }\right)  $ is a permutated polyad and
$\sigma$ is a fixed element of $S_{n}$, the permutation group on $n$ elements.
If (\ref{ms}) holds for all $\sigma\in S_{n}$, then a polyadic system is
\textit{commutative}.
\end{definition}

A special type of the $\sigma$-commutativity%
\begin{equation}
\mu_{n}\left[  g,\mathbf{t},h\right]  =\mu_{n}\left[  h,\mathbf{t},g\right]  ,
\label{mth}%
\end{equation}
where $\mathbf{t}$ is any fixed $\left(  n-2\right)  $-polyad, is called
\textit{semicommutativity}. So for a $n$-ary semicommutative system we have%
\begin{equation}
\mu_{n}\left[  g,h^{n-1}\right]  =\mu_{n}\left[  h^{n-1},g\right]  .
\end{equation}

If a $n$-ary semigroup $\mathfrak{G}^{semigrp}$ is iterated from a commutative
binary semigroup with identity, then $\mathfrak{G}^{semigrp}$ is semicommutative.

\begin{example}
\label{exam-monoid} Let $G$ be the set of natural numbers $\mathbb{N}$, and
the 5-ary multiplication is defined by%
\begin{equation}
\mu_{5}\left[  \mathbf{g}\right]  =g_{1}-g_{2}+g_{3}-g_{4}+g_{5},
\end{equation}
then $\mathfrak{G}_{5}^{\mathbb{N}}=\left\langle \mathbb{N},\mu_{5}%
\right\rangle $ is a semicommutative 5-ary monoid having the identity
$e_{g}=\mu_{5}\left[  g^{5}\right]  =g$ for each $g\in\mathbb{N}$. Therefore,
$\mathfrak{G}_{5}^{\mathbb{N}}$ is the idempotent monoid.
\end{example}

Another possibility is to generalize the binary \textit{mediality} in
semigroups%
\begin{equation}
\left(  g_{11}\cdot g_{12}\right)  \cdot\left(  g_{21}\cdot g_{22}\right)
=\left(  g_{11}\cdot g_{21}\right)  \cdot\left(  g_{12}\cdot g_{22}\right)  ,
\end{equation}
which, obviously, follows from binary commutativity. But for $n$-ary systems
they are different. It is seen that the mediality should contain $\left(
n+1\right)  $ multiplications, it is a relation between $n\times n$ elements,
and therefore can be presented in a matrix from. The latter can be achieved by
placing the arguments of the external multiplication in a column.

\begin{definition}
A polyadic system is \textit{medial} (or \textit{entropic}), if
\cite{eva7,belousov}%
\begin{equation}
\mu_{n}\left[
\begin{array}
[c]{c}%
\mu_{n}\left[  g_{11},\ldots,g_{1n}\right] \\
\vdots\\
\mu_{n}\left[  g_{n1},\ldots,g_{nn}\right]
\end{array}
\right]  =\mu_{n}\left[
\begin{array}
[c]{c}%
\mu_{n}\left[  g_{11},\ldots,g_{n1}\right] \\
\vdots\\
\mu_{n}\left[  g_{1n},\ldots,g_{nn}\right]
\end{array}
\right]  . \label{mmn}%
\end{equation}

\end{definition}

For polyadic semigroups we use the notation (\ref{mn}) and can present the
mediality as follows%
\begin{equation}
\mathbf{\mu}_{n}^{n}\left[  \mathbf{G}\right]  =\mathbf{\mu}_{n}^{n}\left[
\mathbf{G}^{T}\right]  ,
\end{equation}
where $\mathbf{G}=\left\Vert g_{ij}\right\Vert $ is the $n\times n$ matrix of
elements and $\mathbf{G}^{T}$ is its transpose. The semicommutative polyadic
semigroups are medial, as in the binary case, but, in general (except $n=3$)
not vice versa \cite{gla/gle}. A more general concept is $\sigma
$-permutability \cite{sto/dud}, such that the mediality is its particular case
with $\sigma=\left(  1,n\right)  $.

\begin{definition}
\label{def-cancel}A polyadic system is \textit{cancellative}, if%
\begin{equation}
\mu_{n}\left[  g,\mathbf{t}\right]  =\mu_{n}\left[  h,\mathbf{t}\right]
\Longrightarrow g=h,
\end{equation}
where $g,h$ can be on any place. This means that the mapping $\mu_{n}$ is
one-to-one in each variable. If $g,h$ are on the same $i$-th place on both
sides, the polyadic system is called $i$-\textit{cancellative}.
\end{definition}

The \textit{left} and \textit{right} cancellativity are $1$-cancellativity and
$n$-cancellativity respectively. A right and left cancellative $n$-ary
semigroup is \textit{cancellative} (with respect to the same subset).

\begin{definition}
A polyadic system is called (uniquely) $i$-\textit{solvable}, if for all
polyads $\mathbf{t}$, $\mathbf{u}$ and element $h$, one can (uniquely) resolve
the equation (with respect to $h$) for the fundamental operation%
\begin{equation}
\mu_{n}\left[  \mathbf{u},h,\mathbf{t}\right]  =g \label{mug}%
\end{equation}
where $h$ can be on any $i$-th place.
\end{definition}

\begin{definition}
A polyadic system which is uniquely $i$-solvable for \emph{all} places $i$ is
called a $n$-\textit{ary }(or \textit{polyadic})\textit{ quasigroup}.
\end{definition}

It follows, that, if (\ref{mug}) uniquely $i$-solvable for all places, than%
\begin{equation}
\mathbf{\mu}_{n}^{\ell_{\mu}}\left[  \mathbf{u},h,\mathbf{t}\right]  =g
\label{mugl}%
\end{equation}
can be (uniquely) resolved with respect to $h$ on any place.

\begin{definition}
\label{def-grp}An associative polyadic quasigroup is called a $n$-\textit{ary}
(or \textit{polyadic})\textit{ group}.
\end{definition}

The above definition is the most general one, but it is overdetermined. Much
work on polyadic groups was done \cite{rusakov1} to minimize the set of axioms
(solvability not in all places \cite{pos,cel1}, decreasing or increasing the
number of unknowns in determining equations \cite{galmak1}) or construction in
terms of additionally defined objects (various analogs of the identity and
sequences \cite{usa}), as well as using not total associativity, but instead
various partial ones \cite{sok,sok1,yur}.

In a polyadic group the only solution of (\ref{mug}) is called a
\textit{querelement} of $g$ and denoted by $\bar{g}$ \cite{dor3}, such that%
\begin{equation}
\mu_{n}\left[  \mathbf{h},\bar{g}\right]  =g, \label{mgg}%
\end{equation}
where $\bar{g}$ can be on any place. So, any idempotent $g$ coincides with its
querelement $\bar{g}=g$. It follows from (\ref{mgg}) and (\ref{mng}), that the
polyad%
\begin{equation}
\mathbf{n}_{g}=\left(  g^{n-2}\bar{g}\right)  \label{ng}%
\end{equation}
is neutral for any element of a polyadic group, where $\bar{g}$ can be on any
place. If this $i$-th place is important, then we write $\mathbf{n}_{g;i}$.
The number of relations in (\ref{mgg}) can be reduced from $n$ (the number of
possible places) to only $2$ (when $g$ is on the first and last places
\cite{dor3,tim72}, or on some other 2 places). In a polyadic group the
\textit{D\"{o}rnte relations}%
\begin{equation}
\mu_{n}\left[  g,\mathbf{n}_{h;i}\right]  =\mu_{n}\left[  \mathbf{n}%
_{h;j},g\right]  =g \label{mgnn}%
\end{equation}
hold true for any allowable $i,j$. In the case of a binary group the relations
(\ref{mgnn}) become $g\cdot h\cdot h^{-1}=h\cdot h^{-1}\cdot g=g$.

The relation (\ref{mgg}) can be treated as a definition of the unary
\textit{queroperation}%
\begin{equation}
\bar{\mu}_{1}\left[  g\right]  =\bar{g}, \label{m1g}%
\end{equation}
such that the diagram%
\begin{equation}
\begin{diagram} G^{\times n} &\rTo^{\mu_{n}} & G\\ \uTo^{\operatorname*{id}^{\times \left( n-1 \right)}\times \bar{\mu}_{1}}&\;\;\ruTo(3,2)_{\mathsf{Pr}_n}&\\ G^{\times n}& & \end{diagram} \label{diagm4}%
\end{equation}
commutes. Then, using the queroperation (\ref{m1g}) one can give a
\textsl{diagrammatic definition} of a polyadic group (cf. \cite{gle/gla}).

\begin{definition}
A \textit{polyadic group} is a universal algebra
\begin{equation}
\mathfrak{G}_{n}^{grp}=\left\langle G\mid\mu_{n},\bar{\mu}_{1}\mid
\text{associativity, D\"{o}rnte relations}\right\rangle ,
\end{equation}
where $\mu_{n}$ is a $n$-ary associative operation and $\bar{\mu}_{1}$ is the
queroperation, such that the following diagram%
\begin{equation}
\begin{diagram} G^{\times \left( n \right)} & \rTo^{\operatorname*{id}^{\times \left( n-1 \right)} \times\bar{\mu}_{1}}& & G^{\times n}&& \lTo^{\;\;\bar{\mu}_{1}\times\operatorname*{id}^{\times \left( n-1 \right)} } & G^{\times n} \\ \uTo^{\operatorname*{id}\times\mathsf{Diag}_{ \left( n-1 \right)}} && & \dTo_{\mu_{n}}&& & \uTo_{\mathsf{Diag}_{ \left( n-1 \right)}\times\operatorname*{id}} \\ G\times G & \rTo^{ \mathsf{Pr}_1} && G& &\lTo^{ \mathsf{Pr}_2} & G\times G\\ \end{diagram} \label{diam5}%
\end{equation}
commutes, where $\bar{\mu}_{1}$ can be only on first and second places from
the right (resp. left) on the left (resp. right) part of the diagram.
\end{definition}

A straightforward generalization of the queroperation concept and
corresponding definitions can be made by substituting in the above formulas
(\ref{mgg})--(\ref{m1g}) the $n$-ary multiplication $\mu_{n}$ by iterating the
multiplication $\mathbf{\mu}_{n}^{\ell_{\mu}}$ (\ref{mn}) (cf. \cite{dud2} for
$\ell_{\mu}=2$).

\begin{definition}
Let us define the \textit{querpower} $k$ of $g$ recursively%
\begin{equation}
\bar{g}^{\left\langle \left\langle k\right\rangle \right\rangle }%
=\overline{\left(  \bar{g}^{\left\langle \left\langle k-1\right\rangle
\right\rangle }\right)  },
\end{equation}
where $\bar{g}^{\left\langle \left\langle 0\right\rangle \right\rangle }=g$,
$\bar{g}^{\left\langle \left\langle 1\right\rangle \right\rangle }=\bar{g}$,
or as the $k$ composition $\bar{\mu}_{1}^{\circ k}=\overset{k}{\overbrace
{\bar{\mu}_{1}\circ\bar{\mu}_{1}\circ\ldots\circ\bar{\mu}_{1}}}$ of the
queroperation (\ref{m1g}).
\end{definition}

For instance \cite{galmak1}, $\bar{\mu}_{1}^{\circ2}=\mathbf{\mu}_{n}^{n-3}$,
such that for any ternary group $\bar{\mu}_{1}^{\circ2}=\operatorname*{id}$,
i.e. one has $\overline{\bar{g}}=g$. Using the queroperation in polyadic
groups we can define the \textit{negative polyadic power} of an element $g$ by
the following recursive relation%
\begin{equation}
\mu_{n}\left[  g^{\left\langle \ell_{\mu}-1\right\rangle },g^{n-2}%
,g^{\left\langle -\ell_{\mu}\right\rangle }\right]  =g,
\end{equation}
or (after use of (\ref{pp})) as a solution of the equation%
\begin{equation}
\mathbf{\mu}_{n}^{\ell_{\mu}}\left[  g^{\ell_{\mu}\left(  n-1\right)
},g^{\left\langle -\ell_{\mu}\right\rangle }\right]  =g.
\end{equation}

It is known that the querpower and the polyadic power are mutually connected
\cite{dud1}. Here, we reformulate this connection using the so called Heine
numbers \cite{heine} or $q$-deformed numbers \cite{kac/che}%
\begin{equation}
\left[  \left[  k\right]  \right]  _{q}=\dfrac{q^{k}-1}{q-1},
\end{equation}
which have the \textquotedblleft nondeformed\textquotedblright\ limit
$q\rightarrow1$ as $\left[  k\right]  _{q}\rightarrow k$. Then%
\begin{equation}
\bar{g}^{\left\langle \left\langle k\right\rangle \right\rangle }%
=g^{\left\langle -\left[  \left[  k\right]  \right]  _{2-n}\right\rangle },
\end{equation}
which can be treated as follows: the querpower coincides with the negative
polyadic \textit{deformed}\emph{ }power with a \textquotedblleft
deformation\textquotedblright\ parameter $q$ which is equal to the
\textquotedblleft deviation\textquotedblright\ $\left(  2-n\right)  $ from the
binary group.

\section{Homomorphisms of polyadic systems}

Let $\mathfrak{G}$$_{n}=\left\langle G;\mu_{n}\right\rangle $ and
$\mathfrak{G}$$_{n^{\prime}}^{\prime}=\left\langle G^{\prime};\mu_{n^{\prime}%
}^{\prime}\right\rangle $ be two polyadic systems of any kind (quasigroup,
semigroup, group, etc.). If they have the multiplications of the same arity
$n=n^{\prime}$, then one can define the mappings from $\mathfrak{G}$$_{n}$ to
$\mathfrak{G}$$_{n}^{\prime}$. Usually such polyadic systems are
\textit{similar}, and we call mappings between them \textit{the}
\textit{equiary mappings}.

Let us take $n+1$ mappings $\varphi_{i}^{GG^{\prime}}:G\rightarrow G^{\prime}%
$, $i=1,\ldots,n+1$. An ordered system of mappings $\left\{  \varphi
_{i}^{GG^{\prime}}\right\}  $ is called a \textit{homotopy} from
$\mathfrak{G}$$_{n}$ to $\mathfrak{G}$$_{n}^{\prime}$, if \cite{belousov}
\begin{equation}
\varphi_{n+1}^{GG^{\prime}}\left(  \mu_{n}\left[  g_{1},\ldots,g_{n}\right]
\right)  =\mu_{n}^{\prime}\left[  \varphi_{1}^{GG^{\prime}}\left(
g_{1}\right)  ,\ldots,\varphi_{n}^{GG^{\prime}}\left(  g_{n}\right)  \right]
,\ \ \ \ \ g_{i}\in G. \label{fm}%
\end{equation}

In general, one should add to this definition the \textquotedblleft
mapping\textquotedblright\ of the multiplications%
\begin{equation}
\mu_{n}\overset{\psi_{nn^{\prime}}^{\left(  \mu\mu^{\prime}\right)  }}%
{\mapsto}\mu_{n^{\prime}}^{\prime}. \label{mm}%
\end{equation}
In such a way, homotopy can be defined as the extended system of mappings
$\left\{  \varphi_{i}^{GG^{\prime}};\psi_{nn}^{\left(  \mu\mu^{\prime}\right)
}\right\}  $. The corresponding commutative (equiary) diagram is%
\begin{equation}
\begin{diagram} G & \rTo^{\varphi_{n+1}^{GG^{\prime}}} & G^{\prime} \\ \uTo^{\mu_{n}} & \rDotsto~{\psi_{nn}^{\left( \mu\right) }} & \uTo_{\mu_{n}^{\prime}} \\ G^{\times n} & \rTo^{ \varphi_{1}^{GG^{\prime}}\times \ldots \times \varphi_{n}^{GG^{\prime}}} & \left( G^{\prime}\right)^{\times n} \\ \end{diagram} \label{dia1}%
\end{equation}

The existence of the additional \textquotedblleft mapping\textquotedblright%
\ $\psi_{nn}^{\left(  \mu\mu^{\prime}\right)  }$ acting on the second
component of $\left\langle G;\mu_{n}\right\rangle $ is tacitly implied. We
will write/mention the \textquotedblleft mappings\textquotedblright%
\ $\psi_{nn^{\prime}}^{\left(  \mu\mu^{\prime}\right)  }$ manifestly, e.g.,%
\begin{equation}
\mathfrak{G}_{n}\overset{\left\{  \varphi_{i}^{GG^{\prime}};\psi_{nn}^{\left(
\mu\mu^{\prime}\right)  }\right\}  }{\rightrightarrows}\mathfrak{G}%
_{n^{\prime}}^{\prime}, \label{gg1}%
\end{equation}
only as needed. If all the components $\varphi_{i}^{GG^{\prime}}$ of a
homotopy are bijections, it is called an \textit{isotopy}. In case of polyadic
quasigroups \cite{belousov} all mappings $\varphi_{i}^{GG^{\prime}}$ are
usually taken as permutations of the same underlying set $G=G^{\prime}$. If
the multiplications are also coincide $\mu_{n}=\mu_{n}^{\prime}$, then
$\left\{  \varphi_{i}^{GG};\operatorname*{id}\right\}  $ is called an
\textit{autotopy} of the polyadic system $\mathfrak{G}$$_{n}$. Various
properties of homotopy in universal algebras were studied, e.g. in
\cite{pet77,hal94}.

The homotopy, isotopy and autotopy are widely used equiary mappings in the
study of polyadic quasigroups and loops, while their diagonal equiary
counterparts (all $\varphi_{i}^{GG^{\prime}}$ coincide), the homomorphism,
isomorphism and automorphism, are more suitable in investigation of polyadic
semigroups, groups and rings and their wide applications in physics. Usually,
it is written about the latter between similar (equiary) polyadic systems:
they \textquotedblleft...are so well known that we shall not bother to define
them carefully\textquotedblright\ \cite{hob/mck}. Nevertheless, we give a
diagrammatic definition of the standard homomorphism between similar polyadic
systems in our notation, which will be convenient to explain the clear way of
its generalization.

A \textit{homomorphism} from $\mathfrak{G}$$_{n}$ to $\mathfrak{G}$%
$_{n}^{\prime}$ is given, if there exists a mapping $\varphi^{GG^{\prime}%
}:G\rightarrow G^{\prime}$ satisfying
\begin{equation}
\varphi^{GG^{\prime}}\left(  \mu_{n}\left[  g_{1},\ldots,g_{n}\right]
\right)  =\mu_{n}^{\prime}\left[  \varphi^{GG^{\prime}}\left(  g_{1}\right)
,\ldots,\varphi^{GG^{\prime}}\left(  g_{n}\right)  \right]  ,\ \ \ \ \ g_{i}%
\in G, \label{fm1}%
\end{equation}
which means that the corresponding (equiary) diagram is commutative%

\begin{equation}
\begin{diagram} G & \rTo^{\varphi^{GG^{\prime}}} & G^{\prime} \\ \uTo^{\mu_{n}} & \rDotsto~{\psi_{nn}^{\left( \mu\mu^{\prime}\right) }} & \uTo_{\mu_{n}^{\prime}} \\ G^{\times n} & \rTo^{ \left( \varphi^{GG^{\prime}}\right) ^{\times n}} & \left( G^{\prime}\right)^{\times n} \\ \end{diagram} \label{dia2}%
\end{equation}

Usually the homomorphism is denoted by the same one letter $\varphi
^{GG^{\prime}}$, while it would be more consistent to use for its notation the
extended \textsl{pair} of mappings $\left\{  \varphi^{GG^{\prime}};\psi
_{nn}^{\left(  \mu\mu^{\prime}\right)  }\right\}  $. We will use both
notations on a par.

We first mention a small subset of known generalizations of the homomorphism
(for bibliography till 1982 see, e.g., \cite{gla/gle1}) and then introduce a
concrete construction for an analogous mapping which can change the arity of
the multiplication (fundamental operation) without introducing additional
(term) operations. A general approach to mappings between free algebraic
systems was initiated in \cite{fuj59}, where the so-called basic mapping
formulas for generators were introduced, and its generalization to many-sorted
algebras was given in \cite{vidal-tur}. In \cite{nov02} it was shown that the
construction of all homomorphisms between similar polyadic systems can be
reduced to some homomorphisms between corresponding mono-unary algebras
\cite{nov90}. The notion of $n$-ary homomorphism is realized as a sequence of
$n$ consequent homomorphisms $\varphi_{i}$, $i=1,\ldots,n$, of $n$ similar
polyadic systems%
\begin{equation}
\overset{n}{\overbrace{\mathfrak{G}_{n}\overset{\varphi_{1}}{\rightarrow
}\mathfrak{G}_{n}^{\prime}\overset{\varphi_{2}}{\rightarrow}\ldots
\overset{\varphi_{n-1}}{\rightarrow}\mathfrak{G}_{n}^{\prime\prime}%
\overset{\varphi_{n}}{\rightarrow}\mathfrak{G}_{n}^{\prime\prime\prime}}}
\label{gfg}%
\end{equation}
(generalizing Post's $n$-adic substitutions \cite{pos}) was introduced in
\cite{galm98}, and studied in \cite{galm01,galmak2}.

The above constructions do not change the arity of polyadic systems, because
they are based on the corresponding diagram which gives a definition of an
\textit{equiary} mapping. To change arity one has to:

1) add another equiary diagram with \textsl{additional} operations using
\textsl{the same} formula (\ref{fm1}), where both do not change arity;

2) use \textsl{one} modified (and not equiary) diagram and the underlying
formula (\ref{fm1}) by themselves, which will allow us to change arity
\textsl{without} introducing additional operations.

The first way leads to the concept of \textsl{weak homomorphism} which was
introduced in \cite{goe66,mar66,gla/mic1} for non-indexed algebras and in
\cite{gla80} for indexed algebras, then developed in \cite{tra65} for Boolean
and Post algebras, in \cite{den/wis} for coalgebras and $F$-algebras
\cite{dan/sae} (see also \cite{chu/smi}). To define the weak homomorphism in
our notation we should incorporate into the polyadic systems $\left\langle
G;\mu_{n}\right\rangle $ and $\left\langle G^{\prime};\mu_{n^{\prime}}%
^{\prime}\right\rangle $ the following additional term operations of opposite
arity $\nu_{n^{\prime}}:G^{\times n^{\prime}}\rightarrow G$ and $\nu
_{n}^{\prime}:G^{\prime\times n}\rightarrow G^{\prime}$ and consider
\textsl{two} equiary mappings between $\left\langle G;\mu_{n},\nu_{n^{\prime}%
}\right\rangle $ and $\left\langle G^{\prime};\mu_{n^{\prime}}^{\prime}%
,\nu_{n}^{\prime}\right\rangle $.

A \textit{weak homomorphism} from $\left\langle G;\mu_{n},\nu_{n^{\prime}%
}\right\rangle $ to $\left\langle G^{\prime},\mu_{n^{\prime}}^{\prime},\nu
_{n}^{\prime}\right\rangle $ is given, if there exists a mapping
$\varphi^{GG^{\prime}}:G\rightarrow G^{\prime}$ satisfying \textsl{two}
relations simultaneously%
\begin{align}
\varphi^{GG^{\prime}}\left(  \mu_{n}\left[  g_{1},\ldots,g_{n}\right]
\right)   &  =\nu_{n}^{\prime}\left[  \varphi^{GG^{\prime}}\left(
g_{1}\right)  ,\ldots,\varphi^{GG^{\prime}}\left(  g_{n}\right)  \right]
,\label{wh1}\\
\varphi^{GG^{\prime}}\left(  \nu_{n^{\prime}}\left[  g_{1},\ldots
,g_{n^{\prime}}\right]  \right)   &  =\mu_{n^{\prime}}^{\prime}\left[
\varphi^{GG^{\prime}}\left(  g_{1}\right)  ,\ldots,\varphi^{GG^{\prime}%
}\left(  g_{n^{\prime}}\right)  \right]  , \label{wh2}%
\end{align}
which means that two \textsl{equiary} diagrams commute%

\begin{equation}
\begin{diagram} G & \rTo^{\varphi^{GG^{\prime}}} & G^{\prime} \\ \uTo^{\mu_{n}} & \rDotsto~{\psi_{nn}^{\left( \mu\nu^{\prime}\right) }} & \uTo_{\nu_{n}^{\prime}} \\ G^{\times n} & \rTo^{ \left( \varphi^{GG^{\prime}}\right) ^{\times n}} & \left( G^{\prime}\right)^{\times n} \\ \end{diagram}\ \ \ \ \ \ \ \ \ \ \ \ \ \ \ \begin{diagram} G & \rTo^{\varphi^{GG^{\prime}}} & G^{\prime} \\ \uTo^{\nu_{n^{\prime}}} & \rDotsto~{\psi_{n^{\prime}n^{\prime}}^{\left( \nu\mu^{\prime}\right) }} & \uTo_{\mu_{n^{\prime}}^{\prime}} \\ G^{\times n^{\prime}} & \rTo^{ \left( \varphi^{GG^{\prime}}\right) ^{\times n^{\prime}}} & \left( G^{\prime}\right)^{\times n^{\prime}} \\ \end{diagram}
\end{equation}

If only one of the relations (\ref{wh1}) or (\ref{wh2}) holds, such a mapping
is called a \textit{semi-weak homomorphism} \cite{kol84}. If $\varphi
^{GG^{\prime}}$ is bijective, then it defines a \textit{weak isomorphism}. Any
weak epimorphism can be decomposed into a homomorphism and a weak isomorphism
\cite{gla/mic}, and therefore the study of weak homomorphisms reduces to weak
isomorphisms (see also \cite{cza62,mal57,mal58}).

\section{Multiplace mappings of polyadic systems and heteromorphisms}

Let us turn to the second way of changing the arity of the multiplication and
use only one relation which we then modify in some natural manner. First,
recall that in any set $G$ there always exists the additional distinguished
mapping, viz. the identity $\operatorname*{id}\nolimits_{G}$. We use the
multiplication $\mu_{n}$ with its combination of $\operatorname*{id}%
\nolimits_{G}$. We define an ($\ell_{\operatorname*{id}}$-\textit{intact}%
)\textit{ id-product} for the polyadic system $\left\langle G;\mu
_{n}\right\rangle $ as%
\begin{align}
\mu_{n}^{\left(  \ell_{\operatorname*{id}}\right)  }  &  =\mu_{n}\times\left(
\operatorname*{id}\nolimits_{G}\right)  ^{\times\ell_{\operatorname*{id}}%
},\label{mid}\\
\mu_{n}^{\left(  \ell_{\operatorname*{id}}\right)  }  &  :G^{\times\left(
n+\ell_{\operatorname*{id}}\right)  }\rightarrow G^{\times\left(
1+\ell_{\operatorname*{id}}\right)  }. \label{mid1}%
\end{align}

To indicate the exact $i$-th place of $\mu_{n}$ in the r.h.s. of (\ref{mid}),
we write $\mu_{n}^{\left(  \ell_{\operatorname*{id}}\right)  }\left(
i\right)  $, as needed. Here we use the id-product to generalize the
homomorphism and consider mappings between polyadic systems of different
arity. It follows from (\ref{mid1}) that, if the image of the id-product is
$G$ alone, than $\ell_{\operatorname*{id}}=0$. Let us introduce a
\textit{multiplace mapping} $\Phi_{k}^{\left(  n,n^{\prime}\right)  }$ acting
as follows%
\begin{equation}
\Phi_{k}^{\left(  n,n^{\prime}\right)  }:G^{\times k}\rightarrow G^{\prime}.
\label{fgg}%
\end{equation}

While constructing the corresponding diagram, we are allowed to take only one
upper $\Phi_{k}^{\left(  n,n^{\prime}\right)  }$, because of one $G^{\prime}$
in the upper right corner. Since we already know that the lower right corner
is exactly $G^{\prime\times n^{\prime}}$ (as a pre-image of one multiplication
$\mu_{n^{\prime}}^{\prime}$), the lower horizontal arrow should be a product
of $n^{\prime}$ multiplace mappings $\Phi_{k}^{\left(  n,n^{\prime}\right)  }%
$. So we can write a definition of a multiplace analog of homomorphisms which
changes the arity of the multiplication using one relation.

\begin{definition}
A $k$-\textit{place} \textit{heteromorphism} from $\mathfrak{G}$$_{n}$ to
$\mathfrak{G}$$_{n^{\prime}}^{\prime}$ is given, if there exists a $k$-place
mapping $\Phi_{k}^{\left(  n,n^{\prime}\right)  }$ (\ref{fgg}) such that the
following (arity changing or \textit{unequiary}) diagram is commutative%

\begin{equation}
\begin{diagram} G^{\times k} & \rTo^{\Phi_{k}} & G^{\prime} \\ \uTo^{\mu_{n}^{\left( \ell_{\operatorname*{id}}\right) }} & & \uTo_{\mu_{n^{\prime}}^{\prime}} \\ G^{\times kn^{\prime}} & \rTo^{ \left( \Phi_{k}\right) ^{\times n^{\prime}}} & \left( G^{\prime}\right)^{\times n^{\prime}} \\ \end{diagram} \label{dia3}%
\end{equation}
and the corresponding defining equation (a modification of (\ref{fm1}))
depends on the place $i$ of $\mu_{n}$ in (\ref{mid}).
\end{definition}

For $i=1$ a heteromorphism is defined by the formula%
\begin{equation}
\Phi_{k}^{\left(  n,n^{\prime}\right)  }\left(
\begin{array}
[c]{c}%
\mu_{n}\left[  g_{1},\ldots,g_{n}\right] \\
g_{n+1}\\
\vdots\\
g_{n+\ell_{\operatorname*{id}}}%
\end{array}
\right)  =\mu_{n^{\prime}}^{\prime}\left[  \Phi_{k}^{\left(  n,n^{\prime
}\right)  }\left(
\begin{array}
[c]{c}%
g_{1}\\
\vdots\\
g_{k}%
\end{array}
\right)  ,\ldots,\Phi_{k}^{\left(  n,n^{\prime}\right)  }\left(
\begin{array}
[c]{c}%
g_{k\left(  n^{\prime}-1\right)  }\\
\vdots\\
g_{kn^{\prime}}%
\end{array}
\right)  \right]  . \label{fff}%
\end{equation}

The notion ``heteromorphism'' is motivated by \cite{ell06,ell07a}, where
mappings between objects from different categories were considered and called
``chimera morphisms''. See, also, \cite{pec11}.

In the particular case $n=3$, $n^{\prime}=2$, $k=2$, $\ell_{\operatorname*{id}%
}=1$ we have%
\begin{equation}
\Phi_{2}^{\left(  3,2\right)  }\left(
\begin{array}
[c]{c}%
\mu_{3}\left[  g_{1},g_{2},g_{3}\right] \\
g_{4}%
\end{array}
\right)  =\mu_{2}^{\prime}\left[  \Phi_{2}^{\left(  3,2\right)  }\left(
\begin{array}
[c]{c}%
g_{1}\\
g_{2}%
\end{array}
\right)  ,\Phi_{2}^{\left(  3,2\right)  }\left(
\begin{array}
[c]{c}%
g_{3}\\
g_{4}%
\end{array}
\right)  \right]  . \label{f2}%
\end{equation}

This formula was used in the construction of the bi-element representations of
ternary groups \cite{bor/dud/dup3}. Consider the example.

\begin{example}
\label{exam-matr}Let $G=M_{2}^{adiag}\left(  \mathbb{K}\right)  $, a set of
antidiagonal $2\times2$ matrices over the field $\mathbb{K}$ and $G^{\prime
}=\mathbb{K}$, where $\mathbb{K}=\mathbb{R},\mathbb{C},\mathbb{Q},\mathbb{H}$.
The ternary multiplication $\mu_{3}$ is a product of 3 matrices. Obviously,
$\mu_{3}$ is not derived from a binary multiplication. For the elements
$g_{i}=\left(
\begin{array}
[c]{cc}%
0 & a_{i}\\
b_{i} & 0
\end{array}
\right)  $, $i=1,2$, we construct a 2-place mapping $G\times G\rightarrow
G^{\prime}$ as%
\begin{equation}
\Phi_{2}^{\left(  3,2\right)  }\left(
\begin{array}
[c]{c}%
g_{1}\\
g_{2}%
\end{array}
\right)  =a_{1}a_{2}b_{1}b_{2}, \label{f22}%
\end{equation}
which is a heteromorphism, because it satisfies (\ref{f2}). Let us introduce a
standard 1-place mapping by $\varphi\left(  g_{i}\right)  =a_{i}b_{i}, i=1,2$.
It is important to note, that $\varphi\left(  g_{i}\right)  $ is not a
homomorphism, because the product $g_{1}g_{2}$ belongs to diagonal matrices.
Consider the product of mappings
\begin{equation}
\varphi\left(  g_{1}\right)  \cdot\varphi\left(  g_{2}\right)  =a_{1}%
b_{1}a_{2}b_{2} , \label{f2ff}%
\end{equation}
where the product $\left(  \cdot\right)  $ in l.h.s. is taken in $\mathbb{K}$.
We observe that (\ref{f22}) and (\ref{f2ff}) coincide for the commutative
field $\mathbb{K}$ only ($=\mathbb{R},\mathbb{C}$) only, and in this case we
can have the relation between the heteromorhism $\Phi_{2}^{\left(  3,2\right)
}$ and the 1-place mapping $\varphi$
\begin{equation}
\Phi_{2}^{\left(  3,2\right)  }\left(
\begin{array}
[c]{c}%
g_{1}\\
g_{2}%
\end{array}
\right)  =\varphi\left(  g_{1}\right)  \cdot\varphi\left(  g_{2}\right)  ,
\label{f2ff1}%
\end{equation}
while for the noncommutative field $\mathbb{K}$ ($=\mathbb{Q}$ or $\mathbb{H}%
$) there is no such relation.
\end{example}

A heteromorphism is called \textit{derived}, if it can be expressed through a
1-place mapping (not necessary a homomorphism). So, in the above
\emph{Example} \ref{exam-matr} the heteromorphism is derived (by formula
(\ref{f2ff1})) for the commutative field $\mathbb{K}$ and nonderived for the
noncommutative $\mathbb{K}$.

For arbitrary $n$ a slightly modified construction (\ref{f2}) with still
binary final arity, defined by $n^{\prime}=2$, $k=n-1$, $\ell
_{\operatorname*{id}}=n-2$,%
\begin{equation}
\Phi_{n-1}^{\left(  n,2\right)  }\left(
\begin{array}
[c]{c}%
\mu_{n}\left[  g_{1},\ldots,g_{n-1},h_{1}\right] \\
h_{2}\\
\vdots\\
h_{n-1}%
\end{array}
\right)  =\mu_{2}^{\prime}\left[  \Phi_{n-1}^{\left(  n,2\right)  }\left(
\begin{array}
[c]{c}%
g_{1}\\
\vdots\\
g_{n-1}%
\end{array}
\right)  ,\Phi_{n-1}^{\left(  n,2\right)  }\left(
\begin{array}
[c]{c}%
h_{1}\\
\vdots\\
h_{n-1}%
\end{array}
\right)  \right]  . \label{fn}%
\end{equation}
was used in \cite{dud07} to study representations of $n$-ary groups. However,
no new results compared with \cite{bor/dud/dup3} (other than changing $3$ to
$n$ in some formulas) were obtained. This reflects the fact that a major role
is played by the final arity $n^{\prime}$ and the number of $n$-ary
multiplications in the l.h.s. of (\ref{f2}) and (\ref{fn}). In the above
cases, the latter number was one, but can make it arbitrary below $n$.

\begin{definition}
A heteromorphism is called a $\ell_{\mu}$-\textit{ple heteromorphism}, if it
contains $\ell_{\mu}$ multiplications in the argument of $\Phi_{k}^{\left(
n,n^{\prime}\right)  }$ in its defining relation.
\end{definition}

According this definition the mapping defined by (\ref{fff}) is the 1$_{\mu}%
$-ple heteromorphism. So by analogy with (\ref{mid})--(\ref{mid1}) we define a
$\ell_{\mu}$-\textit{ple }$\ell_{\operatorname*{id}}$-\textit{intact
id-product} for the polyadic system $\left\langle G;\mu_{n}\right\rangle $ as%
\begin{align}
\mu_{n}^{\left(  \ell_{\mu},\ell_{\operatorname*{id}}\right)  }  &  =\left(
\mu_{n}\right)  ^{\times\ell_{\mu}}\times\left(  \operatorname*{id}%
\nolimits_{G}\right)  ^{\times\ell_{\operatorname*{id}}},\label{midl}\\
\mu_{n}^{\left(  \ell_{\mu},\ell_{\operatorname*{id}}\right)  }  &
:G^{\times\left(  n\ell_{\mu}+\ell_{\operatorname*{id}}\right)  }\rightarrow
G^{\times\left(  \ell_{\mu}+\ell_{\operatorname*{id}}\right)  }. \label{midl1}%
\end{align}

\begin{definition}
A $\ell_{\mu}$-\textit{ple }$k$-\textit{place} \textit{heteromorphism} from
$\mathfrak{G}$$_{n}$ to $\mathfrak{G}$$_{n^{\prime}}^{\prime}$ is given, if
there exists a $k$-place mapping $\Phi_{k}^{\left(  n,n^{\prime}\right)  }$
(\ref{fgg}) such that the following \textit{unequiary} diagram is commutative%
\begin{equation}
\begin{diagram} G^{\times k} & \rTo^{\Phi_{k}} & G^{\prime} \\ \uTo^{\mu_{n}^{\left( \ell_{\mu},\ell_{\operatorname*{id}}\right) }} & & \uTo_{\mu_{n^{\prime}}^{\prime}} \\ G^{\times kn^{\prime}} & \rTo^{ \left( \Phi_{k}\right) ^{\times n^{\prime}}} & \left( G^{\prime}\right)^{\times n^{\prime}} \\ \end{diagram} \label{dia4}%
\end{equation}
The corresponding \textit{main heteromorphism equation} is%
\begin{equation}
\Phi_{k}^{\left(  n,n^{\prime}\right)  }\left(
\genfrac{}{}{0pt}{}{\left.
\begin{array}
[c]{c}%
\mu_{n}\left[  g_{1},\ldots,g_{n}\right]  ,\\
\vdots\\
\mu_{n}\left[  g_{n\left(  \ell_{\mu}-1\right)  },\ldots,g_{n\ell_{\mu}%
}\right]
\end{array}
\right\}  \ell_{\mu}}{\left.
\begin{array}
[c]{c}%
g_{n\ell_{\mu}+1},\\
\vdots\\
g_{n\ell_{\mu}+\ell_{\operatorname*{id}}}%
\end{array}
\right\}  \ell_{\operatorname*{id}}}%
\right)  =\mu_{n^{\prime}}^{\prime}\left[  \Phi_{k}^{\left(  n,n^{\prime
}\right)  }\left(
\begin{array}
[c]{c}%
g_{1}\\
\vdots\\
g_{k}%
\end{array}
\right)  ,\ldots,\Phi_{k}^{\left(  n,n^{\prime}\right)  }\left(
\begin{array}
[c]{c}%
g_{k\left(  n^{\prime}-1\right)  }\\
\vdots\\
g_{kn^{\prime}}%
\end{array}
\right)  \right]  . \label{he}%
\end{equation}

\end{definition}

Obviously, we can consider various permutations of the multiplications on both
sides, as further additional demands (associativity, commutativity, etc.), are
introduced, which will be considered below. The commutativity of the diagram
(\ref{dia4}) leads to the system of equation connecting initial and final
arities%
\begin{align}
kn^{\prime}  &  =n\ell_{\mu}+\ell_{\operatorname*{id}},\label{k1}\\
k  &  =\ell_{\mu}+\ell_{\operatorname*{id}}. \label{k2}%
\end{align}

Excluding $\ell_{\mu}$ or $\ell_{\operatorname*{id}}$, we obtain two
\textit{arity changing formulas}, respectively%
\begin{align}
n^{\prime}  &  =n-\dfrac{n-1}{k}\ell_{\operatorname*{id}},\label{n1}\\
n^{\prime}  &  =\dfrac{n-1}{k}\ell_{\mu}+1, \label{n2}%
\end{align}
where $\tfrac{n-1}{k}\ell_{\operatorname*{id}}\geq1$ and $\tfrac{n-1}{k}%
\ell_{\mu}\geq1$ are integer.

\begin{figure}[ptb]
\centering
\includegraphics[width=0.4\textwidth,keepaspectratio]{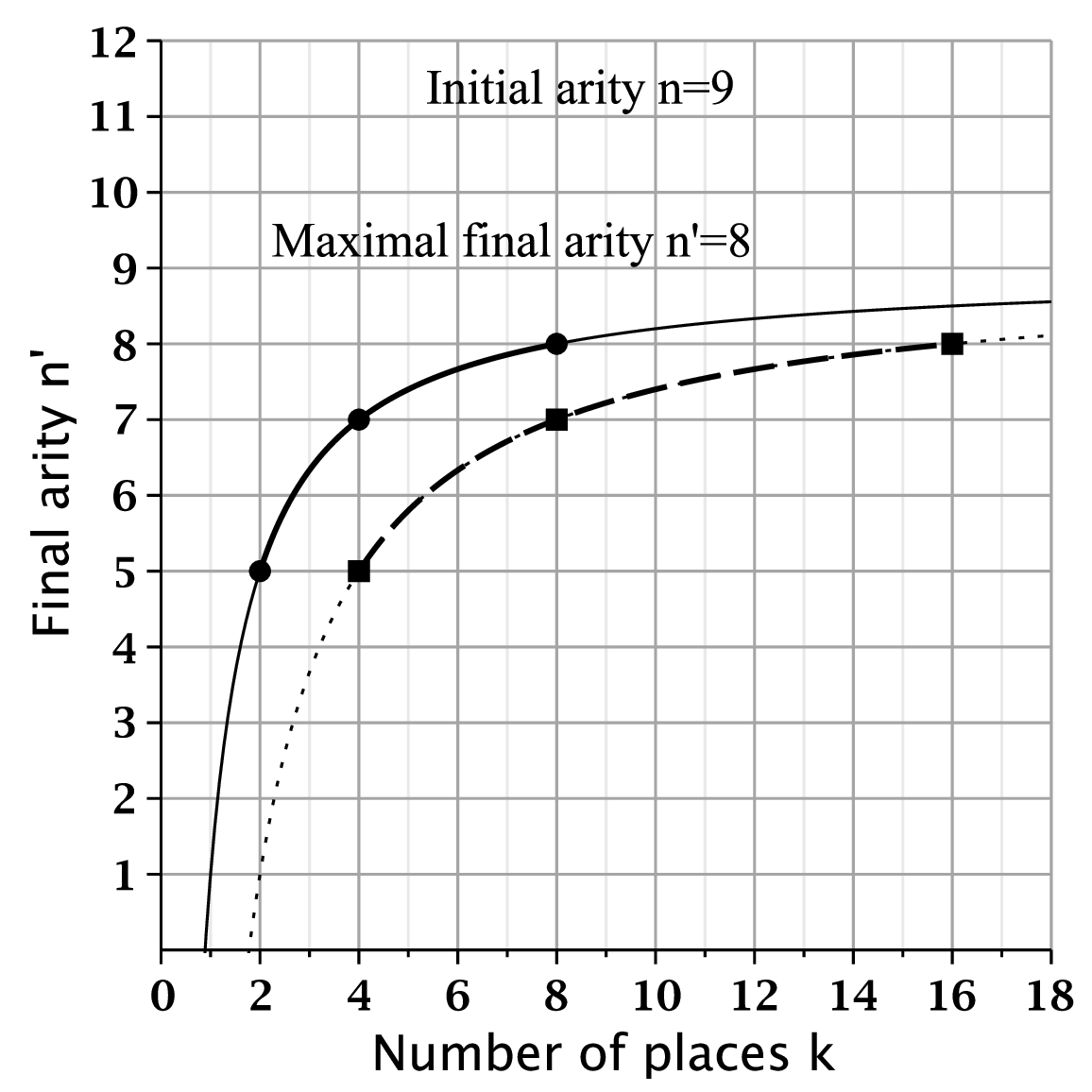}
\hskip 1.2cm
\includegraphics[width=0.4\textwidth,keepaspectratio]{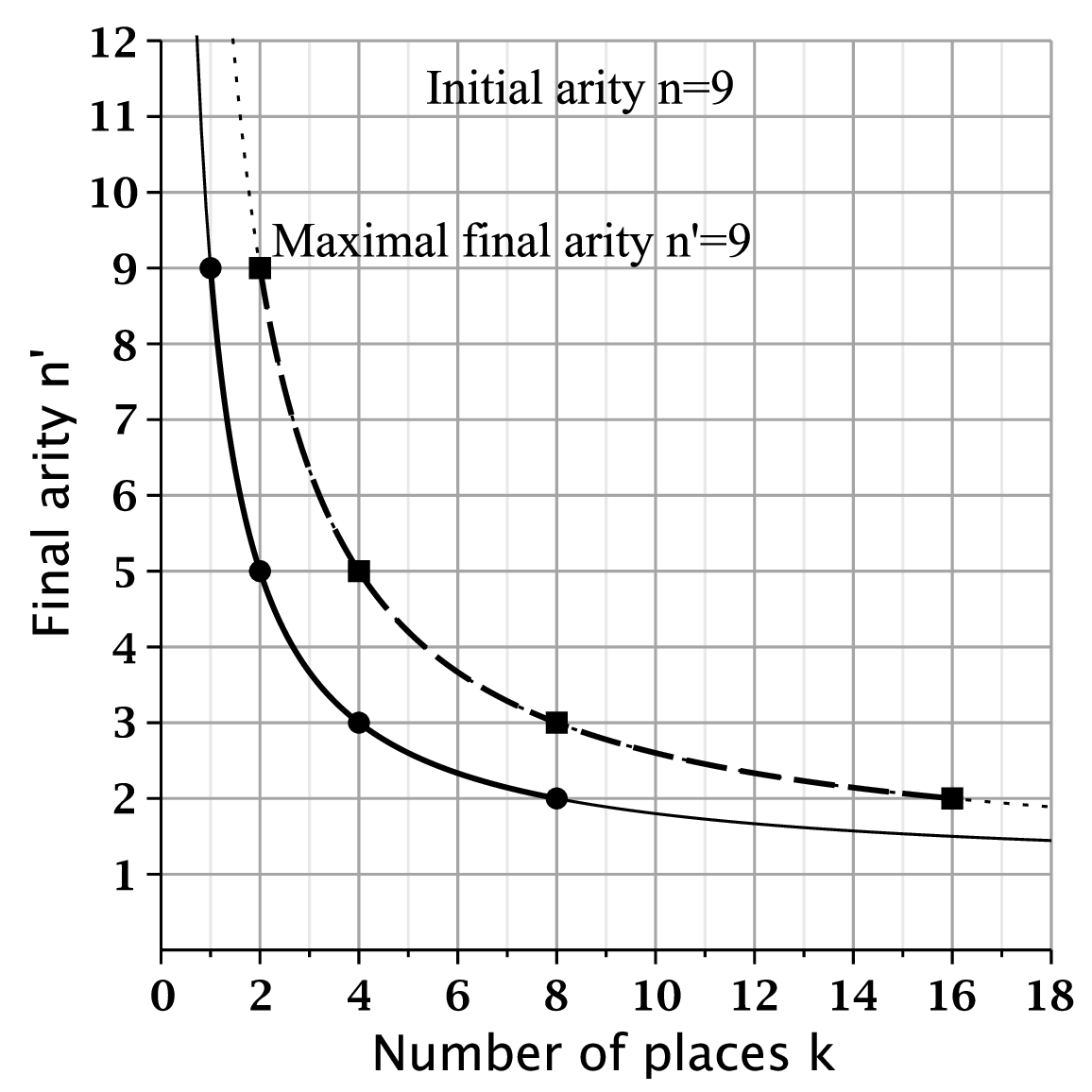}
\caption{{\protect\small Dependence of the final arity $n^{\prime}$ through
the number of heteromorphism places $k$ for the initial arity $n=9$ with the
fixed number of intact elements $\ell_{\operatorname*{id}}$ (left) and the
fixed number of multiplications $\ell_{\mu}$ (right): =1 (solid curves), =2
(dash curves) }}%
\label{F1}%
\end{figure}

As an example, the dependences $n^{\prime}\left(  k\right)  $ for the fixed
$\ell_{\mu}=1,2$ and $\ell_{\operatorname*{id}}=1,2$ with $n=9$ are presented
on Figure \ref{F1}.

The following inequalities hold valid%
\begin{align}
1  &  \leq\ell_{\mu}\leq k,\label{l1}\\
0  &  \leq\ell_{\operatorname*{id}}\leq k-1,\label{l2}\\
\ell_{\mu}  &  \leq k\leq\left(  n-1\right)  \ell_{\mu},\label{lk}\\
2  &  \leq n^{\prime}\leq n, \label{nn}%
\end{align}
which are important for the further classification of heteromorphisms. The
main statement follows from (\ref{nn})

\begin{proposition}
The heteromorphism $\Phi_{k}^{\left(  n,n^{\prime}\right)  }$ defined by the
general relation (\ref{he}) always decreases the arity of polyadic multiplication.
\end{proposition}

Another important observation is the fact that only the id-product
(\ref{midl}) with $\ell_{\operatorname*{id}}\neq0$ leads to a change of the
arity. In the extreme case, when $k$ approaches its minimum, $k=k_{\min}%
=\ell_{\mu}$, the final arity approaches its maximum $n_{\max}^{\prime}=n$,
and the id-product becomes a product of $\ell_{\mu}$ initial multiplications
$\mu_{n}$ without id's, since now $\ell_{\operatorname*{id}}=0$ in (\ref{he}).
Therefore, we call a heteromorphism defined by (\ref{he}) with $\ell
_{\operatorname*{id}}=0$ \textit{a }$k\left(  =\ell_{\mu}\right)
$\textit{-place homomorphism}. The ordinary homomorphism (\ref{fm})
corresponds to $k=\ell_{\mu}=1$, and so it is really a 1-place homomorphism.
An opposite extreme case, when the final arity approaches its minimum
$n_{\min}^{\prime}=2$ (the final operation is binary), corresponds to the
maximal value of $k$, that is $k=k_{\max}=\left(  n-1\right)  \ell_{\mu}$. The
number of id's now is $\ell_{\operatorname*{id}}=\left(  n-2\right)  \ell
_{\mu}\geq0$, which vanishes, when the initial operation is binary as well.
This is the case of the ordinary homomorphism (\ref{fm}) for both binary
operations $n^{\prime}=n=2$ and $k=\ell_{\mu}=1$. We conclude that:

\emph{Any polyadic system can be mapped into a binary system} by means of the
special $k$-place $\ell_{\mu}$-ple heteromorphism $\Phi_{k}^{\left(
n,n^{\prime}\right)  }$, where $k=\left(  n-1\right)  \ell_{\mu}$ (we call it
a \textit{binarizing heteromorphism}) which is defined by (\ref{he}) with
$\ell_{\operatorname*{id}}=\left(  n-2\right)  \ell_{\mu}$.

In relation to the Gluskin-Hossz\'{u} theorem \cite{glu1} (any $n$-ary group
can be constructed from the special binary group and its homomorphism) our
statement can be treated as:

\begin{theorem}
Any $n$-ary system can be mapped into a binary system, using a suitable
binarizing heteromorphism $\Phi_{k}^{\left(  n,2\right)  }$ (\ref{he}).
\end{theorem}

The case of $1$-ple binarizing heteromorphism ($\ell_{\mu}=1$) corresponds to
the formula (\ref{fn}). Further requirements (associativity, commutativity,
etc.) will give additional relations between multiplications and $\Phi
_{k}^{\left(  n,n^{\prime}\right)  }$, and fix the exact structure of
(\ref{he}). Thus, we arrive to the following

\begin{proposition}
Classification of $\ell_{\mu}$-ple heteromorphisms:

\begin{enumerate}
\item $n^{\prime}=n_{\max}^{\prime}=n\Longrightarrow\Phi_{k}^{\left(
n,n\right)  }$ is the $\ell_{\mu}$-place or \textit{multiplace homomorphism,}
i.e.,%
\begin{equation}
k=k_{\min}=\ell_{\mu}. \label{kmin}%
\end{equation}

\item $2<n^{\prime}<n\Longrightarrow\Phi_{k}^{\left(  n,n^{\prime}\right)  }$
is the \textit{intermediate heteromorphism} with
\begin{equation}
k=\dfrac{n-1}{n^{\prime}-1}\ell_{\mu}. \label{kl}%
\end{equation}
In this case the number of intact elements is proportional to the number of
multiplications%
\begin{equation}
\ell_{\operatorname*{id}}=\dfrac{n-n^{\prime}}{n^{\prime}-1}\ell_{\mu}.
\label{ll}%
\end{equation}

\item $n^{\prime}=n_{\min}^{\prime}=2\Longrightarrow\Phi_{k}^{\left(
n,2\right)  }$ is the $\left(  n-1\right)  \ell_{\mu}$-place (multiplace)
\textit{binarizing heteromorphism}, i.e.,%
\begin{equation}
k=k_{\max}=\left(  n-1\right)  \ell_{\mu}. \label{kmax}%
\end{equation}

\end{enumerate}
\end{proposition}

In the extreme (first and third) cases there are no restrictions on the
initial arity $n$, while in the intermediate case $n$ is \textquotedblleft
quantized\textquotedblright\ due to the fact that fractions in (\ref{n1}) and
(\ref{n2}) should be integers.

Observe, that in the extreme (first and third) cases there are no restrictions
on the initial arity $n$, while in the intermediate case $n$ is
\textquotedblleft quantized\textquotedblright\ due to the fact that fractions
in (\ref{n1}) and (\ref{n2}) should be integer. In this way, we obtain the
\textsc{Table \ref{T1}} for the series of $n$ and $n^{\prime}$ (we list only
first ones, just for $2\leq k\leq4$ and include the binarizing case
$n^{\prime}=2$ for completeness).

\begin{table}[h]
\caption{\textquotedblleft Quantization\textquotedblright\ of heteromorphisms}%
\label{T1}
\begin{center}%
\begin{tabular}
[c]{|c|c|c|c|}\hline
$k$ & $\ell_{\mu}$ & $\ell_{\operatorname*{id}}$ & $n/n^{\prime}%
$\\\hline\hline
$2$ & $1$ & $1$ & $%
\begin{array}
[c]{ccccc}%
n= & 3, & 5, & 7, & \ldots\\
n^{\prime}= & 2, & 3, & 4, & \ldots
\end{array}
$\\\hline
$3$ & $1$ & $2$ & $%
\begin{array}
[c]{ccccc}%
n= & 4, & 7, & 10, & \ldots\\
n^{\prime}= & 2, & 3, & 4, & \ldots
\end{array}
$\\\hline
$3$ & $2$ & $1$ & $%
\begin{array}
[c]{ccccc}%
n= & 4, & 7, & 10, & \ldots\\
n^{\prime}= & 3, & 5, & 7, & \ldots
\end{array}
$\\\hline
$4$ & $1$ & $3$ & $%
\begin{array}
[c]{ccccc}%
n= & 5, & 9, & 13, & \ldots\\
n^{\prime}= & 2, & 3, & 4, & \ldots
\end{array}
$\\\hline
$4$ & $2$ & $2$ & $%
\begin{array}
[c]{ccccc}%
n= & 3, & 5, & 7, & \ldots\\
n^{\prime}= & 2, & 3, & 4, & \ldots
\end{array}
$\\\hline
$4$ & $3$ & $1$ & $%
\begin{array}
[c]{ccccc}%
n= & 5, & 9, & 13, & \ldots\\
n^{\prime}= & 4, & 7, & 10, & \ldots
\end{array}
$\\\hline
\end{tabular}
\end{center}
\end{table}

Thus, we have established a general structure and classification of
heteromorphisms defined by (\ref{he}). The next important issue is the
preservation of special properties (associativity, commutativity, etc.), while
passing from $\mu_{n}$ to $\mu_{n^{\prime}}^{\prime}$, which will further
restrict the concrete shape of the main relation (\ref{he}) for each choice of
the \textit{heteromorphism parameters}: arities $n$, $n^{\prime}$, places $k$,
number of intacts $\ell_{\operatorname*{id}}$ and multiplications $\ell_{\mu}$.

\section{Associativity quivers and heteromorphisms}

The most important property of the heteromorphism, which is needed for its
next applications to representation theory, is the associativity of the final
operation $\mu_{n^{\prime}}^{\prime}$, when the initial operation $\mu_{n}$ is
associative. In other words, we consider here the concrete form of semigroup
heteromorphisms. In general, this is a complicated task, because it is not
clear from (\ref{he}), which permutation in the l.h.s. should be taken to get
an associative product in its r.h.s. for each set of the heteromorphism
parameters. Straightforward checking of the associativity of the final
operation $\mu_{n^{\prime}}^{\prime}$ for each permutation in the l.h.s. of
(\ref{he}) is almost impossible, especially for higher $n$. To solve this
difficulty we introduce the concept of the associative polyadic quiver and
special rules to construct the associative final operation $\mu_{n^{\prime}%
}^{\prime}$.

\begin{definition}
A \textit{polyadic quiver of products} is the set of elements from
$\mathfrak{G}_{n}$ (presented as several copies of some matrix of the elements
glued together) and arrows, such that the elements along arrows form $n$-ary
products $\mu_{n}$.
\end{definition}

For instance, for the 4-ary multiplication $\mu_{4}\left[  g_{1},h_{2}%
,g_{2},u_{1}\right]  $ (elements from $\mathfrak{G}$$_{n}$ are arbitrary here)
a corresponding 4-adic quiver will be denoted by $\left\{  g_{1}\rightarrow
h_{2}\rightarrow g_{2}\rightarrow u_{1}\right\}  $, and graphically this
4-adic quiver is%
\begin{equation}
{\scriptsize
\begin{array}
[c]{c}%
\scriptsize\xymatrix@=10pt{
g_1\ar@{->}[dr] & h_1 &u _1 &g_1 & h_1& u_1 \\
g_2 &h_2 \ar@/^/[rr]& u_2&g_2\ar@{->}[urr] & h_2&u_2
\save"1,1"."2,3"*+[F.]\frm{}
\restore\save"1,4"."2,6"*+[F.]\frm{}
\restore}
\\[-4pt]
\xymatrix@=7pt{\\
*+[F-:<4pt>]\txt{corr}
\ar@2{->}[u]\ar@2{->}[d]
\\
\\
}%
\\[-3pt]%
\mu_{4}\left[  g_{1},h_{2},g_{2},u_{1}\right]  .
\end{array}
} \label{qv1}%
\end{equation}

Next we define \textit{polyadic quivers} which are related to the main
heteromorphism equation (\ref{he}) in the following way:

1) the matrix of elements is the transposed matrix from the r.h.s. of
(\ref{he}), such that different letters correspond to their place in $\Phi
_{k}^{\left(  n,n^{\prime}\right)  }$ and the low index of an element is
related to its position in the $\mu_{n^{\prime}}^{\prime}$ product;

2) the number of polyadic quivers is $\ell_{\mu}$, which corresponds to
$\ell_{\mu}$ multiplications in the l.h.s. of (\ref{he});

3) the heteromorphism parameters ($n$, $n^{\prime}$, $k$, $\ell
_{\operatorname*{id}}$ and $\ell_{\mu}$) are not arbitrary, but satisfy the
arity changing formulas (\ref{n1})-(\ref{n2});

4) the intact elements will be placed after a semicolon.

In this way, a polyadic quiver makes a clear visualization of the main
heteromorphism equation (\ref{he}), and later on it will allow us to
distinguish associativity preserving heteromorphisms by precise graphical rules.

For example, the polyadic quiver $\left\{  g_{1}\rightarrow h_{2}\rightarrow
g_{2}\rightarrow u_{1};h_{1},u_{2}\right\}  $ corresponds to the unequiary
heteromorphism with $n=4$, $n^{\prime}=2$, $k=3$, $\ell_{\operatorname*{id}%
}=2$ and $\ell_{\mu}=1$ is%
\begin{equation}
{\scriptsize
\begin{array}
[c]{c}%
\scriptsize\xymatrix@=10pt{
g_1\ar@{->}[dr] & h_1 &u _1 &g_1 & *+[F]{h_1}& u_1 \\
g_2 &h_2 \ar@/^/[rr]& u_2&g_2\ar@{->}[urr] & h_2&*+[F]{u_2}
\save"1,1"."2,3"*++[F.]\frm{}
\restore\save"1,4"."2,6"*++[F.]\frm{}
\restore}
\\[-4pt]
\xymatrix@=7pt{\\
*+[F-:<4pt>]\txt{corr}
\ar@2{->}[u]\ar@2{->}[d]
\\
\\
}%
\\
\Phi_{3}^{(4,2)}\left(
\begin{array}
[c]{c}%
\mu_{4}\left[  g_{1},h_{2},g_{2},u_{1}\right] \\
h_{1}\\
u_{2}%
\end{array}
\right)  =\mu_{2}^{\prime}\left[  \Phi_{3}^{(4,2)}\left(
\begin{array}
[c]{c}%
g_{1}\\
h_{1}\\
u_{1}%
\end{array}
\right)  ,\Phi_{3}^{(4,2)}\left(
\begin{array}
[c]{c}%
g_{2}\\
h_{2}\\
u_{2}%
\end{array}
\right)  \right]  ,
\end{array}
} \label{qv2}%
\end{equation}
where the intact elements $h_{1}$, $u_{2}$ are boxed in squares. As it is seen
from (\ref{qv2}), the product $\mu_{2}^{\prime}$ is not associative, if
$\mu_{4}$ is associative. So, not all polyadic quivers preserve associativity.

\begin{definition}
An \textit{associative polyadic quiver} is a polyadic quiver which ensures the
final associativity of $\mu_{n^{\prime}}^{\prime}$ in the main heteromorphism
equation (\ref{he}), when the initial multiplication $\mu_{n}$ is associative.
\end{definition}

One of the associative polyadic quivers which corresponds to the same
heteromorphism parameters, as the non-associative quiver (\ref{qv2}), is
$\left\{  g_{1}\rightarrow h_{2}\rightarrow u_{1}\rightarrow g_{2};h_{1}%
,u_{2}\right\}  $ which corresponds to%
\begin{equation}
{\scriptsize
\begin{array}
[c]{c}%
\scriptsize\xymatrix@=10pt{
g_1\ar@{->}[dr] &h_1 &u _1\ar@{->}[dr]  &g_1 &  *+[F]{h_1}& u_1\\
g_2 &h_2 \ar@{->}[ur]& u_2&g_2 & h_2&*+[F]{u_2}
\save"1,1"."2,3"*++[F.]\frm{}
\restore\save"1,4"."2,6"*++[F.]\frm{}
\restore}
\\[-4pt]
\xymatrix@=7pt{\\
*+[F-:<4pt>]\txt{corr}
\ar@2{->}[u]\ar@2{->}[d]
\\
\\
}%
\\[-3pt]%
\Phi_{3}^{(4,2)}\left(
\begin{array}
[c]{c}%
\mu_{4}\left[  g_{1},h_{2},u_{1},g_{2}\right] \\
h_{1}\\
u_{2}%
\end{array}
\right)  =\mu_{2}^{\prime}\left[  \Phi_{3}^{(4,2)}\left(
\begin{array}
[c]{c}%
g_{1}\\
h_{1}\\
u_{1}%
\end{array}
\right)  ,\Phi_{3}^{(4,2)}\left(
\begin{array}
[c]{c}%
g_{2}\\
h_{2}\\
u_{2}%
\end{array}
\right)  \right]  .
\end{array}
} \label{qv3}%
\end{equation}

Here we propose a classification of associative polyadic quivers and the rules
of construction of the corresponding heteromorphism equations, and then use
the heteromorphism parameters for the classification of $\ell_{\mu}$-ple
heteromorphisms (\ref{kl}). In other words, we describe a consistent procedure
for building the semigroup heteromorphisms.

Let us consider the first class of heteromorphisms (without intact elements
$\ell_{\operatorname*{id}}=0$ or \textit{intactless}), that is $\ell_{\mu}%
$-place (multiplace) homomorphisms. In the simplest case, associativity can be
achieved, when all elements in a product are taken from the same row. The
number of places $k$ is not fixed by the arity relation (\ref{n1}) and can be
arbitrary, while the arrows can have various directions. There are $2^{k}$
such combinations which preserve associativity. If the arrows have the same
direction, this kind of mapping is also called a homomorphism. As an example,
for $n=n^{\prime}=3$, $k=2$, $\ell_{\mu}=2$ we have%
\begin{equation}
{\scriptsize
\begin{array}
[c]{c}%
\scriptsize\xymatrix@=10pt{
g_1\ar@{->}[d] & h_1 \ar@{->}[d] & \\
g_2\ar@{->}[d]  &h_2\ar@{->}[d] \\
g_3  &h_3
\save"1,1"."3,2"*++[F.]\frm{}
\restore}
\\[-4pt]
\xymatrix@=7pt{\\
*+[F-:<4pt>]\txt{corr}
\ar@2{->}[u]\ar@2{->}[d]
\\
\\
}%
\\[-3pt]%
\Phi_{2}^{(3,3)}\left(
\begin{array}
[c]{c}%
\mu_{3}\left[  g_{1},g_{2},g_{3}\right] \\
\mu_{3}\left[  h_{1},h_{2},h_{3}\right]
\end{array}
\right)  =\mu_{3}^{\prime}\left[  \Phi_{2}^{(3,3)}\left(
\begin{array}
[c]{c}%
g_{1}\\
h_{1}%
\end{array}
\right)  ,\Phi_{2}^{(3,3)}\left(
\begin{array}
[c]{c}%
g_{2}\\
h_{2}%
\end{array}
\right)  ,\Phi_{2}^{(3,3)}\left(
\begin{array}
[c]{c}%
g_{3}\\
h_{3}%
\end{array}
\right)  \right]  .
\end{array}
} \label{qv4}%
\end{equation}
Note that the analogous quiver with opposite arrow directions is%
\begin{equation}
{\scriptsize
\begin{array}
[c]{c}%
\scriptsize\xymatrix@=10pt{
g_1\ar@{->}[d] & h_1 & \\
g_2\ar@{->}[d]  &h_2\ar@{->}[u] \\
g_3  &h_3 \ar@{->}[u]
\save"1,1"."3,2"*+[F.]\frm{}
\restore}
\\[-4pt]
\xymatrix@=7pt{\\
*+[F-:<4pt>]\txt{corr}
\ar@2{->}[u]\ar@2{->}[d]
\\
\\
}%
\\[-3pt]%
\Phi_{2}^{(3,3)}\left(
\begin{array}
[c]{c}%
\mu_{3}\left[  g_{1},g_{2},g_{3}\right] \\
\mu_{3}\left[  h_{3},h_{2},h_{1}\right]
\end{array}
\right)  =\mu_{3}^{\prime}\left[  \Phi_{2}^{(3,3)}\left(
\begin{array}
[c]{c}%
g_{1}\\
h_{1}%
\end{array}
\right)  ,\Phi_{2}^{(3,3)}\left(
\begin{array}
[c]{c}%
g_{2}\\
h_{2}%
\end{array}
\right)  ,\Phi_{2}^{(3,3)}\left(
\begin{array}
[c]{c}%
g_{3}\\
h_{3}%
\end{array}
\right)  \right]  ,
\end{array}
} \label{qv5}%
\end{equation}
The latter mapping and the corresponding vertical quiver were used in
constructing the middle representations of ternary groups \cite{bor/dud/dup3}.

For nonvertical quivers the main rule is the following: \emph{all arrows of an
associative quiver should have direction from left to right or vertical, and
they should not instersect}. Also, we start always from the upper left corner,
because of the permutation symmetry of (\ref{he}), we can rearrange and rename
variables in the necessary way.

An important class of intactless heteromorphisms (with $\ell
_{\operatorname*{id}}=0$) preserving associativity can be constructed using an
analogy with the Post substitutions \cite{pos}, and therefore we call it the
\textit{Post-like associative quiver}. The number of places $k$ is now fixed
by $k=n-1$, while $n^{\prime}=n$ and $\ell_{\mu}=k=n-1$. An example of the
Post-like associative quiver with the same heteromorphisms parameters as in
(\ref{qv4})-(\ref{qv5}) is%
\begin{equation}
{\scriptsize
\begin{array}
[c]{c}%
\scriptsize\xymatrix@=10pt{
g_1\ar@{->}[dr] & h_1 \ar@{->}[dr] &g_1 & h_1 \\
g_2  &h_2\ar@{->}[dr] &g_2  \ar@{->}[dr]&h_2\\
g_3  &h_3&g_3  &h_3\\
\save"1,1"."3,2"*++[F.]\frm{}
\restore\save"1,3"."3,4"*++[F.]\frm{}
\restore}
\\[-4pt]
\xymatrix@=7pt{\\
*+[F-:<4pt>]\txt{corr}
\ar@2{->}[u]\ar@2{->}[d]
\\
\\
}%
\\[-3pt]%
\Phi_{2}^{(3,3)}\left(
\begin{array}
[c]{c}%
\mu_{3}\left[  g_{1},h_{2},g_{3}\right] \\
\mu_{3}\left[  h_{1},g_{2},h_{3}\right]
\end{array}
\right)  =\mu_{3}^{\prime}\left[  \Phi_{2}^{(3,3)}\left(
\begin{array}
[c]{c}%
g_{1}\\
h_{1}%
\end{array}
\right)  ,\Phi_{2}^{(3,3)}\left(
\begin{array}
[c]{c}%
g_{2}\\
h_{2}%
\end{array}
\right)  ,\Phi_{2}^{(3,3)}\left(
\begin{array}
[c]{c}%
g_{3}\\
h_{3}%
\end{array}
\right)  \right]  .
\end{array}
} \label{qv6}%
\end{equation}
This construction appeared in the study of ternary semigroups of morphisms
\cite{chr94,chr94a,chr/nov}. Its $n$-ary generalization was used in the
consideration of polyadic operations on Cartesian powers \cite{galm08},
polyadic analogs of the Cayley and Birkhoff theorems \cite{galm01a,galmak2}
and special representations of $n$-groups \cite{gle/wan/wan,wan/wan} (where
the $n$-group with the multiplication $\mu_{2}^{\prime}$ was called the
\textit{diagonal }$n$\textit{-group}). Consider the following example.

\begin{example}
Let $\Lambda$ be the Grassmann algebra consisting of even and odd parts
$\Lambda=\Lambda_{\bar{0}}\oplus\Lambda_{\bar{1}}$ (see e.g., \cite{berezin}).
The odd part can be considered as a ternary semigroup $\mathfrak{G}$%
$_{3}^{\left(  \bar{1}\right)  }=\left\langle \Lambda_{\bar{1}},\mu
_{3}\right\rangle $, its multiplication $\mu_{3}:\Lambda_{\bar{1}}%
\times\Lambda_{\bar{1}}\times\Lambda_{\bar{1}}\rightarrow\Lambda_{\bar{1}}$ is
defined by $\mu_{3}\left[  \alpha,\beta,\gamma\right]  =\alpha\cdot\beta
\cdot\gamma$, where $\left(  \cdot\right)  $ is multiplication in $\Lambda$
and $\alpha,\beta,\gamma\in\Lambda_{\bar{1}}$, so $\mathfrak{G}$$_{3}^{\left(
\bar{1}\right)  }$ is nonderived and contains no unity. The even part can be
treated as a ternary group $\mathfrak{G}$$_{3}^{\left(  \bar{0}\right)
}=\left\langle \Lambda_{\bar{0}},\mu_{3}^{\prime}\right\rangle $ with the
multiplication $\mu_{3}^{\prime}:\Lambda_{\bar{0}}\times\Lambda_{\bar{0}%
}\times\Lambda_{\bar{0}}\rightarrow\Lambda_{\bar{0}}$, defined by $\mu
_{3}\left[  a,b,c\right]  =a\cdot b\cdot c$, where $a,b,c\in\Lambda_{\bar{0}}%
$, thus $\mathfrak{G}$$_{3}^{\left(  \bar{0}\right)  }$ is derived and
contains unity. We introduce the heteromorphism $\mathfrak{G}$$_{3}^{\left(
\bar{1}\right)  }\rightarrow\mathfrak{G}_{3}^{\left(  \bar{0}\right)  }$ as a
mapping (2-place homomorphism) $\Phi_{2}^{\left(  3,3\right)  }:\Lambda
_{\bar{1}}\times\Lambda_{\bar{1}}\rightarrow\Lambda_{\bar{0}}$ by the formula%
\begin{equation}
\Phi_{2}^{\left(  3,3\right)  }\left(
\begin{array}
[c]{c}%
\alpha\\
\beta
\end{array}
\right)  =\alpha\cdot\beta, \label{fab}%
\end{equation}
where $\alpha,\beta\in\Lambda_{\bar{1}}$. It is seen that $\Phi_{2}^{\left(
3,3\right)  }$ defined by (\ref{fab}) satisfies the Post-like heteromorphism
equation (\ref{qv6}), but not the \textquotedblleft vertical\textquotedblright%
\ one (\ref{qv4}), due to the anticommutativity of odd elements from
$\Lambda_{\bar{1}}$. In other words, $\mathfrak{G}$$_{3}^{\left(  \bar
{0}\right)  }$ can be treated as a nontrivial example of the \textquotedblleft
diagonal\textquotedblright\ semigroup of $\mathfrak{G}$$_{3}^{\left(  \bar
{1}\right)  }$ (according to the notation of \cite{gle/wan/wan,wan/wan}).
\end{example}

Note that for the number of places $k\geq3$ there exist additional (to the
above) associative quivers having the same heteromorphism parameters. For
instance, when $n^{\prime}=n=4$ and $k=3$ we have the Post-like associative
quiver%
\begin{equation}
{\scriptsize
\begin{array}
[c]{c}%
\scriptsize\xymatrix@=10pt{
g_1\ar@{->}[dr] &h_1 \ar@{->}[dr]&u _1\ar@{->}[dr]  &g_1 & h_1& u_1\\
g_2 &h_2 \ar@{->}[dr]& u_2\ar@{->}[dr]&g_2 \ar@{->}[dr]& h_2&u_2\\
g_3 &h_3 &u _3\ar@{->}[dr]  &g_3 \ar@{->}[dr]& h_3\ar@{->}[dr]& u_3\\
g_4 &h_4 & u_4&g_4 & h_4&u_4
\save"1,1"."4,3"*++[F.]\frm{}
\restore\save"1,4"."4,6"*++[F.]\frm{}
\restore}
\\[-4pt]
\xymatrix@=7pt{\\
*+[F-:<4pt>]\txt{corr}
\ar@2{->}[u]\ar@2{->}[d]
\\
\\
}%
\\[-3pt]%
\Phi_{3}^{(4,4)}\left(
\begin{array}
[c]{c}%
\mu_{4}\left[  g_{1},h_{2},u_{3},g_{4}\right] \\
\mu_{4}\left[  h_{1},u_{2},g_{3},h_{4}\right] \\
\mu_{4}\left[  u_{1},g_{2},h_{3},u_{4}\right]
\end{array}
\right)  =\mu_{4}^{\prime}\left[  \Phi_{3}^{(4,4)}\left(
\begin{array}
[c]{c}%
g_{1}\\
h_{1}\\
u_{1}%
\end{array}
\right)  ,\Phi_{3}^{(4,4)}\left(
\begin{array}
[c]{c}%
g_{2}\\
h_{2}\\
u_{2}%
\end{array}
\right)  ,\Phi_{3}^{(4,4)}\left(
\begin{array}
[c]{c}%
g_{1}\\
h_{1}\\
u_{1}%
\end{array}
\right)  ,\Phi_{3}^{(4,4)}\left(
\begin{array}
[c]{c}%
g_{2}\\
h_{2}\\
u_{2}%
\end{array}
\right)  \right]  .
\end{array}
} \label{qv7}%
\end{equation}
Also, we have one intermediate non-Post associative quiver%
\begin{equation}
{\scriptsize
\begin{array}
[c]{c}%
\scriptsize\xymatrix@=10pt{
g_1\ar@{->}[drr] &h_1 \ar@{->}[drr]&u _1\ar@{->}%
[drr]  &g_1 & h_1& u_1&g_1 & h_1& u_1\\
g_2 &h_2 & u_2\ar@{->}[drr]&g_2 \ar@{->}[drr]& h_2\ar@{->}%
[drr]&u_2&g_2 & h_2& u_2\\
g_3 &h_3 &u _3 &g_3& h_3\ar@{->}[drr]& u_3\ar@{->}[drr]&g_3 \ar@
{->}[drr]& h_3& u_3\\
g_4 &h_4 & u_4&g_4 & h_4&u_4&g_4 & h_4& u_4
\save"1,1"."4,3"*++[F.]\frm{}
\restore\save"1,4"."4,6"*++[F.]\frm{}
\restore\save"1,7"."4,9"*++[F.]\frm{}
\restore}
\\[-4pt]
\xymatrix@=7pt{\\
*+[F-:<4pt>]\txt{corr}
\ar@2{->}[u]\ar@2{->}[d]
\\
\\
}%
\\[-3pt]%
\Phi_{3}^{(4,4)}\left(
\begin{array}
[c]{c}%
\mu_{4}\left[  g_{1},u_{2},h_{3},g_{4}\right] \\
\mu_{4}\left[  h_{1},g_{2},u_{3},h_{4}\right] \\
\mu_{4}\left[  u_{1},h_{2},g_{3},u_{4}\right]
\end{array}
\right)  =\mu_{4}^{\prime}\left[  \Phi_{3}^{(4,4)}\left(
\begin{array}
[c]{c}%
g_{1}\\
h_{1}\\
u_{1}%
\end{array}
\right)  ,\Phi_{3}^{(4,4)}\left(
\begin{array}
[c]{c}%
g_{2}\\
h_{2}\\
u_{2}%
\end{array}
\right)  ,\Phi_{3}^{(4,4)}\left(
\begin{array}
[c]{c}%
g_{1}\\
h_{1}\\
u_{1}%
\end{array}
\right)  ,\Phi_{3}^{(4,4)}\left(
\begin{array}
[c]{c}%
g_{2}\\
h_{2}\\
u_{2}%
\end{array}
\right)  \right]  .
\end{array}
} \label{qv8}%
\end{equation}
A general method of constructing associative quivers for $n^{\prime}=n$,
$\ell_{\operatorname*{id}}=0$ and $k=n-1$ can be illustrated from the
following more complicated example with $n=5$. First, we draw $k=n-1$ ($=4$)
copies of element matrices. Then we go from the first element in the first
column $g_{1}$ to the last element in this column $g_{n^{\prime}}$ ($=g_{5}$)
by different $k=n-1$ ($=4$) ways: by the vertical quiver, by the Post-like
quiver (going to the second copy of the element matrix) and by the remaining
$k-2$ ($=2$) non-Post associative quivers, as%
\begin{equation}%
\begin{scriptsize}\xymatrix@=10pt{
g_1\ar@{->}[d] \ar@{->}[dr]\ar@{->}[drr]\ar@{->}%
[drrr]&h_1 &u _1&v_1  &g_1 & h_1& u_1&v_1&g_1 & h_1& u_1&v_1&g_1 & h_1& u_1&v_1\\
g_2\ar@{->}[d] &h_2 \ar@{->}[dr]& u_2\ar@{->}[drr]&v_2\ar@{->}%
[drrr]&g_2 & h_2&u_2&v_2&g_2 & h_2& u_2&v_2&g_2 & h_2& u_2&v_2\\
g_3\ar@{->}[d] &h_3 &u _3 \ar@{->}[dr]&v_3&g_3\ar@{->}[drr]& h_3& u_3\ar@
{->}[drrr]&v_3&g_3 & h_3& u_3&v_3&g_3 & h_3& u_3&v_3\\
g_4 \ar@{->}[d]&h_4 & u_4&v_4\ar@{->}[dr]&g_4 & h_4&u_4\ar@{->}%
[drr]&v_4&g_4 & h_4\ar@{->}[drrr]& u_4&v_4&g_4 & h_4& u_4&v_4\\
g_5 &h_5 & u_5&v_5&g_5 & h_5&u_5&v_5&g_5 & h_5& u_5&v_5&g_5 & h_5& u_5&v_5\\
\txt{vertical}&&&&\txt{Post}&&&&\txt{non-Post}&&&&\txt{non-Post}
\save"1,1"."5,4"*++[F.]\frm{}
\restore\save"1,5"."5,8"*++[F.]\frm{}
\restore\save"1,9"."5,12"*++[F.]\frm{}
\restore\save"1,13"."5,16"*++[F.]\frm{}
\restore}
\end{scriptsize}%
. \label{qv9}%
\end{equation}
Here we show, for short, only the quiver itself without the corresponding
heteromorphism equation and only the arrows corresponding to the first
product, while other arrows (starting from $h_{1}$, $u_{1}$, $v_{1}$) are
parallel to it, as in (\ref{qv8}).

The next type of heteromorphisms (intermediate) is described by the equations
(\ref{k1})-(\ref{ll}), it contains intact elements ($\ell_{\operatorname*{id}%
}\geq1$) and changes (decreases) arity $n^{\prime}<n$. For each fixed $k$ the
arities are not arbitrary and presented in \textsc{Table \ref{T1}}. The first
general rule is: the associative quivers are nondecreasing in both,vertical
(from up to down) and horizontal (from left to right), directions. Second, if
there are several multiplications ($\ell_{\mu}\geq2$), the corresponding
associative quivers do not intersect.

Let us present some examples and start from the smallest number of
heteromorphism places in $\Phi_{k}^{\left(  n,n^{\prime}\right)  }$. For
$k=2$, the first (nonbinarizing $n^{\prime}\geq3$) case is $n=5$, $n^{\prime
}=3$, $\ell_{\operatorname*{id}}=1$ (see the first row and second
$n/n^{\prime}$ pair of \textsc{Table \ref{T1}}). The corresponding associative
quivers are%
\begin{equation}
{\scriptsize
\begin{array}
[c]{c}%
\scriptsize\xymatrix@=10pt{
g_1\ar@{->}[r] & h_1 \ar@{->}[dr] & g_1 & h_1 & g_1 &h_1\\
g_2  &h_2&g_2\ar@{->}[r]  &h_2\ar@{->}[dr]&g_2  &h_2 \\
g_3  &h_3&g_3  &h_3&g_3  &*+[F]{h_3}
\save"1,1"."3,2"*++[F.]\frm{}
\restore\save"1,3"."3,4"*++[F.]\frm{}
\restore\save"1,5"."3,6"*++[F.]\frm{}
\restore}
\\[-4pt]
\xymatrix@=7pt{\\
*+[F-:<4pt>]\txt{corr}
\ar@2{->}[u]\ar@2{->}[d]
\\
\\
}%
\\[-3pt]%
\Phi_{2}^{(5,3)}\left(
\begin{array}
[c]{c}%
\mu_{5}\left[  g_{1},h_{1,}g_{2},h_{2,}g_{3}\right] \\
h_{3}%
\end{array}
\right)
\end{array}
\ \ \
\begin{array}
[c]{c}%
\scriptsize\xymatrix@=10pt{
g_1\ar@{->}[dr] & h_1 & g_1 & h_1 & g_1 & *+[F]{h_1}\\
g_2  &h_2\ar@{->}[r]&g_2\ar@{->}[dr]  &h_2&g_2  &h_2 \\
g_3  &h_3&g_3  &h_3\ar@{->}[r]&g_3  &h_3
\save"1,1"."3,2"*++[F.]\frm{}
\restore\save"1,3"."3,4"*++[F.]\frm{}
\restore\save"1,5"."3,6"*++[F.]\frm{}
\restore}
\\[-4pt]
\xymatrix@=7pt{\\
*+[F-:<4pt>]\txt{corr}
\ar@2{->}[u]\ar@2{->}[d]
\\
\\
}%
\\[-3pt]%
\Phi_{2}^{(5,3)}\left(
\begin{array}
[c]{c}%
\mu_{5}\left[  g_{1},h_{2,}g_{2},h_{3,}g_{3}\right] \\
h_{1}%
\end{array}
\right)  ,
\end{array}
} \label{qv10}%
\end{equation}
where we do not write the r.h.s. of the heteromorphism equation (\ref{he}),
because it is simply related to the transposed quiver matrix of the size
$n^{\prime}\times k$.

More complicated examples can be given for $k=3$, which corresponds to the
second and third lines of \textsc{Table \ref{T1}}. That is we can obtain a
ternary final product ($n^{\prime}=3$) by using one or two multiplications
$\ell_{\mu}=1,2$. Examples of the corresponding associative quivers are
($\ell_{\mu}=1$)%
\begin{equation}
{\scriptsize
\begin{array}
[c]{c}%
\scriptsize\xymatrix@=10pt{
g_1\ar@{->}[r] & h_1\ar@{->}[dr] &u_1 & g_1 & h_1 & u_1&g_1 &h_1&*+[F]{u_1}\\
g_2  &h_2&u_2\ar@{->}[r] &g_2\ar@{->}[dr]  &h_2&u_2&g_2  &*+[F]{h_2}&u_2 \\
g_3  &h_3&u_3&g_3  &h_3\ar@{->}[r] &u_3\ar@{->}[r] &g_3  &h_3&u_3
\save"1,1"."3,3"*++[F.]\frm{}
\restore\save"1,4"."3,6"*++[F.]\frm{}
\restore\save"1,7"."3,9"*++[F.]\frm{}
\restore}
\\[-4pt]
\xymatrix@=7pt{\\
*+[F-:<4pt>]\txt{corr}
\ar@2{->}[u]\ar@2{->}[d]
\\
\\
}%
\\[-3pt]%
\Phi_{3}^{(7,3)}\left(
\begin{array}
[c]{c}%
\mu_{7}\left[  g_{1},h_{1,}u_{2},g_{2},h_{3},u_{3},g_{3}\right] \\
h_{2}\\
u_{1}%
\end{array}
\right)
\end{array}
} \label{qv12}%
\end{equation}
and ($\ell_{\mu}=2$)%
\begin{equation}
{\scriptsize
\begin{array}
[c]{c}%
\scriptsize\xymatrix@=10pt{
g_1\ar@{->}[drr] & h_1\ar@/_/[rrrr] &u_1 & g_1 & h_1 & u_1\ar@{->}%
[dr] &g_1 &h_1&u_1\\
g_2  &h_2&u_2\ar@/_/[rr] &g_2 &h_2\ar@{->}[drr] &u_2&g_2 \ar@{->}%
[dr]  &h_2&u_2 \\
g_3  &h_3&u_3&g_3  &h_3 &u_3 &g_3  &h_3&*+[F]{u_3}
\save"1,1"."3,3"*++[F.]\frm{}
\restore\save"1,4"."3,6"*++[F.]\frm{}
\restore\save"1,7"."3,9"*++[F.]\frm{}
\restore}
\\[-4pt]
\xymatrix@=7pt{\\
*+[F-:<4pt>]\txt{corr}
\ar@2{->}[u]\ar@2{->}[d]
\\
\\
}%
\\[-3pt]%
\Phi_{3}^{(4,3)}\left(
\begin{array}
[c]{c}%
\mu_{4}\left[  g_{1},u_{2},h_{2},g_{3}\right] \\
\mu_{4}\left[  h_{1,}u_{1},g_{2},h_{3}\right] \\
u_{3}%
\end{array}
\right)
\end{array}
} \label{qv13}%
\end{equation}
respectively. Finally, for the case $k=4$ one can construct the associative
quiver corresponding to the first pair of the last line in \textsc{Table
\ref{T1}}. It has three multiplications and one intact element, and the
corresponding quiver is, e.g.,%
\begin{equation}
{\scriptsize
\begin{array}
[c]{c}%
\scriptsize\xymatrix@=10pt{
g_1\ar@{->}[dr] & h_1\ar@/_/[drr] &u_1&v_1\ar@/_/[rrr]  & g_1 & h_1 & u_1\ar@
{->}[drr] &v_1&g_1 &h_1&u_1&v_1\\
g_2  &h_2\ar@{->}[dr] &u_2 &v_2\ar@/_/[rrr]&g_2 &h_2 &u_2\ar@{->}%
[dr]&v_2&g_2 \ar@{->}[dr]  &h_2&u_2&v_2 \\
g_3  &h_3&u_3\ar@/_/[rr] &v_3&g_3\ar@{->}[drrrr]   &h_3 &u_3&v_3\ar@
{->}[drr]  &g_3  &h_3\ar@{->}[drr] &u_3&v_3\\
g_4  &h_4&u_4&v_4&g_4  &h_4 &u_4&v_4 &g_4  &h_4&*+[F]{u_4}&v_4
\save"1,1"."4,4"*++[F.]\frm{}
\restore\save"1,5"."4,8"*++[F.]\frm{}
\restore\save"1,9"."4,12"*++[F.]\frm{}
\restore}
\\[-4pt]
\xymatrix@=7pt{\\
*+[F-:<4pt>]\txt{corr}
\ar@2{->}[u]\ar@2{->}[d]
\\
\\
}%
\\[-3pt]%
\Phi_{4}^{(5,4)}\left(
\begin{array}
[c]{c}%
\mu_{5}\left[  g_{1},h_{2},u_{3},g_{3},g_{4}\right] \\
\mu_{5}\left[  h_{1},v_{2},u_{2},v_{3},h_{4}\right] \\
\mu_{5}\left[  v_{1},u_{1},g_{2},h_{3},v_{4}\right] \\
u_{4}%
\end{array}
\right)  .
\end{array}
} \label{qv14}%
\end{equation}

There are many other possibilities (using permutations and different variants
of quivers) to obtain an associative final product $\mu_{n^{\prime}}^{\prime}$
corresponding the same heteromorphism parameters, and therefore we do list
them all. The above examples are sufficient to understand the rules of the
associative quiver construction and obtain the polyadic semigroup heteromorphisms.

\section{Multiplace representations of polyadic systems}

Representation theory \cite{kirillov} deals with mappings from abstract
algebraic systems into linear systems, such as, e.g. linear operators in
vector spaces, or into general (semi)groups of transformations of some set. In
our notation, this means that in the mapping of polyadic systems (\ref{gg1})
the final multiplication $\mu_{n^{\prime}}^{\prime}$ is a linear map. This
leads to some restrictions on the final polyadic structure $\mathfrak{G}%
$$_{n^{\prime}}^{\prime}$, which are considered below.

Let $V$ be a vector space over a field $\mathbb{K}$ (usually algebraically
closed) and $\operatorname*{End}V$ be a set of linear endomorphisms of $V$,
which is in fact a binary group. In the standard way, a linear representation
of a binary semigroup $\mathfrak{G}$$_{2}=\left\langle G;\mu_{2}\right\rangle
$ is a (1-place) map $\Pi_{1}:\mathfrak{G}_{2}\rightarrow\operatorname*{End}%
V$, such that $\Pi_{1}$ is a homomorphism%
\begin{equation}
\Pi_{1}\left(  \mu_{2}\left[  g,h\right]  \right)  =\Pi_{1}\left(  g\right)
\ast\Pi_{1}\left(  h\right)  , \label{ppp}%
\end{equation}
where $g,h\in G$ and $\left(  \ast\right)  $ is the binary multiplication in
$\operatorname*{End}V$ (usually, it is a (semi)group with multiplication as
composition of operators or product of matrices, if a basis is chosen). If
$\mathfrak{G}$$_{2}$ is a binary group with the unity $e$, then we have the
additional condition%
\begin{equation}
\Pi_{1}\left(  e\right)  =\operatorname*{id}\nolimits_{V}. \label{pid}%
\end{equation}

We will generalize these known formulas to the corresponding polyadic systems
along with the heteromorphism concept introduced above. Our general idea is to
use the heteromorphism equation (\ref{he}) instead of the standard
homomorphism equation (\ref{ppp}), such that the arity of the representation
will be different from the arity of the initial polyadic system $n^{\prime
}\neq n$.

Consider the structure of the final $n^{\prime}$-ary multiplication
$\mu_{n^{\prime}}^{\prime}$ in (\ref{he}), taking into account that the final
polyadic system $\mathfrak{G}$$_{n^{\prime}}^{\prime}$ should be constructed
from $\operatorname*{End}V$. The most natural and physically applicable way is
to consider the binary $\operatorname*{End}V$ and to put $\mathfrak{G}%
$$_{n^{\prime}}^{\prime}=\operatorname*{der}\nolimits_{n^{\prime}}\left(
\operatorname*{End}V\right)  $, as it was proposed for the ternary case in
\cite{bor/dud/dup3}. In this way $\mathfrak{G}$$_{n^{\prime}}^{\prime}$
becomes a derived $n^{\prime}$-ary (semi)group of endomorphisms of $V$ with
the multiplication $\mu_{n^{\prime}}^{\prime}:$ $\left(  \operatorname*{End}%
V\right)  ^{\times n^{\prime}}\rightarrow\operatorname*{End}V$, where%
\begin{equation}
\mu_{n^{\prime}}^{\prime}\left[  v_{1},\ldots,v_{n^{\prime}}\right]
=v_{1}\ast\ldots\ast v_{n^{\prime}},\ \ \ v_{i}\in\operatorname*{End}V.
\label{mv}%
\end{equation}

Because the multiplication $\mu_{n^{\prime}}^{\prime}$ (\ref{mv}) is derived
and is therefore associative by definition, we may consider the associative
initial polyadic systems (semigroups and groups) and the associativity
preserving mappings that are the special heteromorphisms constructed in the
previous section.

Let $\mathfrak{G}$$_{n}=\left\langle G;\mu_{n}\right\rangle $ be an
associative $n$-ary polyadic system. By analogy with (\ref{fgg}), we introduce
the following $k$-place mapping%
\begin{equation}
\Pi_{k}^{\left(  n,n^{\prime}\right)  }:G^{\times k}\rightarrow
\operatorname*{End}V. \label{pe}%
\end{equation}

\textit{A multiplace representation} of an associative polyadic system
$\mathfrak{G}_{n}$ in a vector space $V$ is given, if there exists a $k$-place
mapping (\ref{pe}) which satisfies the (associativity preserving)
heteromorphism equation (\ref{he}), that is
\begin{equation}
\Pi_{k}^{\left(  n,n^{\prime}\right)  }\left(
\genfrac{}{}{0pt}{}{\left.
\begin{array}
[c]{c}%
\mu_{n}\left[  g_{1},\ldots,g_{n}\right]  ,\\
\vdots\\
\mu_{n}\left[  g_{n\left(  \ell_{\mu}-1\right)  },\ldots,g_{n\ell_{\mu}%
}\right]
\end{array}
\right\}  \ell_{\mu}}{\left.
\begin{array}
[c]{c}%
g_{n\ell_{\mu}+1},\\
\vdots\\
g_{n\ell_{\mu}+\ell_{\operatorname*{id}}}%
\end{array}
\right\}  \ell_{\operatorname*{id}}}%
\right)  =\overset{n^{\prime}}{\overbrace{\Pi_{k}^{\left(  n,n^{\prime
}\right)  }\left(
\begin{array}
[c]{c}%
g_{1}\\
\vdots\\
g_{k}%
\end{array}
\right)  \ast\ldots\ast\Pi_{k}^{\left(  n,n^{\prime}\right)  }\left(
\begin{array}
[c]{c}%
g_{k\left(  n^{\prime}-1\right)  }\\
\vdots\\
g_{kn^{\prime}}%
\end{array}
\right)  }}, \label{ppe}%
\end{equation}
{and the following diagram commutes}%
\begin{equation}
\begin{diagram} G^{\times k} & \rTo^{\Pi_{k}} & \operatorname*{End}V \\ \uTo^{\mu_{n}^{\left( \ell_{\mu},\ell_{\operatorname*{id}}\right) }} & & \uTo_{(\ast)^{n^{\prime}}} \\ G^{\times kn^{\prime}} & \rTo^{ \left( \Pi_{k}\right) ^{\times n^{\prime}}} & \left( \operatorname*{End}V\right)^{\times n^{\prime}} \\ \end{diagram} \label{dia5}%
\end{equation}
where $\mu_{n}^{\left(  \ell_{\mu},\ell_{\operatorname*{id}}\right)  }$ is
given by (\ref{midl}), $\ell_{\mu}$ and $\ell_{\operatorname*{id}}$ are the
numbers of multiplications and intact elements in the l.h.s. of (\ref{ppe}), respectively.

The exact permutation in the l.h.s. of (\ref{ppe}) is given by the associative
quiver presented in the previous section. The representation parameters ($n$,
$n^{\prime}$, $k$, $\ell_{\mu}$\textit{ and }$\ell_{\operatorname*{id}}$) in
(\ref{ppe}) are the same as the heteromorphism parameters, and they satisfy
the same arity changing formulas (\ref{n1}) and (\ref{n2}). Therefore, a
general classification of multiplace representations can be done by analogy
with that of the heteromorphisms (\ref{kmin})--(\ref{kmax}) as follows:

\begin{enumerate}
\item The \textit{hom-like multiplace representation} which is a multiplace
homomorphism with $n^{\prime}=n_{\max}^{\prime}=n$\textit{,} without intact
elements $l_{\operatorname*{id}}=l_{\operatorname*{id}}^{(\min)}=0$, and
minimal number of places%
\begin{equation}
k=k_{\min}=\ell_{\mu}.
\end{equation}

\item The \textit{intact element multiplace representation} which is the
intermediate heteromorphism with $2<n^{\prime}<n$ and the number of intact
elements is%
\begin{equation}
l_{\operatorname*{id}}=\dfrac{n-n^{\prime}}{n^{\prime}-1}\ell_{\mu}.
\end{equation}

\item The \textit{binary multiplace representation} which is a binarizing
heteromorphism (\ref{kmax}) with $n^{\prime}=n_{\min}^{\prime}=2$, the maximal
number of intact elements $l_{\operatorname*{id}}^{(\max)}=\left(  n-2\right)
\ell_{\mu}$ and maximal number of places%
\begin{equation}
k=k_{\max}=\left(  n-1\right)  \ell_{\mu}.
\end{equation}

\end{enumerate}

The multiplace representations for $n$-ary semigroups have no additional
defining relations, as compared with (\ref{ppe}). In case of $n$-ary groups,
we need an analog of the \textquotedblleft normalizing\textquotedblright%
\ relation (\ref{pid}). If the $n$-ary group has the unity $e$, then one can
put%
\begin{equation}
\Pi_{k}^{\left(  n,n^{\prime}\right)  }\left(  \left.
\begin{array}
[c]{c}%
e\\
\vdots\\
e
\end{array}
\right\}  k\right)  =\operatorname*{id}\nolimits_{V}. \label{pke}%
\end{equation}
If there is no unity at all, one can \textquotedblleft
normalize\textquotedblright\ the multiplace representation, using analogy with
(\ref{pid}) in the form%
\begin{equation}
\Pi_{1}\left(  h^{-1}\ast h\right)  =\operatorname*{id}\nolimits_{V},
\end{equation}
as follows%
\begin{equation}
\Pi_{k}^{\left(  n,n^{\prime}\right)  }\left(
\begin{array}
[c]{l}%
\left.
\begin{array}
[c]{c}%
\bar{h}\\
\vdots\\
\bar{h}%
\end{array}
\right\}  \ell_{\mu}\\[20pt]%
\left.
\begin{array}
[c]{c}%
h\\
\vdots\\
h
\end{array}
\right\}  \ell_{\operatorname*{id}}%
\end{array}
\right)  =\operatorname*{id}\nolimits_{V}, \label{pkg}%
\end{equation}
for all $h\in$$\mathfrak{G}$$_{n}$, where $\bar{h}$ is the querelement of $h$.
The latter ones can be placed on any places in the l.h.s. of (\ref{pkg}) due
to the D\"{o}rnte identities. Also, the multiplications in the l.h.s. of
(\ref{ppe}) can change their place due to the same reason.

A general form of multiplace representations can be found by applying the
D\"{o}rnte identities to each $n$-ary product in the l.h.s. of (\ref{ppe}).
Then, using (\ref{pkg}) we have schematically%
\begin{equation}
\Pi_{k}^{\left(  n,n^{\prime}\right)  }\left(
\begin{array}
[c]{c}%
g_{1}\\
\vdots\\
g_{k}%
\end{array}
\right)  =\Pi_{k}^{\left(  n,n^{\prime}\right)  }\left(
\begin{array}
[c]{l}%
\begin{array}
[c]{c}%
t_{1}\\
\vdots\\
t_{\ell_{\mu}}%
\end{array}
\\[20pt]%
\left.
\begin{array}
[c]{c}%
g\\
\vdots\\
g
\end{array}
\right\}  \ell_{\operatorname*{id}}%
\end{array}
\right)  , \label{ppu}%
\end{equation}
where $g$ is an arbitrary fixed element of the $n$-ary group and%
\begin{equation}
t_{a}=\mu_{n}\left[  g_{a1},\ldots,g_{an-1},\bar{g}\right]  ,\ \ \ a=1,\ldots
,\ell_{\mu}. \label{ua}%
\end{equation}

This is the special shape of some multiplace representations, while the
concrete formulas should be obtained in each case separately. Nevertheless,
some conclusions can be drawn from (\ref{ppu}). Firstly, the equivalence
classes on which $\Pi_{k}^{\left(  n,n^{\prime}\right)  }$ is constant are
determined by fixing $\ell_{\mu}+1$ elements, i.e. by the surface
$t_{a}=const$, $g=const$. Secondly, some $k$-place representations of a
$n$-ary group can be reduced to $\ell_{\mu}$-place representations of its
retract. In the case $\ell_{\mu}=1$, multiplace representations of a $n$-ary
group derived from a binary group correspond to ordinary representations of
the latter (see \cite{bor/dud/dup3,dud07}).

\begin{example}
Let us consider the case of a binary multiplace representation $\Pi
_{2n-2}^{\left(  n,2\right)  }$ of $n$-ary group $\mathfrak{G}$$_{n}%
=\left\langle G;\mu_{n}\right\rangle $ with two multiplications $\ell_{\mu}=2$
defined by the associativity preserving equation{\scriptsize
\begin{equation}
\Pi_{2n-2}^{\left(  n,2\right)  }\left(
\begin{array}
[c]{c}%
\mu_{n}\left[  g_{1},u_{1},\ldots u_{n-2},g_{2}\right] \\
u_{1}^{\prime}\\
\vdots\\
u_{n-2}^{\prime}\\
\mu_{n}\left[  h_{1},v_{1},\ldots v_{n-2},h_{2}\right] \\
v_{1}^{\prime}\\
\vdots\\
v_{n-2}^{\prime}%
\end{array}
\right)  =\Pi_{2n-2}^{\left(  n,2\right)  }\left(
\begin{array}
[c]{c}%
g_{1}\\
u_{1}\\
\vdots\\
u_{n-2}\\
h_{1}\\
v_{1}\\
\vdots\\
v_{n-2}%
\end{array}
\right)  \ast\Pi_{2n-2}^{\left(  n,2\right)  }\left(
\begin{array}
[c]{c}%
g_{2}\\
u_{1}^{\prime}\\
\vdots\\
u_{n-2}^{\prime}\\
h_{2}\\
v_{1}^{\prime}\\
\vdots\\
v_{n-2}^{\prime}%
\end{array}
\right)  \label{ppn}%
\end{equation}
} and normalizing condition{\scriptsize
\begin{equation}
\Pi_{2n-2}^{\left(  n,2\right)  }\left(
\begin{array}
[c]{l}%
\left.
\begin{array}
[c]{c}%
h\\
\vdots\\
h
\end{array}
\right\}  n-2\\
\ \ \ \bar{h}\\[2pt]%
\left.
\begin{array}
[c]{c}%
h\\
\vdots\\
h
\end{array}
\right\}  n-2\\
\ \ \ \bar{h}%
\end{array}
\right)  =\operatorname*{id}\nolimits_{V}, \label{p2}%
\end{equation}
} where $h\in$$\mathfrak{G}$$_{n}$ is arbitrary. Using (\ref{ppn}) and
(\ref{p2}) for a general form of this $\left(  2n-2\right)  $-place
representation we have{\scriptsize
\begin{equation}
\Pi_{2n-2}^{\left(  n,2\right)  }=\Pi_{2n-2}^{\left(  n,2\right)  }\left(
\begin{array}
[c]{c}%
g_{1}\\
u_{1}\\
\vdots\\
u_{n-2}\\
h_{1}\\
v_{1}\\
\vdots\\
v_{n-2}%
\end{array}
\right)  =\Pi_{2n-2}^{\left(  n,2\right)  }\left(
\begin{array}
[c]{c}%
g_{1}\\
u_{1}\\
\vdots\\
u_{n-2}\\
h_{1}\\
v_{1}\\
\vdots\\
v_{n-2}%
\end{array}
\right)  \ast\Pi_{2n-2}^{\left(  n,2\right)  }\left(
\begin{array}
[c]{l}%
\left.
\begin{array}
[c]{c}%
h\\
\vdots\\
h
\end{array}
\right\}  n-2\\
\ \ \ \bar{h}\\[2pt]%
\left.
\begin{array}
[c]{c}%
h\\
\vdots\\
h
\end{array}
\right\}  n-2\\
\ \ \ \bar{h}%
\end{array}
\right)  =\Pi_{2n-2}^{\left(  n,2\right)  }\left(
\begin{array}
[c]{l}%
\mu_{n}\left[  g_{1},u_{1},\ldots u_{n-2},h\right] \\
\left.
\begin{array}
[c]{c}%
h\\
\vdots\\
h
\end{array}
\right\}  n-3\\
\ \ \ \bar{h}\\
\mu_{n}\left[  h_{1},v_{1},\ldots v_{n-2},h\right] \\
\left.
\begin{array}
[c]{c}%
h\\
\vdots\\
h
\end{array}
\right\}  n-3\\
\ \ \ \bar{h}%
\end{array}
\right)  . \label{ppm}%
\end{equation}
}

The D\"{o}rnte identity applied to the first elements of the products in the
r.h.s. of (\ref{ppm}) together with associativity of $\mu_{n}$
gives{\scriptsize
\begin{equation}
\Pi_{2n-2}^{\left(  n,2\right)  }=\Pi_{2n-2}^{\left(  n,2\right)  }\left(
\begin{array}
[c]{l}%
\mu_{n}\left[  \overset{n-2}{\overbrace{g,\ldots,g}},t_{g},h\right] \\
\left.
\begin{array}
[c]{c}%
h\\
\vdots\\
h
\end{array}
\right\}  n-3\\
\ \ \ \bar{h}\\
\mu_{n}\left[  \overset{n-2}{\overbrace{g,\ldots,g}},t_{h},h\right] \\
\left.
\begin{array}
[c]{c}%
h\\
\vdots\\
h
\end{array}
\right\}  n-3\\
\ \ \ \bar{h}%
\end{array}
\right)  =\Pi_{2n-2}^{\left(  n,2\right)  }\left(
\begin{array}
[c]{l}%
\left.
\begin{array}
[c]{c}%
g\\
\vdots\\
g
\end{array}
\right\}  n-2\\
\ \ \ t_{g}\\
\left.
\begin{array}
[c]{c}%
g\\
\vdots\\
g
\end{array}
\right\}  n-2\\
\ \ \ t_{h}%
\end{array}
\right)  \ast\Pi_{2n-2}^{\left(  n,2\right)  }\left(
\begin{array}
[c]{l}%
\left.
\begin{array}
[c]{c}%
h\\
\vdots\\
h
\end{array}
\right\}  n-2\\
\ \ \ \bar{h}\\[2pt]%
\left.
\begin{array}
[c]{c}%
h\\
\vdots\\
h
\end{array}
\right\}  n-2\\
\ \ \ \bar{h}%
\end{array}
\right)  =\Pi_{2n-2}^{\left(  n,2\right)  }\left(
\begin{array}
[c]{l}%
\left.
\begin{array}
[c]{c}%
g\\
\vdots\\
g
\end{array}
\right\}  n-2\\
\ \ \ t_{g}\\
\left.
\begin{array}
[c]{c}%
g\\
\vdots\\
g
\end{array}
\right\}  n-2\\
\ \ \ t_{h}%
\end{array}
\right)  , \label{pp2}%
\end{equation}
} where%
\begin{equation}
t_{g}=\mu_{n}\left[  \bar{g},g_{1},u_{1},\ldots u_{n-2}\right]  ,\ \ \ t_{h}%
=\mu_{n}\left[  \bar{g},h_{1},v_{1},\ldots v_{n-2}\right]  .
\end{equation}

Thus, the equivalent classes of the multiplace representation $\Pi
_{2n-2}^{\left(  n,2\right)  }$ (\ref{ppn}) are determined by the $\left(
\ell_{\mu}+1=3\right)  $-element surface%
\begin{equation}
t_{g}=const,\ t_{h}=const,\ g=const.
\end{equation}

Note that the \textquotedblleft$\ell_{\mu}$-place reduction\textquotedblright%
\ of a multiplace representation is possible not for all associativity
preserving heteromorphism equations. For instance, if to exchange
$v_{i}\leftrightarrow v_{i}^{\prime}$ in l.h.s. of (\ref{ppn}), then the
associativity remains, but the \textquotedblleft$\ell_{\mu}$-place
reduction\textquotedblright, analogous to (\ref{pp2}) will not be possible.

In case, when it is possible, the corresponding $\ell_{\mu}$-place
representation can be realized on the binary retract of $\mathfrak{G}$$_{n}$
(for some special $n$-ary groups and $\ell_{\mu}=1$ see see
\cite{bor/dud/dup3,dud07}). Indeed, let $\mathfrak{G}$$_{2}^{ret}=\left\langle
G,\circledast\right\rangle =\operatorname*{ret}\nolimits_{g}\left\langle
G;\mu_{n}\right\rangle $, where $g_{1}\circledast g_{2}\overset{def}{=}\mu
_{n}\left[  g_{1},\overset{n-2}{\overbrace{g,\ldots,g}},g_{2}\right]  $,
$g_{1},g_{2},g\in G$ and we define the reduced $\ell_{\mu}$-place
representation (now $\ell_{\mu}=2$) through (\ref{pp2}) as follows
{\scriptsize
\begin{equation}
\Pi_{2}^{\operatorname*{ret}g}\left(
\begin{array}
[c]{c}%
g_{1}\\
g_{2}%
\end{array}
\right)  \overset{def}{=}\Pi_{2n-2}^{\left(  n,2\right)  }\left(
\begin{array}
[c]{l}%
\left.
\begin{array}
[c]{c}%
g\\
\vdots\\
g
\end{array}
\right\}  n-2\\
\ \ \ g_{1}\\
\left.
\begin{array}
[c]{c}%
g\\
\vdots\\
g
\end{array}
\right\}  n-2\\
\ \ \ g_{2}%
\end{array}
\right)  .
\end{equation}
} From (\ref{ppm}) and (\ref{p2}) we obtain{\scriptsize
\begin{align*}
&  \Pi_{2}^{\operatorname*{ret}g}\left(
\begin{array}
[c]{c}%
g_{1}\\
g_{2}%
\end{array}
\right)  \ast\Pi_{2}^{\operatorname*{ret}g}\left(
\begin{array}
[c]{c}%
g_{1}^{\prime}\\
g_{2}^{\prime}%
\end{array}
\right)  =\Pi_{2n-2}^{\left(  n,2\right)  }\left(
\begin{array}
[c]{l}%
\left.
\begin{array}
[c]{c}%
g\\
\vdots\\
g
\end{array}
\right\}  n-2\\
\ \ \ g_{1}\\
\left.
\begin{array}
[c]{c}%
g\\
\vdots\\
g
\end{array}
\right\}  n-2\\
\ \ \ g_{2}%
\end{array}
\right)  \ast\Pi_{2n-2}^{\left(  n,2\right)  }\left(
\begin{array}
[c]{l}%
\left.
\begin{array}
[c]{c}%
g\\
\vdots\\
g
\end{array}
\right\}  n-2\\
\ \ \ g_{1}^{\prime}\\
\left.
\begin{array}
[c]{c}%
g\\
\vdots\\
g
\end{array}
\right\}  n-2\\
\ \ \ g_{2}^{\prime}%
\end{array}
\right)  =\Pi_{2n-2}^{\left(  n,2\right)  }\left(
\begin{array}
[c]{l}%
\mu_{n}\left[  \overset{n-2}{\overbrace{g,\ldots,g}},g_{1},g\right] \\
\left.
\begin{array}
[c]{c}%
g\\
\vdots\\
g
\end{array}
\right\}  n-3\\
\ \ \ g_{1}^{\prime}\\
\mu_{n}\left[  \overset{n-2}{\overbrace{g,\ldots,g}},g_{2},g\right] \\
\left.
\begin{array}
[c]{c}%
g\\
\vdots\\
g
\end{array}
\right\}  n-3\\
\ \ \ g_{2}^{\prime}%
\end{array}
\right) \\
&  =\Pi_{2n-2}^{\left(  n,2\right)  }\left(
\begin{array}
[c]{l}%
\left.
\begin{array}
[c]{c}%
g\\
\vdots\\
g
\end{array}
\right\}  n-2\\
\ \ \ \mu_{n}\left[  g_{1},\overset{n-2}{\overbrace{g,\ldots,g}},g_{1}%
^{\prime}\right] \\
\left.
\begin{array}
[c]{c}%
g\\
\vdots\\
g
\end{array}
\right\}  n-2\\
\ \ \ \mu_{n}\left[  g_{2},\overset{n-2}{\overbrace{g,\ldots,g}},g_{2}%
^{\prime}\right]
\end{array}
\right)  =\Pi_{2n-2}^{\left(  n,2\right)  }\left(
\begin{array}
[c]{l}%
\left.
\begin{array}
[c]{c}%
g\\
\vdots\\
g
\end{array}
\right\}  n-2\\
\ \ \ g_{1}\circledast g_{1}^{\prime}\\
\left.
\begin{array}
[c]{c}%
g\\
\vdots\\
g
\end{array}
\right\}  n-2\\
\ \ \ g_{2}\circledast g_{2}^{\prime}%
\end{array}
\right)  =\Pi_{2}^{\operatorname*{ret}g}\left(
\begin{array}
[c]{c}%
g_{1}\circledast g_{1}^{\prime}\\
g_{2}\circledast g_{2}^{\prime}%
\end{array}
\right)  .
\end{align*}
}

In the framework of our classification $\Pi_{2}^{\operatorname*{ret}g}\left(
\begin{array}
[c]{c}%
g_{1}\\
g_{2}%
\end{array}
\right)  $ is a hom-like 2-place (binary) representation.
\end{example}

The above formulas describe various properties of multiplace representations,
but they give no idea of how to build representations for concrete polyadic
systems. The most common method of representation construction uses the
concept of a group action on a set (see, e.g., \cite{kirillov}). Below we
extend this concept to the multiplace case and use corresponding
heteromorphisms, as it was done above for homomorphisms and representations.

\section{Multiactions and $G$-spaces}

Let $\mathfrak{G}$$_{n}=\left\langle G;\mu_{n}\right\rangle $ be a polyadic
system and $\mathsf{X}$ be a set. A (left) 1-place action of $\mathfrak{G}%
$$_{n}$ on $\mathsf{X}$ is the external binary operation $\mathbf{\rho}%
_{1}^{(n)}:G\times\mathsf{X}\rightarrow\mathsf{X}$ such that it is consistent
with the multiplication $\mu_{n}$, i.e. composition of the binary operations
$\mathbf{\rho}_{1}\left\{  g|\mathsf{x}\right\}  $ gives the $n$-ary product,
that is,%
\begin{equation}
\mathbf{\rho}_{1}^{(n)}\left\{  \mu_{n}\left[  g_{1},\ldots g_{n}\right]
|\text{$\mathsf{x}$}\right\}  =\mathbf{\rho}_{1}^{(n)}\left\{  g_{1}%
|\mathbf{\rho}_{1}^{(n)}\left\{  g_{2}|\ldots|\mathbf{\rho}_{1}^{(n)}\left\{
g_{n}|\text{$\mathsf{x}$}\right\}  \right\}  \ldots\right\}  ,\ \ \ g_{1}%
,\ldots,g_{n}\in G,\ \text{$\mathsf{x}$}\in\text{$\mathsf{X}$}. \label{rm}%
\end{equation}
If the polyadic system is a $n$-ary group, then in addition to (\ref{rm}) it
is implied the there exist such $e_{x}\in G$ (which may or may not coincide
with the unity of $\mathfrak{G}$$_{n}$) that $\mathbf{\rho}_{1}^{(n)}\left\{
e_{x}|\mathsf{x}\right\}  =\mathsf{x}$ for all $\mathsf{x}\in\mathsf{X}$, and
the mapping $\mathsf{x}\mapsto\mathbf{\rho}_{1}^{(n)}\left\{  e_{x}%
|\mathsf{x}\right\}  $ is a bijection of $\mathsf{X}$. The right 1-place
actions of $\mathfrak{G}$$_{n}$ on $\mathsf{X}$ are defined in a symmetric
way, and therefore we will consider below only one of them. Obviously, we
cannot compose $\mathbf{\rho}_{1}^{(n)}$ and $\mathbf{\rho}_{1}^{(n^{\prime}%
)}$ with $n\neq n^{\prime}$. Usually $\mathsf{X}$ is called a $G$-set or
$G$-space depending on its properties (see, e.g., \cite{hus/joa/jur/sch}).

The application of the 1-place action defined by (\ref{rm}) to the
representation theory of $n$-ary groups gave mostly repetitions of the
ordinary (binary) group representation results (except for trivial $b$-derived
ternary groups) \cite{dud/sha}. Also, it is obviously seen that the
construction (\ref{rm}) with the binary external operation $\mathbf{\rho}_{1}$
cannot be applied for studying the most important regular representations of
polyadic systems, when the $\mathsf{X}$ coincides with $\mathfrak{G}$$_{n}$
itself and the action arises from translations.

Here we introduce the multiplace concept of action for polyadic systems, which
is consistent with heteromorphisms and multiplace representations. Then we
will show how it naturally appears when $\mathsf{X}=\mathfrak{G}_{n}$ and
apply it to construct examples of representations including the regular ones.

For a polyadic system $\mathfrak{G}$$_{n}=\left\langle G;\mu_{n}\right\rangle
$ and a set $\mathsf{X}$ we introduce an external polyadic operation%
\begin{equation}
\mathbf{\rho}_{k}:G^{\times k}\times\text{$\mathsf{X}$}\rightarrow
\text{$\mathsf{X}$}, \label{rk}%
\end{equation}
which is called a (left) $k$-\textit{place action} or \textit{multiaction}. To
generalize the 1-action composition (\ref{rm}), we use the analogy with
multiplication laws of the heteromorphisms (\ref{he}) and the multiplace
representations (\ref{ppe}) and propose (schematically)%
\begin{equation}
\mathbf{\rho}_{k}^{\left(  n\right)  }\left\{  \left.
\genfrac{}{}{0pt}{}{\left.
\begin{array}
[c]{c}%
\mu_{n}\left[  g_{1},\ldots,g_{n}\right]  ,\\
\vdots\\
\mu_{n}\left[  g_{n\left(  \ell_{\mu}-1\right)  },\ldots,g_{n\ell_{\mu}%
}\right]
\end{array}
\right\}  \ell_{\mu}}{\left.
\begin{array}
[c]{c}%
g_{n\ell_{\mu}+1},\\
\vdots\\
g_{n\ell_{\mu}+\ell_{\operatorname*{id}}}%
\end{array}
\right\}  \ell_{\operatorname*{id}}}%
\right\vert \mathsf{x}\right\}  =\mathbf{\rho}_{k}^{(n)}\overset{n^{\prime}%
}{\overbrace{\left\{  \left.
\begin{array}
[c]{c}%
g_{1}\\
\vdots\\
g_{k}%
\end{array}
\right\vert \left.
\begin{array}
[c]{c}%
\ \\
\ldots\\
\
\end{array}
\right\vert \rho_{k}^{(n)}\left\{  \left.
\begin{array}
[c]{c}%
g_{k\left(  n^{\prime}-1\right)  }\\
\vdots\\
g_{kn^{\prime}}%
\end{array}
\right\vert \text{$\mathsf{x}$}\right\}  \ldots\right\}  }}. \label{rrk}%
\end{equation}

The connection between all the parameters here is the same as in the arity
changing formulas (\ref{n1})--(\ref{n2}). Composition of mappings is
associative, and therefore in concrete cases we can use the associative quiver
technique, as it is described in the previous sections. If $\mathfrak{G}$%
$_{n}$ is $n$-ary group, then we should add to (\ref{rrk}) the
\textquotedblleft normalizing\textquotedblright\ relations analogous with
(\ref{pke}) or (\ref{pkg}). So, if there is a unity $e\in$$\mathfrak{G}$$_{n}%
$, then%
\begin{equation}
\mathbf{\rho}_{k}^{(n)}\left\{  \left.
\begin{array}
[c]{c}%
e\\
\vdots\\
e
\end{array}
\right\vert \mathsf{x}\right\}  =\mathsf{x},\ \ \ \ \text{for all }%
\mathsf{x}\in\mathsf{X}. \label{re}%
\end{equation}

In terms of the querelement, the normalization has the form%
\begin{equation}
\mathbf{\rho}_{k}^{(n)}\left\{  \left.
\begin{array}
[c]{l}%
\left.
\begin{array}
[c]{c}%
\bar{h}\\
\vdots\\
\bar{h}%
\end{array}
\right\}  \ell_{\mu}\\[20pt]%
\left.
\begin{array}
[c]{c}%
h\\
\vdots\\
h
\end{array}
\right\}  \ell_{\operatorname*{id}}%
\end{array}
\right\vert \mathsf{x}\right\}  =\mathsf{x},\ \ \ \ \text{for all }%
\mathsf{x}\in\mathsf{X}\ \ \text{and for all }h\in\mathfrak{G}_{n}.
\label{re1}%
\end{equation}

The multiaction $\rho_{k}^{(n)}$ is \textit{transitive}, if any two points
$\mathsf{x}$ and $\mathsf{y}$ in $X$ can be \textquotedblleft
connected\textquotedblright\ by $\rho_{k}^{(n)}$, i.e. there exist
$g_{1},\ldots,g_{k}\in$$G$$_{n}$ such that%
\begin{equation}
\rho_{k}^{(n)}\left\{  \left.
\begin{array}
[c]{c}%
g_{1}\\
\vdots\\
g_{k}%
\end{array}
\right\vert \mathsf{x}\right\}  =\mathsf{y}. \label{rxy}%
\end{equation}

If $g_{1},\ldots,g_{k}$ are unique, then $\rho_{k}^{(n)}$ is \textit{sharply
transitive}. The subset of $\mathsf{X}$, in which any points are connected by
(\ref{rxy}) with fixed $g_{1},\ldots,g_{k}$ can be called the
\textit{multiorbit} of $\mathsf{X}$. If there is only one multiorbit, then we
call $\mathsf{X}$ the \textit{heterogenous} $G$-space (by analogy with the
homogeneous one). By analogy with the (ordinary) 1-place actions, we define a
$G$-equivariant map $\Psi$ between two $G$-sets $\mathsf{X}$ and $\mathsf{Y}$
by (in our notation)%
\begin{equation}
\Psi\left(  \rho_{k}^{(n)}\left\{  \left.
\begin{array}
[c]{c}%
g_{1}\\
\vdots\\
g_{k}%
\end{array}
\right\vert \mathsf{x}\right\}  \right)  =\mathbf{\rho}_{k}^{(n)}\left\{
\left.
\begin{array}
[c]{c}%
g_{1}\\
\vdots\\
g_{k}%
\end{array}
\right\vert \Psi\left(  \mathsf{x}\right)  \right\}  \in\mathsf{Y},
\end{equation}
which makes $G$-space into a category (for details, see, e.g.,
\cite{hus/joa/jur/sch}). In the particular case, when $\mathsf{X}$ is a vector
space over $\mathbb{K}$, the multiaction (\ref{rk}) can be called a
\textit{multi}-$G$-module which satisfies (\ref{re}) and the additional
(linearity) conditions%
\begin{equation}
\rho_{k}^{(n)}\left\{  \left.
\begin{array}
[c]{c}%
g_{1}\\
\vdots\\
g_{k}%
\end{array}
\right\vert a\mathsf{x}+b\mathsf{y}\right\}  =a\mathbf{\rho}_{k}^{(n)}\left\{
\left.
\begin{array}
[c]{c}%
g_{1}\\
\vdots\\
g_{k}%
\end{array}
\right\vert \mathsf{x}\right\}  +b\mathbf{\rho}_{k}^{(n)}\left\{  \left.
\begin{array}
[c]{c}%
g_{1}\\
\vdots\\
g_{k}%
\end{array}
\right\vert \mathsf{y}\right\}  ,
\end{equation}
where $a,b\in\mathbb{K}$. Then, comparing (\ref{ppe}) and (\ref{rrk}) we can
define a multiplace representation as a multi-$G$-module by the following
formula%
\begin{equation}
\Pi_{k}^{(n,n^{\prime})}\left(
\begin{array}
[c]{c}%
g_{1}\\
\vdots\\
g_{k}%
\end{array}
\right)  \left(  \mathsf{x}\right)  =\mathbf{\rho}_{k}^{(n)}\left\{  \left.
\begin{array}
[c]{c}%
g_{1}\\
\vdots\\
g_{k}%
\end{array}
\right\vert \mathsf{x}\right\}  . \label{pr}%
\end{equation}

In a similar way, one can generalize to polyadic systems many other notions
from group action theory \cite{kirillov}.

\section{Regular multiactions}

The most important role in the study of polyadic systems is played by the
case, when $\mathsf{X}=$$\mathfrak{G}$$_{n}$, and the multiaction coincides
with the $n$-ary analog of translations \cite{mal54}, so called $i$%
-translations \cite{belousov}. In the binary case, ordinary translations lead
to regular representations \cite{kirillov}, and therefore we call such an
action a \textit{regular multiaction }$\mathbf{\rho}_{k}^{reg(n)}$. In this
connection, the analog of the Cayley theorem for $n$-ary groups was obtained
in \cite{galm86,galm01a}. Now we will show in examples, how the regular
multiactions can arise from $i$-translations.

\begin{example}
Let $\mathfrak{G}$$_{3}$ be a ternary semigroup, $k=2$, and $\mathsf{X}%
=$$\mathfrak{G}$$_{3}$, then 2-place (left) action can be defined as%
\begin{equation}
\mathbf{\rho}_{2}^{reg\left(  3\right)  }\left\{  \left.
\begin{array}
[c]{c}%
g\\
h
\end{array}
\right\vert u\right\}  \overset{def}{=}\mu_{3}\left[  g,h,u\right]  .
\label{r3}%
\end{equation}
This gives the following composition law for two regular multiactions%
\begin{align}
\mathbf{\rho}_{2}^{reg\left(  3\right)  }\left\{  \left.
\begin{array}
[c]{c}%
g_{1}\\
h_{1}%
\end{array}
\right\vert \mathbf{\rho}_{2}^{reg\left(  3\right)  }\left\{  \left.
\begin{array}
[c]{c}%
g_{2}\\
h_{2}%
\end{array}
\right\vert u\right\}  \right\}   &  =\mu_{3}\left[  g_{1},h_{1},\mu
_{3}\left[  g_{2},h_{2},u\right]  \right] \nonumber\\
&  =\mu_{3}\left[  \mu_{3}\left[  g_{1},h_{1},g_{2}\right]  ,h_{2},u\right]
=\mathbf{\rho}_{2}^{reg\left(  3\right)  }\left\{  \left.
\begin{array}
[c]{c}%
\mu_{3}\left[  g_{1},h_{1},g_{2}\right] \\
h_{2}%
\end{array}
\right\vert u\right\}  . \label{r32}%
\end{align}

Thus, using the regular 2-action (\ref{r3}) we have, in fact, derived the
associative quiver corresponding to (\ref{f2}).
\end{example}

The formula (\ref{r3}) can be simultaneously treated as a 2-translation
\cite{belousov}. In this way, the following \textit{left regular multiaction}%
\begin{equation}
\mathbf{\rho}_{k}^{reg\left(  n\right)  }\left\{  \left.
\begin{array}
[c]{c}%
g_{1}\\
\vdots\\
g_{k}%
\end{array}
\right\vert h\right\}  \overset{def}{=}\mu_{n}\left[  g_{1},\ldots
,g_{k},h\right]  , \label{rh}%
\end{equation}
corresponds to (\ref{fn}), where in the r.h.s. there is the $i$-translation
with $i=n$. The \textit{right regular multiaction} corresponds to the
$i$-translation with $i=1$. The binary composition of the left regular
multiactions corresponds to (\ref{fn}). In general, the value of $i$ fixes the
minimal final arity $n_{reg}^{\prime}$, which differs for even and odd values
of the initial arity $n$.

It follows from (\ref{rh}) that for regular multiactions the number of places
is fixed%
\begin{equation}
k_{reg}=n-1, \label{kr}%
\end{equation}
and the arity changing formulas (\ref{n1})--(\ref{n2}) become%
\begin{align}
n_{reg}^{\prime}  &  =n-\ell_{\operatorname*{id}}\label{nr1}\\
n_{reg}^{\prime}  &  =\ell_{\mu}+1. \label{nr2}%
\end{align}

From (\ref{nr1})--(\ref{nr2}) we conclude that for any $n$ a regular
multiaction having one multiplication $\ell_{\mu}=1$ is binarizing and has
$n-2$ intact elements. For $n=3$ see (\ref{r32}). Also, it follows from
(\ref{nr1}) that for regular multiactions the number of intact elements gives
exactly the difference between initial and final arities.

If the initial arity is odd, then there exists a special \textit{middle
regular multiaction} generated by the $i$-translation with $i=\left(
n+1\right)  \diagup2$. For $n=3$ the corresponding associative quiver is
(\ref{qv5}) and such 2-actions were used in \cite{bor/dud/dup3} to construct
middle representations of ternary groups, which did not change arity
($n^{\prime}=n$). Here we give a more complicated example of a middle regular
multiaction, which can contain intact elements and can therefore change arity.

\begin{example}
Let us consider 5-ary semigroup and the following middle 4-action%
\begin{equation}
\mathbf{\rho}_{4}^{reg\left(  5\right)  }\left\{  \left.
\begin{array}
[c]{c}%
g\\
h\\
u\\
v
\end{array}
\right\vert s\right\}  =\mu_{5}\left[  g,h,\overset{i=3}{\overset{\downarrow
}{s}},u,v\right]  .
\end{equation}
Using (\ref{nr2}) we observe that there are two possibilities for the number
of multiplications $\ell_{\mu}=2,4$. The last case $\ell_{\mu}=4$ is similar
to the vertical associative quiver (\ref{qv5}), but with a more complicated
l.h.s., that is%
\begin{align}
&  \mathbf{\rho}_{4}^{reg\left(  5\right)  }\left\{  \left.
\begin{array}
[c]{c}%
\mu_{5}\left[  g_{1},h_{1},g_{2},h_{2,}g_{3}\right] \\
\mu_{5}\left[  h_{3},g_{4},h_{4},g_{5},h_{5}\right] \\
\mu_{5}\left[  u_{5},v_{5},u_{4},v_{4},u_{3}\right] \\
\mu_{5}\left[  v_{3},u_{2},v_{2},u_{1},v_{1}\right]
\end{array}
\right\vert s\right\}  =\nonumber\\
&  \mathbf{\rho}_{4}^{reg\left(  5\right)  }\left\{  \left.
\begin{array}
[c]{c}%
g_{1}\\
h_{1}\\
u_{1}\\
v_{1}%
\end{array}
\right\vert \mathbf{\rho}_{4}^{reg\left(  5\right)  }\left\{  \left.
\begin{array}
[c]{c}%
g_{2}\\
h_{2}\\
u_{2}\\
v_{2}%
\end{array}
\right\vert \mathbf{\rho}_{4}^{reg\left(  5\right)  }\left\{  \left.
\begin{array}
[c]{c}%
g_{3}\\
h_{3}\\
u_{3}\\
v_{3}%
\end{array}
\right\vert \mathbf{\rho}_{4}^{reg\left(  5\right)  }\left\{  \left.
\begin{array}
[c]{c}%
g_{4}\\
h_{4}\\
u_{4}\\
v_{4}%
\end{array}
\right\vert \mathbf{\rho}_{4}^{reg\left(  5\right)  }\left\{  \left.
\begin{array}
[c]{c}%
g_{5}\\
h_{5}\\
u_{5}\\
v_{5}%
\end{array}
\right\vert s\right\}  \right\}  \right\}  \right\}  \right\}  . \label{rr}%
\end{align}
Now we have an additional case with two intact elements $\ell
_{\operatorname*{id}}$ and two multiplications $\ell_{\mu}=2$ as%
\begin{equation}
\mathbf{\rho}_{4}^{reg\left(  5\right)  }\left\{  \left.
\begin{array}
[c]{c}%
\mu_{5}\left[  g_{1},h_{1},g_{2},h_{2,}g_{3}\right] \\
h_{3}\\
\mu_{5}\left[  h_{3},v_{3},u_{2},v_{2},u_{1}\right] \\
v_{1}%
\end{array}
\right\vert s\right\}  =\mathbf{\rho}_{4}^{reg\left(  5\right)  }\left\{
\left.
\begin{array}
[c]{c}%
g_{1}\\
h_{1}\\
u_{1}\\
v_{1}%
\end{array}
\right\vert \mathbf{\rho}_{4}^{reg\left(  5\right)  }\left\{  \left.
\begin{array}
[c]{c}%
g_{2}\\
h_{2}\\
u_{2}\\
v_{2}%
\end{array}
\right\vert \mathbf{\rho}_{4}^{reg\left(  5\right)  }\left\{  \left.
\begin{array}
[c]{c}%
g_{3}\\
h_{3}\\
u_{3}\\
v_{3}%
\end{array}
\right\vert s\right\}  \right\}  \right\}  , \label{rr3}%
\end{equation}
with arity changing from $n=5$ to $n_{reg}^{\prime}=3$. In addition to
(\ref{rr3}) we have 3 more possible regular multiactions due to the
associativity of $\mu_{5}$, when the multiplication brackets in the sequences
of 6 elements in the first two rows and the second two ones can be shifted independently.
\end{example}

For $n>3$, in addition to left, right and middle multiactions, there exist
intermediate cases. First, observe that the $i$-translations with $i=2$ and
$i=n-1$ immediately fix the final arity $n_{reg}^{\prime}=n$. Therefore, the
composition of multiactions will be similar to (\ref{rr}), but with some
permutations in the l.h.s.

Now we consider some multiplace analogs of regular representations of binary
groups \cite{kirillov}. The straightforward generalization is to consider the
previously introduced regular multiactions (\ref{rh}) in the r.h.s. of
(\ref{pr}). Let $\mathfrak{G}_{n}$ be a finite polyadic associative system and
$\mathbb{K}\mathfrak{G}_{n}$ be a vector space spanned by $\mathfrak{G}_{n}$
(some properties of $n$-ary group rings were considered in
\cite{zek/art1,zek/art2}). This means that any element of $\mathbb{K}%
\mathfrak{G}_{n}$ can be uniquely presented in the form $w=\sum_{l}a_{l}\cdot
h_{l}$, $a_{l}\in\mathbb{K}$, $h_{l}\in G$. Then, using (\ref{rh}) and
(\ref{pr}) we define the $i$-\textit{regular }$k$-\textit{place
representation} by%
\begin{equation}
\Pi_{k}^{reg(i)}\left(
\begin{array}
[c]{c}%
g_{1}\\
\vdots\\
g_{k}%
\end{array}
\right)  \left(  w\right)  =\sum_{l}a_{l}\cdot\mu_{k+1}\left[  g_{1}\ldots
g_{i-1}h_{l}g_{i+1}\ldots g_{k}\right]  . \label{pa}%
\end{equation}

Comparing (\ref{rh}) and (\ref{pa}) one can conclude that all the general
properties of multiplace regular representations are similar to those of the
regular multiactions. If $i=1$ or $i=k$, the multiplace representation is
called a \textit{right} or \textit{left regular representation} respectively.
If $k$ is even, the representation with $i=k\diagup2+1$ is called a
\textit{middle regular representation}. The case $k=2$ was considered in
\cite{bor/dud/dup3} for ternary groups.

\section{Multiplace representations of ternary groups}

Let us consider the case $n=3$, $k=2$ in more detail, paying attention to its
special peculiarities, which corresponds to the $2$-place (bi-element)
representations of ternary groups \cite{bor/dud/dup3}. Let $V$ be a vector
space over $\mathbb{K}$ and $\operatorname*{End}V$ be a set of linear
endomorphisms of $V$. From now on we denote the ternary multiplication by
square brackets only, as follows $\mu_{3}\left[  g_{1},g_{2},g_{3}\right]
\equiv\left[  g_{1}g_{2}g_{3}\right]  $, and use the \textquotedblleft
horizontal\textquotedblright\ notation $\Pi\left(
\begin{array}
[c]{c}%
g_{1}\\
g_{2}%
\end{array}
\right)  \equiv\Pi\left(  g_{1},g_{2}\right)  $.

\begin{definition}
\label{def-left} A \textit{left representation} of a ternary group $G,[\ ])$
in $V$ is a map $\Pi^{L}:G\times G\rightarrow\operatorname*{End}V$ such that
\begin{align}
\Pi^{L}\left(  g_{1},g_{2}\right)  \circ\Pi^{L}\left(  g_{3},g_{4}\right)   &
=\Pi^{L}\left(  \left[  g_{1}g_{2}g_{3}\right]  ,g_{4}\right)  ,\label{p1}\\
\Pi^{L}\left(  g,\overline{g}\right)   &  =\operatorname*{id}\nolimits_{V},
\label{p2a}%
\end{align}
where $g,g_{1},g_{2},g_{3},g_{4}\in G.$
\end{definition}

Replacing in (\ref{p2a}) $g$ by $\overline{g}$ we obtain $\Pi^{L}\left(
\overline{g},g\right)  =id_{V}$, which means that in fact (\ref{p2a}) has the
form $\Pi^{L}\left(  \overline{g},g\right)  =\Pi^{L}\left(  g,\overline
{g}\right)  =\operatorname*{id}\nolimits_{V},\;\;\forall g\in G$. Note that
the axioms considered in the above definition are the natural ones satisfied
by left multiplications $g\mapsto\left[  abg\right]  $. For all $g_{1}%
,g_{2},g_{3},g_{4}\in G$ we have
\[
\Pi^{L}\left(  \left[  g_{1}g_{2}g_{3}\right]  ,g_{4}\right)  =\Pi^{L}\left(
g_{1},\left[  g_{2}g_{3}g_{4}\right]  \right)  .
\]

For all $g,h,u\in G$ we have
\begin{equation}
\Pi^{L}\left(  g,h\right)  =\Pi^{L}\left(  [gu\overline{u}],h\right)  =\Pi
^{L}\left(  g,u\right)  \circ\Pi^{L}\left(  \overline{u},h\right)  \label{pp0}%
\end{equation}
and
\begin{equation}
\Pi^{L}\left(  g,u\right)  \circ\Pi^{L}\left(  \overline{u},\overline
{g}\right)  =\Pi^{L}\left(  \overline{u},\overline{g}\right)  \circ\Pi
^{L}\left(  g,u\right)  =\operatorname*{id}\nolimits_{V},\label{ppa}%
\end{equation}
and therefore every $\Pi^{L}\left(  g,u\right)  $ is invertible and $\left(
\Pi^{L}\left(  g,u\right)  \right)  ^{-1}=\Pi^{L}\left(  \overline
{u},\overline{g}\right)  $. This means that any left representation gives a
representation of a ternary group by a binary group \cite{bor/dud/dup3}. If
the ternary group is medial, then
\[
\Pi^{L}\left(  g_{1},g_{2}\right)  \circ\Pi^{L}\left(  g_{3},g_{4}\right)
=\Pi^{L}\left(  g_{3},g_{4}\right)  \circ\Pi^{L}\left(  g_{1},g_{2}\right)  ,
\]
i.e. the group so obtained is commutative. If the ternary group $\left\langle
G,\left[  \ \right]  \right\rangle $ is commutative, then also $\Pi^{L}\left(
g,h\right)  =\Pi^{L}\left(  h,g\right)  $, because
\begin{align*}
\Pi^{L}\left(  g,h\right)    & =\Pi^{L}\left(  g,h\right)  \circ\Pi^{L}\left(
g,\overline{g}\right)  =\Pi^{L}\left(  \left[  g\,h\,g\right]  ,\,\overline
{g}\right)  \\
& =\Pi^{L}\left(  \left[  h\,g\,g\right]  ,\,\overline{g}\right)  =\Pi
^{L}\left(  h,g\right)  \circ\Pi^{L}\left(  g,\overline{g}\right)  =\Pi
^{L}\left(  h,g\right)  .
\end{align*}
In the case of a commutative and idempotent ternary group any of its left
representations is idempotent and $\left(  \Pi^{L}\left(  g,h\right)  \right)
^{-1}=\Pi^{L}\left(  g,h\right)  $, so that commutative and idempotent ternary
groups are represented by Boolean groups.

\begin{assertion}
\label{xx} Let $\left\langle G,\left[  \ \right]  \right\rangle
=\mathrm{der\,}\left(  G,\odot\right)  $ be a ternary group derived from a
binary group $\left\langle G,\odot\right\rangle $, then there is one-to-one
correspondence between representations of $\left(  G,\odot\right)  $ and left
representations of $\left(  G,\left[  \ \right]  \right)  $.
\end{assertion}

Indeed, because $\left(  G,\left[  \ \right]  \right)  =\mathrm{der}\,\left(
G,\odot\right)  $, then $g\odot h=[geh]$ and $\overline{e}=e$, where $e$ is
unity of the binary group $\left(  G,\odot\right)  $. If $\pi\in
\mathrm{Rep}\left(  G,\odot\right)  $, then (as it is not difficult to see)
$\Pi^{L}\left(  g,h\right)  =\pi\left(  g\right)  \circ\pi\left(  h\right)  $
is a left representation of $\left\langle G,\left[  \ \right]  \right\rangle
$. Conversely, if $\Pi^{L}$ is a left representation of $\left\langle
G,\left[  \ \right]  \right\rangle $ then $\pi\left(  g\right)  =\Pi
^{L}\left(  g,e\right)  $ is a representation of $\left(  G,\odot\right)  $.
Moreover, in this case $\Pi^{L}\left(  g,h\right)  =\pi\left(  g\right)
\circ\pi\left(  h\right)  $, because we have
\[
\Pi^{L}\left(  g,h\right)  =\Pi^{L}\left(  g,[ehe]\right)  =\Pi^{L}\left(
[geh],e\right)  =\Pi^{L}\left(  g,e\right)  \circ\Pi^{L}\left(  h,e\right)
=\pi\left(  g\right)  \circ\pi\left(  h\right)  .
\]

Let $(G,[\ ])$ be a ternary group and $(G\times G,\ast)$ be a semigroup used
in the construction of left representations. According to Post \cite{pos} one
says that two pairs $(a,b)$, $(c,d)$ of elements of $G$ are equivalent, if
there exists an element $g\in G$ such that $[abg]=[cdg]$. Using a covering
group we can see that if this equation holds for some $g\in G$, then it holds
also for all $g\in G$. This means that
\[
\Pi^{L}(a,b)=\Pi^{L}(c,d)\Longleftrightarrow(a,b)\sim(c,d),
\]
i.e.
\[
\Pi^{L}(a,b)=\Pi^{L}(c,d)\Longleftrightarrow\lbrack abg]=[cdg]
\]
for some $g\in G$. Indeed, if $[abg]=[cdg]$ holds for some $g\in G$, then
\begin{align*}
\Pi^{L}(a,b)  &  =\Pi^{L}(a,b)\circ\Pi^{L}(g,\overline{g})=\Pi^{L}%
([abg],\,\overline{g})\\
&  =\Pi^{L}([cdg],\,\overline{g})=\Pi^{L}(c,d)\circ\Pi^{L}(g,\overline{g}%
)=\Pi^{L}(c,d).
\end{align*}

By analogy we can define

\begin{definition}
A \textit{right representation} of a ternary group $(G,[\ ])$ in $V$ is a map
$\Pi^{R}:G\times G\rightarrow\operatorname*{End}\,V$ such that
\begin{align}
\Pi^{R}\left(  g_{3},g_{4}\right)  \circ\Pi^{R}\left(  g_{1},g_{2}\right)   &
=\Pi^{R}\left(  g_{1},\left[  g_{2}g_{3}g_{4}\right]  \right)  ,\label{r1}\\
\Pi^{R}\left(  g,\overline{g}\right)   &  =\operatorname*{id}\nolimits_{V},
\label{r2}%
\end{align}
where $g,g_{1},g_{2},g_{3},g_{4}\in G.$
\end{definition}

From (\ref{r1})-(\ref{r2}) it follows that
\begin{equation}
\Pi^{R}\left(  g,h\right)  =\Pi^{R}\left(  g,\left[  u\,\overline
{u}\,h\right]  \right)  =\Pi^{R}\left(  \overline{u},h\right)  \circ\Pi
^{R}\left(  g,u\right)  . \label{pra}%
\end{equation}

It is easy to check that $\Pi^{R}\left(  g,h\right)  =\Pi^{L}\left(
\overline{h},\overline{g}\right)  =\left(  \Pi^{L}\left(  g,h\right)  \right)
^{-1}$. So it is sufficient to consider only left representations (as in the
binary case). Consider the following example of a group algebra ternary
generalization \cite{bor/dud/dup3}.

\begin{example}
\label{exam-kg}Let $G$ be a ternary group and $\mathbb{K}G$ be a vector space
spanned by $G$, which means that any element of $\mathbb{K}G$ can be uniquely
presented in the form $t=\sum_{i=1}^{n}k_{i}h_{i}$, $k_{i}\in\mathbb{K},$
$h_{i}\in G$, $n\in\mathbb{N}$ (we do not assume that $G$ has finite rank).
Then left and right regular representations are defined by
\begin{align}
\Pi_{reg}^{L}\left(  g_{1},g_{2}\right)  t  &  =\sum_{i=1}^{n}k_{i}\left[
g_{1}g_{2}h_{i}\right]  ,\label{pr1}\\
\Pi_{reg}^{R}\left(  g_{1},g_{2}\right)  t  &  =\sum_{i=1}^{n}k_{i}\left[
h_{i}g_{1}g_{2}\right]  . \label{pr2}%
\end{align}

\end{example}

Let us construct the middle representations as follows.

\begin{definition}
A \textit{middle representation} of a ternary group $\left\langle G,\left[
\ \right]  \right\rangle $ in $V$ is a map $\Pi^{M}:G\times G\rightarrow
\operatorname*{End}\,V$ such that
\begin{align}
\Pi^{M}\left(  g_{3},h_{3}\right)  \circ\Pi^{M}\left(  g_{2},h_{2}\right)
\circ\Pi^{M}\left(  g_{1},h_{1}\right)   &  =\Pi^{M}\left(  \left[  g_{3}%
g_{2}g_{1}\right]  ,\left[  h_{1}h_{2}h_{3}\right]  \right)  ,\label{pm}\\
\Pi^{M}\left(  g,h\right)  \circ\Pi^{M}\left(  \overline{g},\overline
{h}\right)   &  =\Pi^{M}\left(  \overline{g},\overline{h}\right)  \circ\Pi
^{M}\left(  g,h\right)  =\operatorname*{id}\nolimits_{V} \label{pm1}%
\end{align}

\end{definition}

It can be seen that a middle representation is a ternary group homomorphism
$\Pi^{M}:G\times G^{op}\rightarrow\mathrm{der}\operatorname*{End}\,V.$ Note
that instead of (\ref{pm1}) one can use $\Pi^{M}\left(  g,\overline{h}\right)
\circ\Pi^{M}\left(  \overline{g},h\right)  =id_{V}$ after changing $g$ to
$\overline{g}$ and taking into account that $g=\overline{\overline{g}}$. In
the case of idempotent elements $g$ and $h$ we have $\Pi^{M}(g,h)\circ\Pi
^{M}(g,h)=id_{V}$, which means that the matrices $\Pi^{M}$ are Boolean. Thus
all middle representation matrices of idempotent ternary groups are Boolean.
The composition $\Pi^{M}\left(  g_{1},h_{1}\right)  \circ\Pi^{M}\left(
g_{2},h_{2}\right)  $ is not a middle representation, but the following
proposition nevertheless holds.

Let$\;\Pi^{M}$ be a middle representation of a ternary group $\left\langle
G,[\ ]\right\rangle $, then, if$\;\Pi_{u}^{L}(g,h)=\Pi^{M}(g,u)\circ\Pi
^{M}(h,\overline{u})$ is a left representation of $\,\left\langle
G,[\ ]\right\rangle $, then $\;\Pi_{u}^{L}(g,h)\circ\Pi_{u^{\prime}}%
^{L}(g^{\prime},h^{\prime})=\Pi_{u^{\prime}}^{L}(\left[  ghu^{\prime}\right]
,h^{\prime})$, and, if$\;\Pi_{u}^{R}(g,h)=\Pi^{M}(u,h)\circ\Pi^{M}%
(\overline{u},g)$ is a right representation of $\,\left\langle
G,[\ ]\right\rangle $, then $\;\Pi_{u}^{R}(g,h)\circ\Pi_{u^{\prime}}%
^{R}(g^{\prime},h^{\prime})=\Pi_{u}^{R}(g,\left[  hg^{\prime}h^{\prime
}\right]  )$. In particular, $\Pi_{u}^{L}$ ($\Pi_{u}^{R}$) is a family of left
(right) representations.

If a middle representation $\Pi^{M}$ of a ternary group $\left\langle
G,[\ ]\right\rangle $ satisfies $\Pi^{M}\left(  g,\overline{g}\right)
=\mathrm{id}_{V}$ for all $g\in G$, then it is a left and a right
representation and $\Pi^{M}(g,h)=\Pi^{M}(h,g)$ for all $g,h\in G$. Note that
in general $\,\Pi_{reg}^{M}(g,\overline{g})\neq\mathrm{id}$. For regular
representations we have the following commutation relations%
\[
\Pi_{reg}^{L}\left(  g_{1},h_{1}\right)  \circ\Pi_{reg}^{R}\left(  g_{2}%
,h_{2}\right)  =\Pi_{reg}^{R}\left(  g_{2},h_{2}\right)  \circ\Pi_{reg}%
^{L}\left(  g_{1},h_{1}\right)  .
\]

Let $\left\langle G,[\ ]\right\rangle $ be a ternary group and let
$\left\langle G\times G,\left[  \ \right]  ^{\prime}\right\rangle $ be a
ternary group used in the construction of the middle representation. In
$\left\langle G,[\ ]\right\rangle $, and consequently in $\left\langle G\times
G,\left[  \ \right]  ^{\prime}\right\rangle ,$ we define the relation
\[
(a,b)\sim(c,d)\Longleftrightarrow\lbrack aub]=[cud]
\]
for all $u\in G$. It is not difficult to see that this relation is a
congruence in the ternary group $\left\langle G\times G,\left[  \ \right]
^{\prime}\right\rangle $. For regular representations $\Pi_{reg}^{M}%
(a,b)=\Pi_{reg}^{M}(c,d)$ if $(a,b)\sim(c,d)$. We have the following relation
\[
a\eqsim a^{\prime}\Longleftrightarrow a=[\overline{g}a^{\prime}%
g]\;\;\text{for\ some}\;g\in G
\]
or equivalently
\[
a\eqsim a^{\prime}\Longleftrightarrow a^{\prime}=[ga\overline{g}%
]\;\;\text{for\ some}\;g\in G.
\]

It is not difficult to see that it is an equivalence relation on $\left\langle
G,[\ ]\right\rangle $, moreover, if $\left\langle G,[\ ]\right\rangle $ is
medial, then this relation is a congruence.

Let $\left\langle G\times G,\left[  \ \right]  ^{\prime}\right\rangle $ be a
ternary group used in a construction of middle representations, then
\begin{align*}
(a,b)  &  \eqsim(a^{\prime},b)\Longleftrightarrow a^{\prime}=[ga\overline
{g}]\;\;\text{and}\;\;\\
b^{\prime}  &  =[hb\overline{b}\,]\;\;\text{for\ some}\;(g,h)\in G\times G
\end{align*}
is an equivalence relation on $\left\langle G\times G,\left[  \ \right]
^{\prime}\right\rangle $. Moreover, if $(G,[\;])$ is medial, then this
relation is a congruence. Unfortunately, however it is a weak relation. In a
ternary group $\mathbb{Z}_{3}$, where $[ghu]=\left(  g-h+u\right)
(\operatorname{mod}3)$ we have only one class, i.e. all elements are
equivalent. In $\mathbb{Z}_{4}$ with the operation $[ghu]=\left(
g+h+u+1\right)  (\operatorname{mod}\,4)$ we have $a\eqsim a^{\prime
}\Longleftrightarrow a=a^{\prime}$. However, for this relation the following
statement holds. If $(a,b)\eqsim(a^{\prime},b^{\prime})$, then
\[
\operatorname{tr}\Pi^{M}(a,b)=\operatorname{tr}\Pi^{M}(a^{\prime},b^{\prime
}).
\]

We have $\operatorname{tr}(AB)=\operatorname{tr}(BA)$ for all $A,B\in
\mathrm{End}V$, and
\begin{align*}
\operatorname{tr}\Pi^{M}(a,b)  &  =\operatorname{tr}\Pi^{M}([ga^{\prime
}\overline{g}],[hb^{\prime}\overline{h}]\,)=\operatorname{tr}\left(  \Pi
^{M}(g,\overline{h})\circ\Pi^{M}(a^{\prime},b^{\prime})\circ\Pi^{M}%
(\overline{g},h)\right) \\
&  =\operatorname{tr}\left(  \Pi^{M}(g,\overline{h})\circ\Pi^{M}(\overline
{g},h)\circ\Pi^{M}(a^{\prime},b^{\prime})\right)  =\operatorname{tr}\left(
id_{V}\circ\Pi^{M}(a^{\prime}b^{\prime})\right)  =\operatorname{tr}\Pi
^{M}(a^{\prime},b^{\prime})
\end{align*}

In our derived case the connection with standard group representations is
given by the following. Let $\left(  G,\odot\right)  $ be a binary group, and
the ternary derived group as $\left\langle G,[\ ]\right\rangle =\mathrm{der}%
\,\left(  G,\odot\right)  $. There is one-to-one correspondence between a pair
of commuting binary group representations and a middle ternary derived group
representation. Indeed, let $\pi,\rho\in Rep\left(  G,\odot\right)  $,
$\,\pi\left(  g\right)  \circ\rho\left(  h\right)  =\rho\left(  h\right)
\circ\pi\left(  g\right)  \,$ and $\,\Pi^{L}\in Rep\left(  G,\left[
\ \right]  \right)  $. We take
\[
\Pi^{M}\left(  g,h\right)  =\pi\left(  g\right)  \circ\rho\left(
h^{-1}\right)  ,\ \ \pi\left(  g\right)  =\Pi^{M}\left(  g,e\right)
,\ \ \rho\left(  g\right)  =\Pi^{M}\left(  e,\overline{g}\right)  .
\]
Then using (\ref{pm}) we prove the needed representation laws.

Let $\left\langle G,[\ ]\right\rangle $ be a fixed ternary group,
$\left\langle G\times G,\left[  \ \right]  ^{\prime}\right\rangle $ a
corresponding ternary group used in the construction of middle
representations, $(\left(  G\times G\right)  ^{\ast},\circledast)$ a covering
group of $\left\langle G\times G,\left[  \ \right]  ^{\prime}\right\rangle $,
\ $(G\times G,\diamond)=\mathrm{ret}_{(a,b)}(G\times G,\langle\;\rangle)$. If
$\;\Pi^{M}(a,b)\,$ is a middle representation of $\left\langle
G,[\ ]\right\rangle $, then $\pi$ defined by
\[
\pi(g,h,0)=\Pi^{M}(g,h),\ \ \ \pi(g,h,1)=\Pi^{M}(g,h)\circ\Pi^{M}(a,b)
\]
is a representation of the covering group \cite{pos}. Moreover
\[
\rho(g,h)=\Pi^{M}(g,h)\circ\Pi^{M}(a,b)=\pi(g,h,1)
\]
is a representation of the above retract induced by $(a,b)$. Indeed,
$(\overline{a},\overline{b})$ is the identity of this retract and
$\rho(\overline{a},\overline{b})=\Pi^{M}(\overline{a},\overline{b})\circ
\Pi^{M}(a,b)=\operatorname*{id}\nolimits_{V}$. Similarly
\begin{align*}
\rho\left(  (g,h)\diamond(u,u)\right)   &  =\rho\left(  \langle
(g,h),(a,b),(u,u)\rangle\right)  =\rho\left(  \lbrack gau],[ubh]\right)
=\Pi^{M}\left(  [gau],[ubh]\right)  )\circ\Pi^{M}(a,b)\\
&  =\Pi^{M}(g,h)\circ\Pi^{M}(a,b)\circ\Pi^{M}(u,u)\circ\Pi^{M}(a,b)=\rho
(g,h)\circ\rho(u,u)
\end{align*}

But $\tau(g)=(g,\overline{g})$ is an embedding of $(G,[\;])$ into
$\left\langle G\times G,\left[  \ \right]  ^{\prime}\right\rangle $. Hence
$\mu$ defined by $\mu(g,0)=\Pi^{M}(g,\overline{g})$ and $\mu(g,1)=\Pi
^{M}(g,\overline{g})\circ\Pi^{M}(a,\overline{a})$ is a representation of a
covering group $G^{\ast}$ for $(G,[\;])$ (see the Post theorem \cite{pos} for
$a=c$). On the other hand, $\beta(g)=\Pi^{M}(g,\overline{g})\circ\Pi
^{M}(a,\overline{a})$ is a representation of a binary retract $(G,\cdot
\,)=\mathrm{ret}_{a}(G,[\;])$. Thus $\beta$ can induce some middle
representation of $(G,[\;])$ (by the Gluskin-Hossz\'{u} theorem \cite{glu1}).

Note that in the ternary group of quaternions $\left\langle \mathbb{K}%
,[\;]\right\rangle $ (with norm $1)$, where $[ghu]=ghu(-1)=-ghu\,$ and $\,gh$
is the multiplication of quaternions ($-1$ is a central element) we have
$\,\overline{1}=-1$, $\overline{-1}=1\,$ and $\,\overline{g}=g\,$ for others.
In $\left\langle \mathbb{K}\times\mathbb{K},\left[  \ \right]  ^{\prime
}\right\rangle $ we have $(a,b)\sim(-a,-b)$ and $(a,-b)\sim(-a,b)$, which
gives 32 two-element equivalence classes. The embedding $\tau(g)=(g,\overline
{g})$ suggest that $\Pi^{M}(i,i)=\pi(i)\neq\pi(-i)=\Pi^{M}(-i,-i)$. Generally
$\Pi^{M}(a,b)\neq\Pi^{M}(-a,-b)$ and $\Pi^{M}(a,-b)\neq\Pi^{M}(-a,b)$.

The relation $(a,b)\sim(c,d)\Longleftrightarrow\lbrack abg]=[cdg]$ for all
$g\in G$ is a congruence on $(G\times G,\ast)$. Note that this relation can be
defined as "for some $g$". Indeed, using a covering group we can see that if
$\,[abg]=[cdg]\,$ holds for some $g$ then it holds also for all $g$. Thus
$\;\pi^{L}(a,b)=\Pi^{L}(c,d)\Longleftrightarrow(a,b)\sim(c,d)$. Indeed
\begin{align*}
\Pi^{L}(a,b)  &  =\Pi^{L}(a,b)\circ\Pi^{L}(g,\overline{g})=\Pi^{L}%
([a\ b\ g],\overline{g})\\
&  =\Pi^{L}([c\ d\ g],\overline{g})=\Pi^{L}(c,d)\circ\Pi^{L}(g,\overline
{g})=\Pi^{L}(c,d).
\end{align*}

We conclude, that every left representation of a commutative group
$\left\langle G,[\ ]\right\rangle $ is a middle representation. Indeed,
\[
\Pi^{L}(g,h)\circ\Pi^{L}(\overline{g},\overline{h})=\Pi^{L}([g\ h\ \overline
{g}],\overline{h})=\Pi^{L}([g\ \overline{g}\ h],\overline{h})=\Pi
^{L}(h,\overline{h})=\mathrm{id}_{V}%
\]
and
\begin{align*}
&  \Pi^{L}(g_{1},g_{2})\circ\Pi^{L}(g_{3},g_{4})\circ\Pi^{L}(g_{5},g_{6})
=\Pi^{L}([[g_{1}g_{2}g_{3}]g_{4}g_{5}],g_{6})=\Pi^{L}([[g_{1}g_{3}g_{2}%
]g_{4}g_{5}],g_{6})\\
&  =\Pi^{L}([g_{1}g_{3}[g_{2}g_{4}g_{5}]],g_{6}) =\Pi^{L}([g_{1}g_{3}%
[g_{5}g_{4}g_{2}]],g_{6})=\Pi^{L}([g_{1}g_{3}g_{5}],[g_{4}g_{2}g_{6}])=\Pi
^{L}([g_{1}g_{3}g_{5}],[g_{6}g_{4}g_{2}]).
\end{align*}

Note that the converse holds only for the special kind of middle
representations such that $\Pi^{M}(g,\overline{g})=\mathrm{id}_{V}$. Therefore,

\begin{assertion}
There is one-one correspondence between left representations of $\left\langle
G,[\ ]\right\rangle $ and binary representations of the retract $\mathrm{ret}%
_{a}(G,[\ ])$.
\end{assertion}

Indeed, let $\Pi^{L}(g,a)$ be given, then $\rho(g)=\Pi^{L}(g,a)$ is such
representation of the retract, as can be directly shown. Conversely, assume
that $\rho(g)$ is a representation of the retract $\mathrm{ret}_{a}(G,[\;])$.
Define $\,\Pi^{L}(g,h)=\rho(g)\circ\rho(\overline{h})^{-1}$, then $\,\Pi
^{L}(g,h)\circ\Pi^{L}(u,u)=\rho(g)\circ\rho(\overline{h})^{-1}\circ
\rho(u)\circ\rho(\overline{u})^{-1}=\rho(g\circledast(\overline{h})^{-1}%
\circ\circledast u)\circ\rho(\overline{u})^{-1}=\rho([\,[\,g\,a\,[\,\overline
{a}\,h\,\overline{a}\,]\,]\,a\,u\,])\circ\rho(\overline{u})^{-1}%
=\rho([\,g\,h\,g\,])\circ\rho(\overline{u})^{-1}=\Pi^{L}([\,g\,h\,u\,],u),$

\section{Matrix representations of ternary groups}

Here we give several examples of matrix representations for concrete ternary
groups. Let $G=\mathbb{Z}_{3}\ni\left\{  0,1,2\right\}  $ and the ternary
multiplication be $\left[  ghu\right]  =g-h+u$. Then $\left[  ghu\right]
=\left[  uhg\right]  $ and $\,\overline{0}=0,$ $\,\overline{1}=1,$
$\,\overline{2}=2$, therefore $(G,[\ ])$ is an idempotent medial ternary
group. Thus $\,\Pi^{L}(g,h)=\Pi^{R}(h,g)\,$ and
\begin{equation}
\Pi^{L}(a,b)=\Pi^{L}(c,d)\Longleftrightarrow
(a-b)=(c-d)\mathrm{\operatorname{mod}}\,3. \label{mod3}%
\end{equation}
The calculations give the left regular representation in the manifest matrix
form%
\begin{align}
\Pi_{reg}^{L}\left(  0,0\right)   &  =\Pi_{reg}^{L}\left(  2,2\right)
=\Pi_{reg}^{L}\left(  1,1\right)  =\Pi_{reg}^{R}\left(  0,0\right) \nonumber\\
&  =\Pi_{reg}^{R}\left(  2,2\right)  =\Pi_{reg}^{R}\left(  1,1\right)
=\left(
\begin{array}
[c]{ccc}%
1 & 0 & 0\\
0 & 1 & 0\\
0 & 0 & 1
\end{array}
\right)  =[1]\oplus\lbrack1]\oplus\lbrack1],
\end{align}%
\begin{align}
\Pi_{reg}^{L}\left(  2,0\right)   &  =\Pi_{reg}^{L}\left(  1,2\right)
=\Pi_{reg}^{L}\left(  0,1\right)  =\Pi_{reg}^{R}\left(  2,1\right)  =\Pi
_{reg}^{R}\left(  1,0\right)  =\Pi_{reg}^{R}\left(  0,2\right)  =\left(
\begin{array}
[c]{ccc}%
0 & 1 & 0\\
0 & 0 & 1\\
1 & 0 & 0
\end{array}
\right) \nonumber\\
&  =[1]\oplus\left(
\begin{array}
[c]{cc}%
-\dfrac{1}{2} & -\dfrac{\sqrt{3}}{2}\\
\dfrac{\sqrt{3}}{2} & -\dfrac{1}{2}%
\end{array}
\right)  =[1]\oplus\left[  -\dfrac{1}{2}+\dfrac{1}{2}i\sqrt{3}\right]
\oplus\left[  -\dfrac{1}{2}-\dfrac{1}{2}i\sqrt{3}\right]  ,
\end{align}%
\begin{align}
\Pi_{reg}^{L}\left(  2,1\right)   &  =\Pi_{reg}^{L}\left(  1,0\right)
=\Pi_{reg}^{L}\left(  0,2\right)  =\Pi_{reg}^{R}\left(  2,0\right)  =\Pi
_{reg}^{R}\left(  1,2\right)  =\Pi_{reg}^{R}\left(  0,1\right)  =\left(
\begin{array}
[c]{ccc}%
0 & 0 & 1\\
1 & 0 & 0\\
0 & 1 & 0
\end{array}
\right) \nonumber\\
&  =[1]\oplus\left(
\begin{array}
[c]{cc}%
-\dfrac{1}{2} & \dfrac{\sqrt{3}}{2}\\
-\dfrac{\sqrt{3}}{2} & -\dfrac{1}{2}%
\end{array}
\right)  =[1]\oplus\left[  -\dfrac{1}{2}-\dfrac{1}{2}i\sqrt{3}\right]
\oplus\left[  -\dfrac{1}{2}+\dfrac{1}{2}i\sqrt{3}\right]  .
\end{align}

Consider next the middle representation construction. The middle regular
representation is defined by
\[
\Pi_{reg}^{M}\left(  g_{1},g_{2}\right)  t=\sum_{i=1}^{n}k_{i}\left[
g_{1}h_{i}g_{2}\right]  .
\]
For regular representations we have
\begin{align}
\Pi_{reg}^{M}\left(  g_{1},h_{1}\right)  \circ\Pi_{reg}^{R}\left(  g_{2}%
,h_{2}\right)   &  =\Pi_{reg}^{R}\left(  h_{2},h_{1}\right)  \circ\Pi
_{reg}^{M}\left(  g_{1},g_{2}\right)  ,\label{pp30}\\
\Pi_{reg}^{M}\left(  g_{1},h_{1}\right)  \circ\Pi_{reg}^{L}\left(  g_{2}%
,h_{2}\right)   &  =\Pi_{reg}^{L}\left(  g_{1},g_{2}\right)  \circ\Pi
_{reg}^{M}\left(  h_{2},h_{1}\right)  . \label{pp3}%
\end{align}

For the middle regular representation matrices we obtain
\begin{align*}
\Pi_{reg}^{M}\left(  0,0\right)   &  =\Pi_{reg}^{M}\left(  1,2\right)
=\Pi_{reg}^{M}\left(  2,1\right)  =\left(
\begin{array}
[c]{ccc}%
1 & 0 & 0\\
0 & 0 & 1\\
0 & 1 & 0
\end{array}
\right)  ,\\
\Pi_{reg}^{M}\left(  0,1\right)   &  =\Pi_{reg}^{M}\left(  1,0\right)
=\Pi_{reg}^{M}\left(  2,2\right)  =\left(
\begin{array}
[c]{ccc}%
0 & 1 & 0\\
1 & 0 & 0\\
0 & 0 & 1
\end{array}
\right)  ,\\
\Pi_{reg}^{M}\left(  0,2\right)   &  =\Pi_{reg}^{M}\left(  2,0\right)
=\Pi_{reg}^{M}\left(  1,1\right)  =\left(
\begin{array}
[c]{ccc}%
0 & 0 & 1\\
0 & 1 & 0\\
1 & 0 & 0
\end{array}
\right)  .
\end{align*}

The above representation $\Pi_{reg}^{M}$ of $\left\langle \mathbb{Z}%
_{3},\left[  \;\right]  \right\rangle $ is equivalent to the orthogonal direct
sum of two irreducible representations
\begin{align*}
\Pi_{reg}^{M}\left(  0,0\right)   &  =\Pi_{reg}^{M}\left(  1,2\right)
=\Pi_{reg}^{M}\left(  2,1\right)  =\left[  1\right]  \oplus\left[
\begin{array}
[c]{cc}%
-1 & 0\\
0 & 1
\end{array}
\right]  ,\\
\Pi_{reg}^{M}\left(  0,1\right)   &  =\Pi_{reg}^{M}\left(  1,0\right)
=\Pi_{reg}^{M}\left(  2,2\right)  =\left[  1\right]  \oplus\left[
\begin{array}
[c]{cc}%
\dfrac{1}{2} & -\dfrac{\sqrt{3}}{2}\\
-\dfrac{\sqrt{3}}{2} & -\dfrac{1}{2}%
\end{array}
\right]  ,\\
\Pi_{reg}^{M}\left(  0,2\right)   &  =\Pi_{reg}^{M}\left(  2,0\right)
=\Pi_{reg}^{M}\left(  1,1\right)  =\left[  1\right]  \oplus\left[
\begin{array}
[c]{cc}%
\dfrac{1}{2} & \dfrac{\sqrt{3}}{2}\\
\dfrac{\sqrt{3}}{2} & -\dfrac{1}{2}%
\end{array}
\right]  ,
\end{align*}
i.e. one-dimensional trivial $\left[  1\right]  $ and two-dimensional
irreducible. Note, that in this example $\,\Pi^{M}(g,\overline{g})=\Pi
^{M}(g,g)\neq\mathrm{id}_{V}$, but $\,\Pi^{M}(g,h)\circ\Pi^{M}%
(g,h)=\mathrm{id}_{V}$, and so $\,\Pi^{M}$ are of the second degree.

Consider a more complicated example of left representations. \label{exam-z4}%
Let $G=\mathbb{Z}_{4}\ni\left\{  0,1,2,3\right\}  $ and the ternary
multiplication be
\begin{equation}
\left[  ghu\right]  =\left(  g+h+u+1\right)  \operatorname{mod}4.
\label{mult4}%
\end{equation}
We have the multiplication table
\begin{align*}
\left[  g,h,0\right]   &  =\left(
\begin{array}
[c]{cccc}%
1 & 2 & 3 & 0\\
2 & 3 & 0 & 1\\
3 & 0 & 1 & 2\\
0 & 1 & 2 & 3
\end{array}
\right)  \ \ \ \ \ \ \ \left[  g,h,1\right]  =\left(
\begin{array}
[c]{cccc}%
2 & 3 & 0 & 1\\
3 & 0 & 1 & 2\\
0 & 1 & 2 & 3\\
1 & 2 & 3 & 0
\end{array}
\right) \\
\left[  g,h,2\right]   &  =\left(
\begin{array}
[c]{cccc}%
3 & 0 & 1 & 2\\
0 & 1 & 2 & 3\\
1 & 2 & 3 & 0\\
2 & 3 & 0 & 1
\end{array}
\right)  \ \ \ \ \ \ \ \left[  g,h,3\right]  =\left(
\begin{array}
[c]{cccc}%
0 & 1 & 2 & 3\\
1 & 2 & 3 & 0\\
2 & 3 & 0 & 1\\
3 & 0 & 1 & 2
\end{array}
\right)
\end{align*}
Then the skew elements are $\,\overline{0}=3,$ $\,\overline{1}=2,$
$\,\overline{2}=1,$ $\,\overline{3}=0$, and therefore $(G,[\ ])$ is a
(non-idempotent) commutative ternary group. The left representation is defined
by the expansion $\,\Pi_{reg}^{L}\left(  g_{1},g_{2}\right)  t=\sum_{i=1}%
^{n}k_{i}\left[  g_{1}g_{2}h_{i}\right]  $, which means that (see the general
formula (\ref{pa}))
\[
\Pi_{reg}^{L}\left(  g,h\right)  |u>=|\left[  ghu\right]  >.
\]
Analogously, for right and middle representations
\[
\Pi_{reg}^{R}\left(  g,h\right)  |u>=|\left[  ugh\right]  >,\ \ \ \Pi
_{reg}^{M}\left(  g,h\right)  |u>=|\left[  guh\right]  >.
\]
Therefore $\,|\left[  ghu\right]  >=|\left[  ugh\right]  >=|\left[
guh\right]  >$ and
\[
\Pi_{reg}^{L}\left(  g,h\right)  =\Pi_{reg}^{R}\left(  g,h\right)
|u>=\Pi_{reg}^{M}\left(  g,h\right)  |u>,
\]
so $\Pi_{reg}^{L}\left(  g,h\right)  =\Pi_{reg}^{R}\left(  g,h\right)
=\Pi_{reg}^{M}\left(  g,h\right)  $. Thus it is sufficient to consider the
left representation only.

In this case the equivalence is $\,\Pi^{L}(a,b)=\Pi^{L}%
(c,d)\Longleftrightarrow(a+b)=(c+d)$\textrm{$\operatorname{mod}$}$\,4$, and we
obtain the following classes
\begin{align*}
\Pi_{reg}^{L}\left(  0,0\right)   &  =\Pi_{reg}^{L}\left(  1,3\right)
=\Pi_{reg}^{L}\left(  2,2\right)  =\Pi_{reg}^{L}\left(  3,1\right)  =\left(
\begin{array}
[c]{cccc}%
0 & 0 & 0 & 1\\
1 & 0 & 0 & 0\\
0 & 1 & 0 & 0\\
0 & 0 & 1 & 0
\end{array}
\right)  =\left[  1\right]  \oplus\left[  -1\right]  \oplus\left[  -i\right]
\oplus\left[  i\right]  ,\\
\Pi_{reg}^{L}\left(  0,1\right)   &  =\Pi_{reg}^{L}\left(  1,0\right)
=\Pi_{reg}^{L}\left(  2,3\right)  =\Pi_{reg}^{L}\left(  3,2\right)  =\left(
\begin{array}
[c]{cccc}%
0 & 0 & 1 & 0\\
0 & 0 & 0 & 1\\
1 & 0 & 0 & 0\\
0 & 1 & 0 & 0
\end{array}
\right)  =\left[  1\right]  \oplus\left[  -1\right]  \oplus\left[  -1\right]
\oplus\left[  -1\right]  ,\\
\Pi_{reg}^{L}\left(  0,2\right)   &  =\Pi_{reg}^{L}\left(  1,1\right)
=\Pi_{reg}^{L}\left(  2,0\right)  =\Pi_{reg}^{L}\left(  3,3\right)  =\left(
\begin{array}
[c]{cccc}%
0 & 1 & 0 & 0\\
0 & 0 & 1 & 0\\
0 & 0 & 0 & 1\\
1 & 0 & 0 & 0
\end{array}
\right)  =\left[  1\right]  \oplus\left[  -1\right]  \oplus\left[  i\right]
\oplus\left[  -i\right]  ,\\
\Pi_{reg}^{L}\left(  0,3\right)   &  =\Pi_{reg}^{L}\left(  1,2\right)
=\Pi_{reg}^{L}\left(  2,1\right)  =\Pi_{reg}^{L}\left(  3,0\right)  =\left(
\begin{array}
[c]{cccc}%
1 & 0 & 0 & 0\\
0 & 1 & 0 & 0\\
0 & 0 & 1 & 0\\
0 & 0 & 0 & 1
\end{array}
\right)  =\left[  1\right]  \oplus\left[  -1\right]  \oplus\left[  1\right]
\oplus\left[  1\right]  .
\end{align*}
It is seen that, due to the fact that the ternary operation (\ref{mult4}) is
commutative, there are only one-dimensional irreducible left representations.

Let us \textquotedblleft algebralize\textquotedblright\ the above regular
representations in the following way. From (\ref{p1}) we have, for the left
representation
\begin{equation}
\Pi_{reg}^{L}\left(  i,j\right)  \circ\Pi_{reg}^{L}\left(  k,l\right)
=\Pi_{reg}^{L}\left(  i,\left[  jkl\right]  \right)  ,
\end{equation}
where $\,\left[  jkl\right]  =j-k+l$, $\;i,j,k,l\in\mathbb{Z}_{3}$. Denote
$\,\gamma_{i}^{L}=\Pi_{reg}^{L}\left(  0,i\right)  $, $\;i\in\mathbb{Z}_{3}$,
then we obtain the algebra with the relations%
\begin{equation}
\gamma_{i}^{L}\gamma_{j}^{L}=\gamma_{i+j}^{L}.
\end{equation}
Conversely, any matrix representation of $\,\gamma_{i}\gamma_{j}=\gamma
_{i+j}\,$ leads to the left representation by $\,\Pi^{L}\left(  i,j\right)
=\gamma_{j-i}$. In the case of the middle regular representation we introduce
$\gamma_{k+l}^{M}=\Pi_{reg}^{M}\left(  k,l\right)  $, $\;k,l\in\mathbb{Z}_{3}%
$, then we obtain
\begin{equation}
\gamma_{i}^{M}\gamma_{j}^{M}\gamma_{k}^{M}=\gamma_{\left[  ijk\right]  }%
^{M},\;\;\;i,j,k\in\mathbb{Z}_{3}. \label{ggg}%
\end{equation}

In some sense (\ref{ggg}) can be treated as a \textit{ternary analog of the
Clifford algebra}. As before, any matrix representation of (\ref{ggg}) gives
the middle representation $\,\Pi^{M}\left(  k,l\right)  =\gamma_{k+l}$.

\section{Ternary algebras and Hopf algebras}

Let us consider associative ternary algebras \cite{car4,azc/izq}. One can
introduce an autodistributivity property $\left[  \left[  xyz\right]
ab\right]  =\left[  \left[  xab\right]  \left[  yab\right]  \left[
zab\right]  \right]  $ (see \cite{dud1}). If we take 2 ternary operations
$\left\{  \;,\;,\;\right\}  $ and $\left[  \;,\;,\;\right]  $, then
distributivity is given by $\left\{  \left[  xyz\right]  ab\right\}  =\left[
\left\{  xab\right\}  \left\{  yab\right\}  \left\{  zab\right\}  \right]  $.
If $\left(  +\right)  $ is a binary operation (addition), then \textit{left
linearity} is
\begin{equation}
\left[  \left(  x+z\right)  ,a,b\right]  =\left[  xab\right]  +\left[
zab\right]  .
\end{equation}
By analogy one can define central (middle) and right linearity.
\textit{Linearity} is defined, when left, middle and right linearity hold simultaneously.

\begin{definition}
An\textit{ associative ternary algebra} is a triple $\left(  A,\mu_{3}%
,\eta^{\left(  3\right)  }\right)  $, where $A$ is a linear space over a field
$\mathbb{K}$, $\mu_{3}$ is a linear map $A\otimes A\otimes A\rightarrow A$
called \textit{ternary multiplication} $\mu_{3}\left(  a\otimes b\otimes
c\right)  =\left[  abc\right]  $ which is ternary associative $\left[  \left[
abc\right]  de\right]  =\left[  a\left[  bcd\right]  e\right]  =\left[
ab\left[  cde\right]  \right]  $ or
\begin{equation}
\mu_{3}\circ\left(  \mu_{3}\otimes\operatorname*{id}\otimes\operatorname*{id}%
\right)  =\mu_{3}\circ\left(  \operatorname*{id}\otimes\mu_{3}\otimes
\operatorname*{id}\right)  =\mu_{3}\circ\left(  \operatorname*{id}%
\otimes\operatorname*{id}\otimes\mu_{3}\right)  . \label{ass}%
\end{equation}

\end{definition}

There are two types \cite{dup26} of ternary unit maps $\eta^{\left(  3\right)
}:\mathbb{K}\rightarrow A$:

1) One \textit{strong} unit map
\begin{equation}
\mu_{3}\circ\left(  \eta^{\left(  3\right)  }\otimes\eta^{\left(  3\right)
}\otimes\operatorname*{id}\right)  =\mu_{3}\circ\left(  \eta^{\left(
3\right)  }\otimes\operatorname*{id}\otimes\eta^{\left(  3\right)  }\right)
=\mu_{3}\circ\left(  \operatorname*{id}\otimes\eta^{\left(  3\right)  }%
\otimes\eta^{\left(  3\right)  }\right)  =\operatorname*{id};
\end{equation}

2) Two \textit{sequential} units $\eta_{1}^{\left(  3\right)  }$ and $\eta
_{2}^{\left(  3\right)  }$ satisfying
\begin{equation}
\mu_{3}\circ\left(  \eta_{1}^{\left(  3\right)  }\otimes\eta_{2}^{\left(
3\right)  }\otimes\operatorname*{id}\right)  =\mu_{3}\circ\left(  \eta
_{1}^{\left(  3\right)  }\otimes\operatorname*{id}\otimes\eta_{2}^{\left(
3\right)  }\right)  =\mu_{3}\circ\left(  \operatorname*{id}\otimes\eta
_{1}^{\left(  3\right)  }\otimes\eta_{2}^{\left(  3\right)  }\right)
=\operatorname*{id};
\end{equation}

In first case the ternary analog of the binary relation $\eta^{\left(
2\right)  }\left(  x\right)  =x1$, where $x\in\mathbb{K}$, $1\in A$, is%
\begin{equation}
\eta^{\left(  3\right)  }\left(  x\right)  =\left[  x,1,1\right]  =\left[
1,1,x\right]  =\left[  1,x,1\right]  .
\end{equation}

Let $(A,\mu_{A},\eta_{A})$, $(B,\mu_{B},\eta_{B})$ and $(C,\mu_{C},\eta_{C})$
be ternary algebras, then the \textit{ternary tensor product} space
$A\mathbf{\otimes}B\mathbf{\otimes}C$ is naturally endowed with the structure
of an algebra. The multiplication $\mu_{A\mathbf{\otimes}B\mathbf{\otimes}C}$
on $A\mathbf{\otimes}B\mathbf{\otimes}C$ reads
\begin{equation}
\left[  (a_{1}\mathbf{\otimes}b_{1}\mathbf{\otimes}c_{1})(a_{2}\mathbf{\otimes
}b_{2}\mathbf{\otimes}c_{2})(a_{3}\mathbf{\otimes}b_{3}\mathbf{\otimes}%
c_{3})\right]  =\left[  a_{1}a_{2}a_{3}\right]  \mathbf{\otimes}\left[
b_{1}b_{2}b_{3}\right]  \mathbf{\otimes}\left[  c_{1}c_{2}c_{3}\right]  ,
\end{equation}
and so the set of ternary algebras is closed under taking ternary tensor
products. A \textit{ternary algebra map} (homomorphism) is a linear map
between ternary algebras $f:A\rightarrow B$ which respects the ternary algebra
structure
\begin{align}
f\left(  \left[  xyz\right]  \right)   &  =\left[  f\left(  x\right)
,f\left(  y\right)  ,f\left(  z\right)  \right]  ,\\
f\left(  1_{A}\right)   &  =1_{B}.
\end{align}

Let $C$ be a linear space over a field $\mathbb{K}$.

\begin{definition}
A\textit{ ternary comultiplication} $\Delta_{3}$ is a linear map over a field
$\mathbb{K}$ such that
\begin{equation}
\Delta_{3}:C\rightarrow C\otimes C\otimes C. \label{ccc}%
\end{equation}

\end{definition}

In the standard Sweedler notations \cite{sweedler} $\Delta_{3}\left(
a\right)  =\sum_{i=1}^{n}a_{i}^{\prime}\otimes a_{i}^{\prime\prime}\otimes
a_{i}^{\prime\prime\prime}=a_{\left(  1\right)  }\otimes a_{\left(  2\right)
}\otimes a_{\left(  3\right)  }$. Consider different possible types of ternary
coassociativity \cite{dup26,bor/dud/dup1}.

\begin{enumerate}
\item A\textit{ standard} ternary coassociativity
\begin{equation}
(\Delta_{3}\otimes\operatorname*{id}\otimes\operatorname*{id})\circ\Delta
_{3}=(\operatorname*{id}\otimes\Delta_{3}\otimes\operatorname*{id})\circ
\Delta_{3}=(\operatorname*{id}\otimes\operatorname*{id}\otimes\Delta_{3}%
)\circ\Delta_{3}, \label{ast}%
\end{equation}

\item A\textit{ nonstandard} ternary $\Sigma$-coassociativity (Gluskin-type
positional operatives)
\[
(\Delta_{3}\otimes\operatorname*{id}\otimes\operatorname*{id})\circ\Delta
_{3}=(\operatorname*{id}\otimes\left(  \sigma\circ\Delta_{3}\right)
\otimes\operatorname*{id})\circ\Delta_{3},
\]
where $\sigma\circ\Delta_{3}\left(  a\right)  =\Delta_{3}\left(  a\right)
=a_{\left(  \sigma\left(  1\right)  \right)  }\otimes a_{\left(  \sigma\left(
2\right)  \right)  }\otimes a_{\left(  \sigma\left(  3\right)  \right)  }$ and
$\sigma\in\Sigma\subset S_{3}.$

\item A\textit{ permutational} ternary coassociativity
\[
(\Delta_{3}\otimes\operatorname*{id}\otimes\operatorname*{id})\circ\Delta
_{3}=\pi\circ(\operatorname*{id}\otimes\Delta_{3}\otimes\operatorname*{id}%
)\circ\Delta_{3},
\]
where $\pi\in\Pi\subset S_{5}$.
\end{enumerate}

A\textit{ ternary} \textit{comediality} is
\[
\left(  \Delta_{3}\otimes\Delta_{3}\otimes\Delta_{3}\right)  \circ\Delta
_{3}=\sigma_{medial}\circ\left(  \Delta_{3}\otimes\Delta_{3}\otimes\Delta
_{3}\right)  \circ\Delta_{3},
\]
where $\sigma_{medial}=\binom{123456789}{147258369}\in S_{9}$. A\textit{
ternary counit} is defined as a map $\varepsilon^{\left(  3\right)
}:C\rightarrow\mathbb{K}$. In general, $\varepsilon^{\left(  3\right)  }%
\neq\varepsilon^{\left(  2\right)  }$ satisfying one of the conditions below.
If $\Delta_{3}$ is derived, then maybe $\varepsilon^{\left(  3\right)
}=\varepsilon^{\left(  2\right)  }$, but another counits may exist. There are
two types of ternary counits:

\begin{enumerate}
\item Standard (\textit{strong}) ternary counit
\begin{equation}
(\varepsilon^{\left(  3\right)  }\otimes\varepsilon^{\left(  3\right)
}\otimes\operatorname*{id})\circ\Delta_{3}=(\varepsilon^{\left(  3\right)
}\otimes\operatorname*{id}\otimes\varepsilon^{\left(  3\right)  })\circ
\Delta_{3}=(\operatorname*{id}\otimes\varepsilon^{\left(  3\right)  }%
\otimes\varepsilon^{\left(  3\right)  })\circ\Delta_{3}=\operatorname*{id},
\label{e1}%
\end{equation}

\item Two \textit{sequential} (polyadic) counits $\varepsilon_{1}^{\left(
3\right)  }$ and $\varepsilon_{2}^{\left(  3\right)  }$%
\begin{equation}
(\varepsilon_{1}^{\left(  3\right)  }\otimes\varepsilon_{2}^{\left(  3\right)
}\otimes\operatorname*{id})\circ\Delta=(\varepsilon_{1}^{\left(  3\right)
}\otimes\operatorname*{id}\otimes\varepsilon_{2}^{\left(  3\right)  }%
)\circ\Delta=(\operatorname*{id}\otimes\varepsilon_{1}^{\left(  3\right)
}\otimes\varepsilon_{2}^{\left(  3\right)  })\circ\Delta=\operatorname*{id},
\label{e2}%
\end{equation}

\end{enumerate}

Below we will consider only the first standard type of associativity
(\ref{ast}). The $\sigma$-\textit{cocommutativity} is defined as\textbf{\ }%
$\sigma\circ\Delta_{3}=\Delta_{3}$.

\begin{definition}
A\textit{ ternary coalgebra} is a triple $\left(  C,\Delta_{3},\varepsilon
^{\left(  3\right)  }\right)  $, where $C$ is a linear space and $\Delta_{3}$
is a ternary comultiplication (\ref{ccc}) which is coassociative in one of the
above senses and $\varepsilon^{\left(  3\right)  }$ is one of the above counits.
\end{definition}

Let $\left(  A,\mu_{3}\right)  $ be a ternary algebra and $\left(
C,\Delta_{3}\right)  $ be a ternary coalgebra and $f,g,h\in\operatorname*{Hom}%
_{\mathbb{K}}\left(  C,A\right)  $. \textit{Ternary convolution product} is
\begin{equation}
\left[  f,g,h\right]  _{\ast}=\mu_{3}\circ\left(  f\otimes g\otimes h\right)
\circ\Delta_{3} \label{tconv}%
\end{equation}
or in the Sweedler notation $\left[  f,g,h\right]  _{\ast}\left(  a\right)
=\left[  f\left(  a_{\left(  1\right)  }\right)  g\left(  a_{\left(  2\right)
}\right)  h\left(  a_{\left(  3\right)  }\right)  \right]  $.

\begin{definition}
A\textit{ }ternary coalgebra is called \textit{derived}, if there exists a
binary (usual, see e.g. \cite{sweedler}) coalgebra $\Delta_{2}:C\rightarrow
C\otimes C$ such that%
\begin{equation}
\Delta_{3,\operatorname*{der}}=\left(  \operatorname*{id}\otimes\Delta
_{2}\right)  \otimes\Delta_{2}. \label{dder}%
\end{equation}

\end{definition}

\begin{definition}
A\textit{ }ternary bialgebra $B$ is $\left(  B,\mu_{3},\eta^{\left(  3\right)
},\Delta_{3},\varepsilon^{\left(  3\right)  }\right)  $ for which $\left(
B,\mu_{3},\eta^{\left(  3\right)  }\right)  $ is a ternary algebra and
$\left(  B,\Delta_{3},\varepsilon^{\left(  3\right)  }\right)  $ is a ternary
coalgebra and they are compatible
\begin{equation}
\Delta_{3}\circ\mu_{3}=\mu_{3}\circ\left(  \operatorname*{id}\otimes\tau
_{12}\otimes\operatorname*{id}\right)  \circ\Delta_{3}, \label{md}%
\end{equation}
where $\tau_{12}=\binom{12}{21}$.
\end{definition}

One can distinguish four kinds of ternary bialgebras with respect to a
\textquotedblleft being derived\textquotedblright\ property:

\begin{enumerate}
\item A\textit{ }$\Delta$-\textit{derived} ternary bialgebra
\begin{equation}
\Delta_{3}=\Delta_{3,\operatorname*{der}}=\left(  \operatorname*{id}%
\otimes\Delta_{2}\right)  \circ\Delta_{2} \label{d3}%
\end{equation}

\item A\textit{ }$\mu$-\textit{derived} ternary bialgebra
\begin{equation}
\mu_{3,\operatorname*{der}}=\mu_{2}\circ\left(  \mu_{2}\otimes
\operatorname*{id}\right)  \label{m3a}%
\end{equation}

\item A\textit{ derived} ternary bialgebra is simultaneously $\mu$-derived and
$\Delta$-derived ternary bialgebra.

\item A\textit{ non-derived} ternary bialgebra which does not satisfy
(\ref{d3}) and (\ref{m3a}).
\end{enumerate}

Possible types of ternary antipodes can be defined by analogy with binary coalgebras.

\begin{definition}
A\textit{ skew ternary antipode} is
\begin{equation}
\mu_{3}\circ(S_{skew}^{\left(  3\right)  }\otimes\operatorname*{id}%
\otimes\operatorname*{id})\circ\Delta_{3}=\mu_{3}\circ(\operatorname*{id}%
\otimes S_{skew}^{\left(  3\right)  }\otimes\operatorname*{id})\circ\Delta
_{3}=\mu_{3}\circ(\operatorname*{id}\otimes\operatorname*{id}\otimes
S_{skew}^{\left(  3\right)  })\circ\Delta_{3}=\operatorname*{id}. \label{as}%
\end{equation}

\end{definition}

If only one equality from (\ref{as}) is satisfied, the corresponding skew
antipode is called \textit{left}, \textit{middle} or \textit{right}.

\begin{definition}
\textit{Strong ternary antipode} is
\[
\left(  \mu_{2}\otimes\operatorname*{id}\right)  \circ(\operatorname*{id}%
\otimes S_{strong}^{\left(  3\right)  }\otimes\operatorname*{id})\circ
\Delta_{3}=1\otimes\operatorname*{id},\left(  \operatorname*{id}\otimes\mu
_{2}\right)  \circ(\operatorname*{id}\otimes\operatorname*{id}\otimes
S_{strong}^{\left(  3\right)  })\circ\Delta_{3}=\operatorname*{id}\otimes1,
\]
where $1$ is a unit of algebra.
\end{definition}

If in a ternary coalgebra the relation%
\begin{equation}
\Delta_{3}\circ S=\tau_{13}\circ\left(  S\otimes S\otimes S\right)
\circ\Delta_{3}%
\end{equation}
holds true, where $\tau_{13}=\binom{123}{321}$, then it is called
\textit{skew-involutive}.

\begin{definition}
A\textit{ }ternary Hopf algebra $\left(  H,\mu_{3},\eta^{\left(  3\right)
},\Delta_{3},\varepsilon^{\left(  3\right)  },S^{\left(  3\right)  }\right)  $
is a ternary bialgebra with a ternary antipode $S^{\left(  3\right)  }$ of the
corresponding above type .
\end{definition}

Let us consider concrete constructions of ternary comultiplications,
bialgebras and Hopf algebras. A \textit{ternary group-like element} can be
defined by $\Delta_{3}\left(  g\right)  =g\otimes g\otimes g$, and for 3 such
elements we have%
\begin{equation}
\Delta_{3}\left(  \left[  g_{1}g_{2}g_{3}\right]  \right)  =\Delta_{3}\left(
g_{1}\right)  \Delta_{3}\left(  g_{2}\right)  \Delta_{3}\left(  g_{3}\right)
.
\end{equation}
But an analog of the binary primitive element (satisfying $\Delta_{2}\left(
x\right)  =x\otimes1+1\otimes x$) cannot be chosen simply as $\Delta
_{3}\left(  x\right)  =x\otimes e\otimes e+e\otimes x\otimes e+e\otimes
e\otimes x$, since the algebra structure is not preserved. Nevertheless, if we
introduce \textit{two} idempotent units $e_{1},e_{2}$ satisfying
\textquotedblleft semiorthogonality\textquotedblright\ $\left[  e_{1}%
e_{1}e_{2}\right]  =0$, $\left[  e_{2}e_{2}e_{1}\right]  =0$, then
\begin{equation}
\Delta_{3}\left(  x\right)  =x\otimes e_{1}\otimes e_{2}+e_{2}\otimes x\otimes
e_{1}+e_{1}\otimes e_{2}\otimes x \label{exx}%
\end{equation}
and now $\Delta_{3}\left(  \left[  x_{1}x_{2}x_{3}\right]  \right)  =\left[
\Delta_{3}\left(  x_{1}\right)  \Delta_{3}\left(  x_{2}\right)  \Delta
_{3}\left(  x_{3}\right)  \right]  $. Using (\ref{exx}) $\varepsilon\left(
x\right)  =0$, $\varepsilon\left(  e_{1,2}\right)  =1$, and $S^{\left(
3\right)  }\left(  x\right)  =-x$, $S^{\left(  3\right)  }\left(
e_{1,2}\right)  =e_{1,2}$, one can construct a ternary universal enveloping
algebra in full analogy with the binary case (see e.g. \cite{kassel}).

One of the most important examples of noncocommutative Hopf algebras is the
well known Sweedler Hopf algebra \cite{sweedler} which in the binary case has
two generators $x$ and $y$ satisfying
\begin{align}
\mu_{2}\left(  x,x\right)   &  =1,\\
\mu_{2}\left(  y,y\right)   &  =0,\\
\sigma_{+}^{\left(  2\right)  }\left(  xy\right)   &  =-\sigma_{-}^{\left(
2\right)  }\left(  xy\right)  .
\end{align}
It has the following comultiplication
\begin{align}
\Delta_{2}\left(  x\right)   &  =x\otimes x,\\
\Delta_{2}\left(  y\right)   &  =y\otimes x+1\otimes y,
\end{align}
counit $\varepsilon^{\left(  2\right)  }\left(  x\right)  =1$, $\varepsilon
^{\left(  2\right)  }\left(  y\right)  =0$, and antipode $S^{\left(  2\right)
}\left(  x\right)  =x$, $S^{\left(  2\right)  }\left(  y\right)  =-y$, which
respect the algebra structure. In the derived case a \textit{ternary Sweedler
algebra} is generated also by two generators $x$ and $y$ obeying \cite{dup26}
\begin{align}
\mu_{3}\left(  x,e,x\right)   &  =\mu_{3}\left(  e,x,x\right)  =\mu_{3}\left(
x,x,e\right)  =e,\\
\sigma_{+}^{\left(  3\right)  }\left(  \left[  yey\right]  \right)   &  =0,\\
\sigma_{+}^{\left(  3\right)  }\left(  \left[  xey\right]  \right)   &
=-\sigma_{-}^{\left(  3\right)  }\left(  \left[  xey\right]  \right)  .
\end{align}
The derived Hopf algebra structure is given by
\begin{align}
\Delta_{3}\left(  x\right)   &  =x\otimes x\otimes x,\label{h1}\\
\Delta_{3}\left(  y\right)   &  =y\otimes x\otimes x+e\otimes y\otimes
x+e\otimes e\otimes y,\\
\varepsilon^{\left(  3\right)  }\left(  x\right)   &  =\varepsilon^{\left(
2\right)  }\left(  x\right)  =1,\label{h2}\\
\varepsilon^{\left(  3\right)  }\left(  y\right)   &  =\varepsilon^{\left(
2\right)  }\left(  y\right)  =0,\\
S^{\left(  3\right)  }\left(  x\right)   &  =S^{\left(  2\right)  }\left(
x\right)  =x,\label{h3}\\
S^{\left(  3\right)  }\left(  y\right)   &  =S^{\left(  2\right)  }\left(
y\right)  =-y,
\end{align}
and it can be checked that (\ref{h1})-(\ref{h2}) are algebra maps, while
(\ref{h3}) are antialgebra maps. To obtain a non-derived ternary Sweedler
example we have the following possibilities: 1) one \textquotedblleft
even\textquotedblright\ generator $x$, two \textquotedblleft
odd\textquotedblright\ generators $y_{1,2}$ and one ternary unit $e$; 2) two
\textquotedblleft even\textquotedblright\ generators $x_{1,2}$, one
\textquotedblleft odd\textquotedblright\ generator $y$ and two ternary units
$e_{1,2}$. In the first case the ternary algebra structure is (no summation,
$i=1,2$)
\begin{align}
\left[  xxx\right]   &  =e,\;\;\label{ha}\\
\left[  y_{i}y_{i}y_{i}\right]   &  =0,\;\;\\
\sigma_{+}^{\left(  3\right)  }\left(  \left[  y_{i}xy_{i}\right]  \right)
&  =\sigma_{+}^{\left(  3\right)  }\left(  \left[  xy_{i}x\right]  \right)
=0,\\
\left[  xey_{i}\right]   &  =-\left[  xy_{i}e\right]  ,\nonumber\\
\left[  exy_{i}\right]   &  =-\left[  y_{i}xe\right]  ,\\
\left[  ey_{i}x\right]   &  =-\left[  y_{i}ex\right]  ,\\
\sigma_{+}^{\left(  3\right)  }\left(  \left[  y_{1}xy_{2}\right]  \right)
&  =-\sigma_{-}^{\left(  3\right)  }\left(  \left[  y_{1}xy_{2}\right]
\right)  .
\end{align}

The corresponding ternary Hopf algebra structure is
\begin{align}
\Delta_{3}\left(  x\right)   &  =x\otimes x\otimes x,\;\Delta_{3}\left(
y_{1,2}\right)  =y_{1,2}\otimes x\otimes x+e_{1,2}\otimes y_{2,1}\otimes
x+e_{1,2}\otimes e_{2,1}\otimes y_{2,1},\\
\varepsilon^{\left(  3\right)  }\left(  x\right)   &  =1,\;\;\varepsilon
^{\left(  3\right)  }\left(  y_{i}\right)  =0,\;\;\\
S^{\left(  3\right)  }\left(  x\right)   &  =x,\;\;\;\;S^{\left(  3\right)
}\left(  y_{i}\right)  =-y_{i}. \label{hh2}%
\end{align}

In the second case we have for the algebra structure
\begin{align}
\left[  x_{i}x_{j}x_{k}\right]   &  =\delta_{ij}\delta_{ik}\delta_{jk}%
e_{i},\;\;\;\left[  yyy\right]  =0,\;\;\\
\sigma_{+}^{\left(  3\right)  }\left(  \left[  yx_{i}y\right]  \right)   &
=0,\;\sigma_{+}^{\left(  3\right)  }\left(  \left[  x_{i}yx_{i}\right]
\right)  =0,\\
\;\;\;\sigma_{+}^{\left(  3\right)  }\left(  \left[  y_{1}xy_{2}\right]
\right)   &  =0,\;\;\;\sigma_{-}^{\left(  3\right)  }\left(  \left[
y_{1}xy_{2}\right]  \right)  =0, \label{ha1}%
\end{align}
and the ternary Hopf algebra structure is \cite{dup26}
\begin{align}
\Delta_{3}\left(  x_{i}\right)   &  =x_{i}\otimes x_{i}\otimes x_{i}%
,\;\nonumber\\
\Delta_{3}\left(  y\right)   &  =y\otimes x_{1}\otimes x_{1}+e_{1}\otimes
y\otimes x_{2}+e_{1}\otimes e_{2}\otimes y,\\
\varepsilon^{\left(  3\right)  }\left(  x_{i}\right)   &  =1,\label{hhh2}\\
\varepsilon^{\left(  3\right)  }\left(  y\right)   &  =0,\\
S^{\left(  3\right)  }\left(  x_{i}\right)   &  =x_{i},\\
S^{\left(  3\right)  }\left(  y\right)   &  =-y.
\end{align}

\section{Ternary quantum groups}

A ternary commutator can be obtained in different ways \cite{bre/hen}. We will
consider the simplest version called a Nambu bracket (see e.g.
\cite{tak3,azc/izq}). Let us introduce two maps $\omega_{\pm}^{\left(
3\right)  }:A\otimes A\otimes A\rightarrow A\otimes A\otimes A$ by
\begin{align}
\omega_{+}^{\left(  3\right)  }\left(  a\mathbf{\otimes}b\mathbf{\otimes
}c\right)   &  =a\otimes b\otimes c+b\otimes c\otimes a+c\otimes a\otimes
b,\label{w1}\\
\omega_{-}^{\left(  3\right)  }\left(  a\mathbf{\otimes}b\mathbf{\otimes
}c\right)   &  =b\otimes a\otimes c+c\otimes b\otimes a+a\otimes c\otimes b.
\label{w2}%
\end{align}
Thus, obviously $\mu_{3}\circ\omega_{\pm}^{\left(  3\right)  }=\sigma_{\pm
}^{\left(  3\right)  }\circ\mu_{3}$, where $\sigma_{\pm}^{\left(  3\right)
}\in S_{3}$ denotes a sum of terms having even and odd permutations
respectively. In the binary case $\omega_{+}^{\left(  2\right)  }%
=\operatorname*{id}\otimes\operatorname*{id}$ and $\omega_{-}^{\left(
2\right)  }=\tau$ is the twist operator $\tau:a\mathbf{\otimes}b\rightarrow
b\otimes a$, while $\mu_{2}\circ\omega_{-}^{\left(  2\right)  }$ is
permutation $\sigma_{-}^{\left(  2\right)  }\left(  ab\right)  =ba$. So the
Nambu product is $\omega_{N}^{\left(  3\right)  }=\omega_{+}^{\left(
3\right)  }-\omega_{-}^{\left(  3\right)  }$, and the ternary commutator is
$\left[  ,,\right]  _{N}=\sigma_{N}^{\left(  3\right)  }=\sigma_{+}^{\left(
3\right)  }-\sigma_{-}^{\left(  3\right)  }$, or \cite{tak3}%
\begin{equation}
\left[  a,b,c\right]  _{N}=\left[  abc\right]  +\left[  bca\right]  +\left[
cab\right]  -\left[  cba\right]  -\left[  acb\right]  -\left[  bac\right]
\end{equation}

An \textit{abelian ternary algebra} is defined by the vanishing of the Nambu
bracket $\left[  a,b,c\right]  _{N}=0$ or \textit{ternary commutation}
relation $\sigma_{+}^{\left(  3\right)  }=\sigma_{-}^{\left(  3\right)  }$. By
analogy with the binary case a \textit{deformed ternary algebra} can be
defined by
\begin{equation}
\sigma_{+}^{\left(  3\right)  }=q\sigma_{-}^{\left(  3\right)  }\text{ or
}\left[  abc\right]  +\left[  bca\right]  +\left[  cab\right]  =q\left(
\left[  cba\right]  +\left[  acb\right]  +\left[  bac\right]  \right)  ,
\label{sq}%
\end{equation}
where multiplication by $q$ is treated as an external operation.

Let us consider a ternary analog of the Woronowicz example of a bialgebra
construction, which has two generators satisfying $xy=qyx$ (or $\sigma
_{+}^{\left(  2\right)  }\left(  xy\right)  =q\sigma_{-}^{\left(  2\right)
}\left(  xy\right)  $), then the following coproducts
\begin{align}
\Delta_{2}\left(  x\right)   &  =x\otimes x\\
\Delta_{2}\left(  y\right)   &  =y\otimes x+1\otimes y
\end{align}
are the algebra maps. In the derived ternary case using (\ref{sq}) we have%
\begin{equation}
\sigma_{+}^{\left(  3\right)  }\left(  \left[  xey\right]  \right)
=q\sigma_{-}^{\left(  3\right)  }\left(  \left[  xey\right]  \right)  ,
\end{equation}
where $e$ is the ternary unit and ternary coproducts are
\begin{align}
\Delta_{3}\left(  e\right)   &  =e\otimes e\otimes e,\\
\Delta_{3}\left(  x\right)   &  =x\otimes x\otimes x,\\
\Delta_{3}\left(  y\right)   &  =y\otimes x\otimes x+e\otimes y\otimes
x+e\otimes e\otimes y,
\end{align}
which are \textit{ternary algebra maps}, i.e. they satisfy
\begin{equation}
\sigma_{+}^{\left(  3\right)  }\left(  \left[  \Delta_{3}\left(  x\right)
\Delta_{3}\left(  e\right)  \Delta_{3}\left(  y\right)  \right]  \right)
=q\sigma_{-}^{\left(  3\right)  }\left(  \left[  \Delta_{3}\left(  x\right)
\Delta_{3}\left(  e\right)  \Delta_{3}\left(  y\right)  \right]  \right)  .
\end{equation}

Let us consider the group $G=SL\left(  n,\mathbb{K}\right)  $. Then the
algebra generated by $a_{j}^{i}\in SL\left(  n,\mathbb{K}\right)  $ can be
endowed with the structure of a ternary Hopf algebra (see, e.g., \cite{madore}
for the binary case) by choosing the ternary coproduct, counit and antipode as
(here summation is implied)
\begin{equation}
\Delta_{3}\left(  a_{j}^{i}\right)  =a_{k}^{i}\otimes a_{l}^{k}\otimes
a_{j}^{l},\;\;\;\varepsilon\left(  a_{j}^{i}\right)  =\delta_{j}%
^{i},\;\;\;S^{\left(  3\right)  }\left(  a_{j}^{i}\right)  =\left(
a^{-1}\right)  _{j}^{i}. \label{sl2}%
\end{equation}
This antipode is a skew one since from (\ref{as}) it follows that%
\begin{equation}
\mu_{3}\circ(S^{\left(  3\right)  }\otimes\operatorname*{id}\otimes
\operatorname*{id})\circ\Delta_{3}\left(  a_{j}^{i}\right)  =S^{\left(
3\right)  }\left(  a_{k}^{i}\right)  a_{l}^{k}a_{j}^{l}=\left(  a^{-1}\right)
_{k}^{i}a_{l}^{k}a_{j}^{l}=\delta_{l}^{i}a_{j}^{l}=a_{j}^{i}.
\end{equation}
This ternary Hopf algebra is derived since for $\Delta_{2}=a_{j}^{i}\otimes
a_{k}^{j}$ we have
\begin{equation}
\Delta_{3}=\left(  \operatorname*{id}\otimes\Delta_{2}\right)  \otimes
\Delta_{2}\left(  a_{j}^{i}\right)  =\left(  \operatorname*{id}\otimes
\Delta_{2}\right)  \left(  a_{k}^{i}\otimes a_{j}^{k}\right)  =a_{k}%
^{i}\otimes\Delta_{2}\left(  a_{j}^{k}\right)  =a_{k}^{i}\otimes a_{l}%
^{k}\otimes a_{j}^{l}.
\end{equation}

In the most important case $n=2$ we can obtain the manifest action of the
ternary coproduct $\Delta_{3}$ on components. Possible non-derived matrix
representations of the ternary product can be done only by four-rank $n\times
n\times n\times n$ twice covariant and twice contravariant tensors $\left\{
a_{kl}^{ij}\right\}  $. Among all products the non-derived ones are only the
following: $a_{jk}^{oi}b_{oo}^{jl}c_{il}^{ko}$ and $a_{ok}^{ij}b_{io}%
^{ol}c_{il}^{ko}$ (where $o$ is any index). So using e.g. the first choice we
can define the non-derived Hopf algebra structure by%
\begin{align}
\Delta_{3}\left(  a_{kl}^{ij}\right)   &  =a_{v\rho}^{i\mu}\otimes
a_{kl}^{v\sigma}\otimes a_{\mu\sigma}^{\rho j},\\
\varepsilon\left(  a_{kl}^{ij}\right)   &  =\dfrac{1}{2}\left(  \delta_{k}%
^{i}\delta_{l}^{j}+\delta_{l}^{i}\delta_{k}^{j}\right)  ,
\end{align}
and the skew antipode $s_{kl}^{ij}=S^{\left(  3\right)  }\left(  a_{kl}%
^{ij}\right)  $ which is a solution of the equation $s_{v\rho}^{i\mu}%
a_{kl}^{v\sigma}=\delta_{\rho}^{i}\delta_{k}^{\mu}\delta_{l}^{\sigma}$.

Next consider ternary dual pair $kG$ (push-forward) and $\mathcal{F}\left(
G\right)  $ (pull-back) which are related by $\left(  kG\right)  ^{\ast}%
\cong\mathcal{F}\left(  G\right)  $ (see, e.g., \cite{kog/soi}). Here
$kG=\mathrm{span}\left(  G\right)  $ is a ternary group algebra ($G$ has a
ternary product $\left[  \ ,\ ,\ \right]  _{G}$ or $\mu_{G}^{\left(  3\right)
}$) over a field $k$.

If $u\in kG$ ($u=u^{i}x_{i},x_{i}\in G$), then%
\begin{equation}
\left[  uvw\right]  _{k}=u^{i}v^{j}w^{l}\left[  x_{i}x_{j}x_{l}\right]  _{G}%
\end{equation}
is associative, and so $\left(  kG,\left[  \ ,\ ,\ \right]  _{k}\right)  $
becomes a ternary algebra. Define a ternary coproduct $\Delta_{3}%
:kG\rightarrow kG\otimes kG\otimes kG$ by
\begin{equation}
\Delta_{3}\left(  u\right)  =u^{i}x_{i}\otimes x_{i}\otimes x_{i}%
\end{equation}
(derived and associative), then $\Delta_{3}\left(  \left[  uvw\right]
_{k}\right)  =\left[  \Delta_{3}\left(  u\right)  \Delta_{3}\left(  v\right)
\Delta_{3}\left(  w\right)  \right]  _{k}$, and $kG$ is a ternary bialgebra.
If we define a ternary antipode by $S_{k}^{\left(  3\right)  }=u^{i}\bar
{x}_{i}$, where $\bar{x}_{i}$ is a skew element of $x_{i}$, then $kG$ becomes
a ternary Hopf algebra.

In the dual case of functions $\mathcal{F}\left(  G\right)  :\left\{
\varphi:G\rightarrow k\right\}  $ a ternary product $\left[  \ ,\ ,\ \right]
_{\mathcal{F}}$ or $\mu_{3,\mathcal{F}}$ (derived and associative) acts on
$\psi\left(  x,y,z\right)  $ as
\begin{equation}
\left(  \mu_{3,\mathcal{F}}\psi\right)  \left(  x\right)  =\psi\left(
x,x,x\right)  ,
\end{equation}
and so $\mathcal{F}\left(  G\right)  $ is a ternary algebra. Let
$\mathcal{F}\left(  G\right)  \otimes\mathcal{F}\left(  G\right)
\otimes\mathcal{F}\left(  G\right)  \cong\mathcal{F}\left(  G\times G\times
G\right)  $, then we define a ternary coproduct $\Delta_{3}:\mathcal{F}\left(
G\right)  \rightarrow\mathcal{F}\left(  G\right)  \otimes\mathcal{F}\left(
G\right)  \otimes\mathcal{F}\left(  G\right)  $ as
\begin{equation}
\left(  \Delta_{3}\varphi\right)  \left(  x,y,z\right)  =\varphi\left(
\left[  xyz\right]  _{\mathcal{F}}\right)  ,
\end{equation}
which is derive and associative. Thus we can obtain $\Delta_{3}\left(  \left[
\varphi_{1}\varphi_{2}\varphi_{3}\right]  _{\mathcal{F}}\right)  =\left[
\Delta_{3}\left(  \varphi_{1}\right)  \Delta_{3}\left(  \varphi_{2}\right)
\Delta_{3}\left(  \varphi_{3}\right)  \right]  _{\mathcal{F}}$, and therefore
$\mathcal{F}\left(  G\right)  $ is a ternary bialgebra. If we define a ternary
antipode by%
\begin{equation}
S_{\mathcal{F}}^{\left(  3\right)  }\left(  \varphi\right)  =\varphi\left(
\bar{x}\right)  ,
\end{equation}
where $\bar{x}$ is a skew element of $x$, then $\mathcal{F}\left(  G\right)  $
becomes a ternary Hopf algebra.

Let us introduce a ternary analog of the $R$-matrix \cite{dup26}. For a
ternary Hopf algebra $H$ we consider a linear map $R^{\left(  3\right)
}:H\otimes H\otimes H\rightarrow H\otimes H\otimes H$.

\begin{definition}
A ternary Hopf algebra $\left(  H,\mu_{3},\eta^{\left(  3\right)  },\Delta
_{3},\varepsilon^{\left(  3\right)  },S^{\left(  3\right)  }\right)  $ is
called \textit{quasifiveangular}\footnote{The reason for such notation is
clear from (\ref{r5}).} if it satisfies
\begin{align}
\left(  \Delta_{3}\otimes\operatorname*{id}\otimes\operatorname*{id}\right)
&  =R_{145}^{\left(  3\right)  }R_{245}^{\left(  3\right)  }R_{345}^{\left(
3\right)  },\label{r1a}\\
\left(  \operatorname*{id}\otimes\Delta_{3}\otimes\operatorname*{id}\right)
&  =R_{125}^{\left(  3\right)  }R_{145}^{\left(  3\right)  }R_{135}^{\left(
3\right)  },\label{r2a}\\
\left(  \operatorname*{id}\otimes\operatorname*{id}\otimes\Delta_{3}\right)
&  =R_{125}^{\left(  3\right)  }R_{124}^{\left(  3\right)  }R_{123}^{\left(
3\right)  }, \label{r3a}%
\end{align}
where as usual the index of $R$ denotes action component positions.
\end{definition}

Using the standard procedure (see, e.g., \cite{kassel,cha/pre,majid}) we
obtain a set of abstract \textit{ternary quantum Yang-Baxter equations}, one
of which has the form \cite{dup26}
\begin{equation}
R_{243}^{\left(  3\right)  }R_{342}^{\left(  3\right)  }R_{125}^{\left(
3\right)  }R_{145}^{\left(  3\right)  }R_{135}^{\left(  3\right)  }%
=R_{123}^{\left(  3\right)  }R_{132}^{\left(  3\right)  }R_{145}^{\left(
3\right)  }R_{245}^{\left(  3\right)  }R_{345}^{\left(  3\right)  },
\label{r5}%
\end{equation}
and others can be obtained by corresponding permutations. The classical
ternary Yang-Baxter equations form a one parameter family of solutions
$R\left(  t\right)  $ can be obtained by the expansion%
\begin{equation}
R^{\left(  3\right)  }\left(  t\right)  =e\otimes e\otimes e+rt+\mathcal{O}%
\left(  t^{2}\right)  ,
\end{equation}
where $r$ is a ternary classical $R$-matrix, then e.g. for (\ref{r5}) we have
\begin{align*}
&  r_{342}r_{125}r_{145}r_{135}+r_{243}r_{125}r_{145}r_{135}+r_{243}%
r_{342}r_{145}r_{135}+r_{243}r_{342}r_{125}r_{135}+r_{243}r_{342}%
r_{125}r_{145}\\
&  =r_{132}r_{145}r_{245}r_{345}+r_{123}r_{145}r_{245}r_{345}+r_{123}%
r_{132}r_{245}r_{345}+r_{123}r_{132}r_{145}r_{345}+r_{123}r_{132}%
r_{145}r_{245}.
\end{align*}

For three ternary Hopf algebras $\left(  H_{A,B,C},\mu_{A,B,C}^{\left(
3\right)  },\eta_{A,B,C}^{\left(  3\right)  },\Delta_{A,B,C}^{\left(
3\right)  },\varepsilon_{A,B,C}^{\left(  3\right)  },S_{A,B,C}^{\left(
3\right)  }\right)  $ we can introduce a non-degenerate ternary
\textquotedblleft pairing\textquotedblright\ (see, e.g., \cite{cha/pre} for
the binary case) $\left\langle \;,\;,\;\right\rangle ^{\left(  3\right)
}:H_{A}\times H_{B}\times H_{C}\rightarrow\mathbb{K}$, trilinear over
$\mathbb{K}$, satisfying \cite{dup26}
\begin{align*}
\left\langle \eta_{A}^{\left(  3\right)  }\left(  a\right)  ,b,c\right\rangle
^{\left(  3\right)  }  &  =\left\langle a,\varepsilon_{B}^{\left(  3\right)
}\left(  b\right)  ,c\right\rangle ^{\left(  3\right)  },\;\;\left\langle
a,\eta_{B}^{\left(  3\right)  }\left(  b\right)  ,c\right\rangle ^{\left(
3\right)  }=\left\langle \varepsilon_{A}^{\left(  3\right)  }\left(  a\right)
,b,c\right\rangle ^{\left(  3\right)  },\\
\left\langle b,\eta_{B}^{\left(  3\right)  }\left(  b\right)  ,c\right\rangle
^{\left(  3\right)  }  &  =\left\langle a,b,\varepsilon_{C}^{\left(  3\right)
}\left(  c\right)  \right\rangle ^{\left(  3\right)  },\;\;\left\langle
a,b,\eta_{C}^{\left(  3\right)  }\left(  c\right)  \right\rangle ^{\left(
3\right)  }=\left\langle a,\varepsilon_{B}^{\left(  3\right)  }\left(
b\right)  ,c\right\rangle ^{\left(  3\right)  },\\
\left\langle a,b,\eta_{C}^{\left(  3\right)  }\left(  c\right)  \right\rangle
^{\left(  3\right)  }  &  =\left\langle \varepsilon_{A}^{\left(  3\right)
}\left(  a\right)  ,b,c\right\rangle ^{\left(  3\right)  },\;\;\left\langle
\eta_{A}^{\left(  3\right)  }\left(  a\right)  ,b,c\right\rangle ^{\left(
3\right)  }=\left\langle a,b,\varepsilon_{C}^{\left(  3\right)  }\left(
c\right)  \right\rangle ^{\left(  3\right)  },
\end{align*}%
\begin{align*}
\left\langle \mu_{A}^{\left(  3\right)  }\left(  a_{1}\otimes a_{2}\otimes
a_{3}\right)  ,b,c\right\rangle ^{\left(  3\right)  }  &  =\left\langle
a_{1}\otimes a_{2}\otimes a_{3},\Delta_{B}^{\left(  3\right)  }\left(
b\right)  ,c\right\rangle ^{\left(  3\right)  },\\
\left\langle \Delta_{A}^{\left(  3\right)  }\left(  a\right)  ,b_{1}\otimes
b_{2}\otimes b_{3},c\right\rangle ^{\left(  3\right)  }  &  =\left\langle
a,\mu_{B}^{\left(  3\right)  }\left(  b_{1}\otimes b_{2}\otimes b_{3}\right)
,c\right\rangle ^{\left(  3\right)  },\\
\left\langle a,\mu_{B}^{\left(  3\right)  }\left(  b_{1}\otimes b_{2}\otimes
b_{3}\right)  ,c\right\rangle ^{\left(  3\right)  }  &  =\left\langle
a,b_{1}\otimes b_{2}\otimes b_{3},\Delta_{C}^{\left(  3\right)  }\left(
c\right)  \right\rangle ^{\left(  3\right)  },\\
\left\langle a,\Delta_{B}^{\left(  3\right)  }\left(  b\right)  ,c_{1}\otimes
c_{2}\otimes c_{3}\right\rangle ^{\left(  3\right)  }  &  =\left\langle
a,b,\mu_{C}^{\left(  3\right)  }\left(  c_{1}\otimes c_{2}\otimes
c_{3}\right)  \right\rangle ^{\left(  3\right)  },\\
\left\langle a,b,\mu_{C}^{\left(  3\right)  }\left(  c_{1}\otimes c_{2}\otimes
c_{3}\right)  \right\rangle ^{\left(  3\right)  }  &  =\left\langle \Delta
_{A}^{\left(  3\right)  }\left(  a\right)  ,b,c_{1}\otimes c_{2}\otimes
c_{3}\right\rangle ^{\left(  3\right)  },\\
\left\langle a_{1}\otimes a_{2}\otimes a_{3},b,\Delta_{C}^{\left(  3\right)
}\left(  c\right)  \right\rangle ^{\left(  3\right)  }  &  =\left\langle
\mu_{A}^{\left(  3\right)  }\left(  a_{1}\otimes a_{2}\otimes a_{3}\right)
,b,c\right\rangle ^{\left(  3\right)  },
\end{align*}%
\[
\left\langle S_{A}^{\left(  3\right)  }\left(  a\right)  ,b,c\right\rangle
^{\left(  3\right)  }=\left\langle a,S_{B}^{\left(  3\right)  }\left(
b\right)  ,c\right\rangle ^{\left(  3\right)  }=\left\langle a,b,S_{C}%
^{\left(  3\right)  }\left(  c\right)  \right\rangle ^{\left(  3\right)  },
\]
where $a,a_{i}\in H_{A}$, $b,b_{i}\in H_{B}$. The ternary \textquotedblleft
paring\textquotedblright\ between $H_{A}\otimes H_{A}\otimes H_{A}$ and
$H_{B}\otimes H_{B}\otimes H_{B}$ is given by $\left\langle a_{1}\otimes
a_{2}\otimes a_{3},b_{1}\otimes b_{2}\otimes b_{3}\right\rangle ^{\left(
3\right)  }=\left\langle a_{1},b_{1}\right\rangle ^{\left(  3\right)
}\left\langle a_{2},b_{2}\right\rangle ^{\left(  3\right)  }\left\langle
a_{3},b_{3}\right\rangle ^{\left(  3\right)  }$. These constructions can
naturally lead to ternary generalizations of the duality concept and the
quantum double which are key ingredients in the theory of quantum groups
\cite{dri0,kassel,majid}.

\section{Conclusions}

In this paper we presented a review of polyadic systems and their
representations, ternary algebras and Hopf algebras. We have classified
general polyadic systems and considered their homomorphisms and their
multiplace generalizations, paying attention to their associativity. Then, we
defined multiplace representations and multiactions and have given examples of
matrix representations for some ternary groups. We defined and investigated
ternary algebras and Hopf algebras, and have given some examples. We then
considered some ternary generalizations of quantum groups and the Yang-Baxter equation.

\end{document}